\documentclass[11pt]{article}
\usepackage{amsmath,mathrsfs,amsfonts,amssymb,stmaryrd,xy,dsfont,amsthm,pslatex} 
\usepackage[a4paper]{geometry}
\usepackage[pdftex]{hyperref}
\theoremstyle{definition} 
\newtheorem{remark}{Remark}[section]
\newtheorem{remarks}[remark]{Remarks}
\newtheorem{example}[remark]{Example}
\newtheorem{examples}[remark]{Examples} 
\newtheorem*{example*}{Example}
\theoremstyle{plain}
\newtheorem{definition}[remark]{Definition}
\newtheorem{theorem}[remark]{Theorem}
\newtheorem{proposition}[remark]{Proposition}
\newtheorem{corollary}[remark]{Corollary}
\newtheorem{lemma}[remark]{Lemma}
\newtheorem*{assumption*}{Assumption}

\newcommand{\mfrak}[1]{\mathfrak{#1}}
\newcommand{\frakN}{\mfrak{N}}

\newcommand{\complex}{\mathds{C}}
\newcommand{\reals}{\mathds{R}}
\newcommand{\naturals}{\mathds{N}}

\DeclareMathOperator{\Id}{id}
\DeclareMathOperator{\Tot}{Tot}
\newcommand{\TotL}{\Tot{\cal L}}

\DeclareMathOperator{\supp}{supp}
\DeclareMathOperator{\lspan}{span}

\DeclareMathOperator{\Ad}{Ad}

\DeclareMathOperator{\Ind}{Ind}

\DeclareMathOperator{\Mor}{Mor}

\newcommand{\rd}{\,\mathrm{d}}

\newcommand{\bfC}{\mathbf{C}}
\newcommand{\bfD}{\mathbf{D}}
\newcommand{\bfCat}{\mathbf{Cat}}
\newcommand{\bfCatn}{\mathbf{Cat}_{(N,N^{op})}}
\newcommand{\bfEnd}{\mathbf{End}}

\newcommand{\bfF}{\mathbf{F}}
\newcommand{\bfH}{\mathbf{H}}
\newcommand{\bfR}{\mathbf{R}}
\newcommand{\bfU}{\mathbf{U}}

\newcommand{\bfG}{\mathbf{G}}
\newcommand{\bfB}{\mathbf{B}}

\newcommand{\hilb}{\mathbf{Hilb}}
\newcommand{\ob}{\mathrm{ob}\,}

\newcommand{\bfcs}{\mathbf{C}^{\mathbf{*}}}
\newcommand{\bfcsnd}{\mathbf{C}^{\mathbf{*}{s,nd}}}
\newcommand{\bfws}{\mathbf{W^{*}}}
\newcommand{\bfwsn}{\mathbf{W^{*}}_{(N,N^{op})}}

\newcommand{\czalg}{\bfcs_{C_{0}(Z)}}

\newcommand{\cbalg}{\tilde{\mathbf{C}}^{\mathbf{*}}_{\frakb}}

\newcommand{\bfrep}{\mathbf{Rep}}

\newcommand{\bfBund}{\mathbf{Bdl}}

\newcommand{\bfMod}{\mathbf{Mod}}
\newcommand{\bfmod}{\mathbf{mod}}
\newcommand{\bfbimod}{\mathbf{bimod}}

\newcommand{\wmodr}[1]{\bfws\text{-}\bfmod_{#1}}
\newcommand{\wmodrn}{\bfws\text{-}\bfmod_{(N,N^{op})}}
\newcommand{\cbimodr}[1]{\bfcs\text{-}\bfmod_{(#1^{\dag},#1)}}
\newcommand{\wbimodr}[1]{\bfws\text{-}\bfbimod_{(N,N^{op})}}

\newcommand{\calgss}[1]{\mathbf{C}^{* (s)}_{#1}}
\newcommand{\calgs}[1]{\mathbf{C}^{* s}_{#1}}
\newcommand{\calgr}[1]{\mathbf{C}^{* r}_{#1}}
\newcommand{\calg}[1]{\bfcs_{#1}}

\newcommand{\frakA}{\mathfrak{A}}
\newcommand{\frakB}{\mathfrak{B}}
\newcommand{\frakC}{\mathfrak{C}}

\newcommand{\frakAo}{\mathfrak{A}^{\dag}}
\newcommand{\frakBo}{\mathfrak{B}^{\dag}}
\newcommand{\frakCo}{\mathfrak{C}^{\dag}}

\newcommand{\frakU}{\mathfrak{U}}

\newcommand{\fraka}{\mathfrak{a}}
\newcommand{\frakb}{\mathfrak{b}}
\newcommand{\frakc}{\mathfrak{c}}
\newcommand{\frakd}{\mathfrak{d}}

\newcommand{\frakao}{\mathfrak{a}^{\dag}}
\newcommand{\frakbo}{\mathfrak{b}^{\dag}}
\newcommand{\frakco}{\mathfrak{c}^{\dag}}

\newcommand{\frakH}{\mathfrak{H}}
\newcommand{\frakK}{\mathfrak{K}}
\newcommand{\frakL}{\mathfrak{L}}

\newcommand{\cbasel}[2]{(\mathfrak{#2},\mathfrak{#1},  \mathfrak{#1}^{\dag})}
\newcommand{\cbases}[2]{\cbasel{#1}{#2}}
\newcommand{\cbaseos}[2]{(\mathfrak{#2},\mathfrak{#1}^{\dag},  \mathfrak{#1})} 
\newcommand{\cbasesa}{\cbases{A}{H}}
\newcommand{\cbasesb}{\cbases{B}{K}}
\newcommand{\cbasesc}{\cbases{C}{L}}

\newcommand{\cbaseosb}{\cbaseos{B}{K}}

\newcommand{\trivbase}{\mathfrak{t}}

\newcommand{\aHb}{{_{\alpha}H_{\beta}}}

\newcommand{\cKd}{{_{\gamma}K_{\delta}}}

\newcommand{\eLf}{{_{\epsilon}L_{\phi}}}

\newcommand{\aHbn}[1]{{_{\alpha_{#1}}H^{#1}_{\beta_{#1}}}}
\newcommand{\cKdn}[1]{{_{\gamma_{#1}}K^{#1}_{\delta_{#1}}}}

\newcommand{\frakAoHA}{{_{\frakAo}\frakH_{\frakA}}}
\newcommand{\frakBoKB}{{_{\frakBo}\frakK_{\frakB}}}

\newcommand{\lt}{\triangleleft}
\newcommand{\rt}{\triangleright}

\newcommand{\hA}{\widehat{A}}

\newcommand{\frei}{\,\cdot\,}
\newcommand{\mycong}{\xrightarrow{\cong}}

\newcommand{\ftimes}[3]{ {_{#1}\! \underset{#2}{\times}\! {}_{#3}}}
\newcommand{\rtensor}[3]{ {_{#1}\! \underset{#2}{\otimes}\! {}_{#3}}}

\newcommand{\rtensorbc}{\rtensor{\beta}{\frakb}{\gamma}}

\newcommand{\btensor}{\underset{\frakb}{\otimes}}

\newcommand{\ctensor}{\underset{\frakc}{\otimes}}
\newcommand{\rtensorh}{\underset{\frakb}{\otimes}}

\newcommand{\rbtensor}[2]{\rtensor{#1}{\frakb}{#2}}

\newcommand{\fibre}[3]{ {_{#1}\! \underset{#2}{\ast}\!
    {}_{#3}}}
\newcommand{\sfibre}[1]{ {\underset{#1}{\ast}}}

\newcommand{\bfibre}{\underset{\frakb}{\ast}}
\newcommand{\cfibre}{\underset{\frakc}{\ast}}

\newcommand{\tl}{\ensuremath \olessthan}
\newcommand{\tr}{\ensuremath \ogreaterthan}

\newcommand{\HfibreK}{H \rtensorbc K}

\newcommand{\fibrebc}{\fibre{\beta}{\frakb}{\gamma}}

\newcommand{\AfibreB}{A \fibrebc B}

\newcommand{\Ab}{A_{\beta}}
\newcommand{\Cl}{C_{\lambda}}
\newcommand{\gB}{{_{\gamma}B}}
\newcommand{\mD}{{_{\mu}D}}

\newcommand{\ind}[2]{\mathrm{Ind}_{#1}({#2})}

\newcommand{\kalpha}[1]{|\alpha\rangle_{{#1}}}
\newcommand{\balpha}[1]{\langle\alpha|_{{#1}}}
\newcommand{\kbeta}[1]{|\beta{}\rangle_{{#1}}}
\newcommand{\bbeta}[1]{\langle\beta|_{{#1}}}

\newcommand{\kgamma}[1]{|\gamma\rangle_{{#1}}}
\newcommand{\bgamma}[1]{\langle\gamma|_{{#1}}}

\newcommand{\HtensorK}{H \rbtensor{\beta}{\gamma} K}

\title{The relative tensor product and a  minimal fiber
  product in the setting of $C^{*}$-algebras}

\author{Thomas Timmermann\footnote{This work was supported
    by the SFB 478 ``Geometrische Strukturen in der
    Mathematik'' funded by the Deutsche
    Forschungsgemeinschaft (DFG) and initiated
    during a stay at the ``Special Programme on Operator
    Algebras'' at the Fields Institute in Toronto, Canada,
    2007.
  } \\
  University of Muenster, \\ Einsteinstr.\ 62, 48149
  Muenster, Germany, \\ timmermt@math.uni-muenster.de \\
  Phone: ++49 251 8332724, Fax: ++49 251 8332708}

\date{\today}

\begin{document} 

\xyrequire{matrix} \xyrequire{arrow} \xyrequire{curve}
\CompileMatrices

\maketitle

\begin{abstract}
  We introduce a relative tensor product of $C^{*}$-bimodules
  and a spatial fiber product of $C^{*}$-algebras that are
  analogues of Connes' fusion of correspondences and the
  fiber product of von Neumann algebras introduced by
  Sauvageot, respectively.  These new constructions form the
  basis for our approach to quantum groupoids in the setting
  of $C^{*}$-algebras that is published separately.


\end{abstract}

\tableofcontents

\section{Introduction}
\label{section:introduction}

\paragraph{Background} The relative tensor product of
Hilbert modules over von Neumann algebras was 
  introduced  by Connes in an unpublished manuscript
\cite{connes:3,falcone,sauvageot} and later used by
Sauvageot to define a fiber product of von
Neumann algebras relative to a common (commutative) von
Neumann subalgebra \cite{sauvageot:2}. These
constructions and Haagerups theory of operator-valued
weights on von Neumann algebras \cite{haagerup:1,haagerup:2}
form the basis for the theory of measured quantum groupoids
developed by Enock, Lesieur and Vallin
\cite{enock:action,enock:1,lesieur,vallin:1,vallin:2}.

In this article, we introduce a new notion of a bimodule in
the setting of $C^{*}$-algebras, construct relative tensor
products of such bimodules, and define a fiber product of
$C^{*}$-algebras represented on such bimodules. These
constructions form the basis for a series of articles on
quantum groupoids in the setting of $C^{*}$-algebras,
individually addressing fundamental unitaries
\cite{timmermann:cpmu-hopf}, axiomatics of the compact case
\cite{timmermann:leiden}, and coactions of quantum groupoids
on $C^{*}$-algebras \cite{timmermann:coactions}.  Moreover,
our previous approach to quantum groupoids in the setting of
$C^{*}$-algebras \cite{timmermann:hopf} embeds functorially
into this new framework \cite{timmer:comparison}, and the
latter overcomes the serious restrictions of the former one.

Already in the definition of a quantum groupoid, the
relative tensor product and a fiber product appear as
follows.  Roughly, such an object consists of the following
ingredients: an algebra $B$, thought of as the functions on
the unit space, an algebra $A$, thought of as functions on
the total space, a homomorphism $r \colon B \to A$ and an
antihomomorphism $s \colon B \to A$ corresponding to the
range and the source map, and a comultiplication $\Delta
\colon B \to A \underset{B}{\ast} A$ corresponding to the
multiplication of the quantum groupoid. Here, $A
\underset{B}{\ast} A$ is a fiber product whose precise
definition depends on the class of the algebras involved. In
the setting of operator algebras, $A$ acts naturally on some
bimodule $H$ and product $A \underset{B}{\ast} A$ is a
certain subalgebra of operators acting on a relative tensor
product $H \underset{B}{\otimes} H$. This relative tensor
product is important also because it forms the domain or
range of the fundamental unitary of the quantum groupoid.

\paragraph{Overview}

Let us now sketch the problems and constructions studied in
this article.

The first problem is the construction of a tensor product $H
\underset{B}{\otimes} K$ of modules $H,K$ over some algebra
$B$.  In the algebraic setting, $H \underset{B}{\otimes} K$
is simply a quotient of the full tensor product $H \otimes
K$. In the setting of von Neumann algebras, $H$ and $K$ are
Hilbert spaces, and Connes explained that the right tensor
product is not a completion of the algebraic one but
something more complicated. If $B$ is commutative and of the
form $B=L^{\infty}(X,\mu)$, then the modules $H,K$ can be
disintegrated into two measurable fields of Hilbert spaces
in the form $H=\int_{X}^{\oplus} H_{x} d\mu(x)$ and
$K=\int_{X}^{\oplus} K_{x} d\mu(x)$, and $H
\underset{B}{\otimes} K$ is obtained by taking tensor
products of the fibers and integrating again: $H
\underset{B}{\otimes} K =\int^{\oplus}_{X} H_{x}\otimes
K_{x} d\mu(x)$.  For the situation where $B$ is a
$C^{*}$-algebra, we propose an approach that is based on the
internal tensor product of Hilbert $C^{*}$-modules and
essentially consists of an algebraic reformulation of
Connes' fusion. Central to this approach is a new notion of
a bimodule in the setting of $C^{*}$-algebras.

The second problem is the construction of a fiber product $A
\underset{B}{\ast} C$ of two algebras $A,C$ relative to a
subalgebra $B$. If $B$ is central in $A$ and the opposite
$B^{op}$ is central in $C$, this fiber product is just a
relative tensor product. In the algebraic setting, it
coincides with the tensor product of modules; in the setting
of operator algebras, it can be obtained via disintegration
and a fiberwise tensor product again. This approach was
studied by Sauvageot for Neumann algebras
\cite{sauvageot:2}, and by Blanchard \cite{blanchard} for
$C^{*}$-algebras.

The case where the subalgebra $B^{(op)}$ is no longer
central in $A$ or $C$ is more difficult. In the algebraic
setting, the fiber product was introduced by Takeuchi
\cite{takeuchi} and is, roughly, the largest subalgebra of
the relative tensor product $A \underset{B}{\otimes} C$
where componentwise multiplication is still well defined. In
the setting of von Neumann algebras, Sauvageot's definition
of the fiber product carries over to the general case and
takes the form $A \underset{B}{\ast} C = (A'
\underset{B}{\otimes} C')'$, where $A$ and $C$ are
represented on Hilbert spaces $H$ and $K$, respectively, and
$A' \underset{B}{\otimes} C'$ acts on Connes' relative
tensor product $H \underset{B}{\otimes} K$. Here, it is
important to note that $A' \underset{B}{\otimes} C'$ is a
completion of an algebraic tensor product spanned by
elementary tensors, but in general, $A \underset{B}{\ast} C$
is not. Similarly, in the setting of $C^{*}$-algebras, one
can not start from some algebraic tensor product and define
the fiber product to be some completion; rather, a new idea
is needed. We propose such a new fiber product for
$C^{*}$-algebras represented on the new class of modules
mentioned above.  Unfortunately, several important questions
concerning this construction remain open, but the
applications in
\cite{timmermann:leiden,timmermann:coactions,timmermann:cpmu-hopf}
already prove its usefulness.

\paragraph{Plan}

This article is organized as follows.  

The introduction ends with a short summary on terminology
and some background on Hilbert $C^{*}$-modules.  

Section \ref{section:rtp} is devoted to the relative tensor
product in the setting of $C^{*}$-algebras.  It starts with
some motivation, then presents a new notion of modules and
bimodules in the setting of $C^{*}$-algebras, and finally
gives the construction and its formal properties like
functoriality, associativity and unitality.

Section \ref{section:fiber} introduces a fiber product of
$C^{*}$-algebras. It starts with an overview and then
proceeds to $C^{*}$-algebras represented on the class of
modules and bimodules introduced in Section
\ref{section:rtp}.  The fiber product is first defined and
studied for such represented $C^{*}$-algebras, including a
discussion of functoriality, slice maps, lack of
associativity, and unitality.  A natural extension to
non-represented $C^{*}$-algebras is indicated at the end.

Section \ref{section:vn} relates our constructions for the
setting of $C^{*}$-algebras to the corresponding
constructions for the setting of von Neumann algebras.
Adapting our constructions to von Neumann algebras, one
recovers Connes fusion and Sauvageot's fiber product;
moreover, the constructions are related by functors going
from the $C^{*}$-level to the $W^{*}$-level. The section
ends with a categorical interpretation of  Sauvageot's
fiber product.

Section \ref{section:commutative} shows that for a
commutative base $B=C_{0}(X)$, the relative tensor product of
the new class of modules corresponds to the fiberwise tensor
product of continuous Hilbert bundles over $X$, and the
fiber product of represented $C^{*}$-algebras is related to
the relative tensor product of continuous
$C_{0}(X)$-algebras studied by Blanchard.

\paragraph{Preliminaries and notation}

Given a category $\bfC$, we write $A,B \in \bfC$ to indicate
that $A,B$ are objects of $\bfC$, and denote by $\bfC(A,B)$
the associated set of morphisms.

Given a subset $Y$ of a normed space $X$, we denote by $[Y]
\subset X$ the closed linear span of $Y$.

All sesquilinear maps like inner products on Hilbert spaces
are assumed to be conjugate-linear in the first component
and linear in the second one.  

Given a Hilbert space $H$ and an element $\xi \in H$, we
define ket-bra operators $|\xi\rangle \colon \complex \to
H$, $\lambda \mapsto \lambda \xi$, and
$\langle\xi|=|\xi\rangle^{*} \colon H \to \complex$, $\xi'
\mapsto \langle \xi|\xi'\rangle$.

 We shall make extensive use of (right) Hilbert
 $C^{*}$-modules; a standard reference is \cite{lance}.

 Let $A$ and $B$ be $C^{*}$-algebras.  Given Hilbert
 $C^{*}$-modules $E$ and $F$ over $B$, we denote by
 $\mathcal{L}(E,F)$ the space of all adjointable operators
 from $E$ to $F$.  Let $E$ and $F$ be Hilbert
 $C^{*}$-modules over $A$ and $B$, respectively, and let
 $\pi \colon A \to \mathcal{L}(F)$ be a
 $*$-homomorphism. Then the internal tensor product $E
 \otimes_{\pi} F$ is a Hilbert $C^{*}$-module over $B$
 \cite[\S 4]{lance} and the closed linear span of elements
 $\eta \otimes_{\pi} \xi$, where $\eta \in E$ and $\xi \in
 F$ are arbitrary, and $\langle \eta \otimes_{\pi} \xi|\eta'
 \otimes_{\pi} \xi'\rangle = \langle
 \xi|\pi(\langle\eta|\eta'\rangle)\xi'\rangle$ and $(\eta
 \otimes_{\pi} \xi)b=\eta \otimes_{\pi} \xi b$ for all
 $\eta,\eta' \in E$, $\xi,\xi' \in F$, $b \in B$.  We denote
 the internal tensor product by ``$\tr$'' and drop the index
 $\pi$ if the representation is understood; thus, $E \tr F=E
 \tr_{\pi} F=E \otimes_{\pi} F$.

 We define a {\em flipped internal tensor product} $F
 {_{\pi}\tl} E$ as follows. We equip the algebraic tensor
 product $F \odot E$ with the structure maps $\langle \xi
 \odot \eta | \xi' \odot \eta'\rangle := \langle \xi|
 \pi(\langle \eta|\eta'\rangle) \xi'\rangle$, $(\xi \odot
 \eta) b := \xi b \odot \eta$, form the separated
 completion, and obtain a Hilbert $C^{*}$-$B$-module $F
 {_{\pi}\tl} E$ which is the closed linear span of elements
 $\xi {_{\pi}\tl} \eta$, where $\eta \in E$ and $\xi \in F$
 are arbitrary, and $\langle \xi {_{\pi}\tl} \eta|\xi'
 {_{\pi}\tl} \eta'\rangle = \langle
 \xi|\pi(\langle\eta|\eta'\rangle)\xi'\rangle$ and $(\xi
 {_{\pi}\tl} \eta)b=\xi b {_{\pi}\tl} \eta$ for all
 $\eta,\eta' \in E$, $\xi,\xi' \in F$, $b\in B$. As above,
 we usually drop the index $\pi$ and simply write ``$\tl$''
 instead of ``${_{\pi}\tl}$''.  Evidently, there exists a
 unitary $\Sigma \colon F \tr E \mycong E \tl F$, $\eta \tr
 \xi \mapsto \xi \tl \eta$.


 Let $E_{1},E_{2}$ be Hilbert $C^{*}$-modules over $A$, let
 $F_{1}$, $F_{2}$ be Hilbert $C^{*}$-modules over $B$ with
 $*$-homomorphisms $\pi_{i} \colon A \to \mathcal{
   L}(F_{i})$ for $i=1,2$, and let $S \in \mathcal{
   L}(E_{1},E_{2})$, $T \in \mathcal{
   L}(F_{1},F_{2})$ such that $T\pi_{1}(a)=\pi_{2}(a)T$
 for all $a \in A$. Then there exists a unique operator $S
 \tr T \in \mathcal{ L}(E_{1} \tr F_{1}, E_{2} \tr
 F_{2})$ such that $(S \tr T)(\eta \tr \xi)= S\eta \tr T\xi$
 for all $\eta \in E_{1}$, $\xi \in F_{1}$, and $(S\tr
 T)^{*}=S^{*} \tr T^{*}$ \cite[Proposition 1.34]{echter}.

 \section{The relative tensor product in the setting of
$C^{*}$-algebras}

\label{section:rtp}

\subsection{Motivation}

\label{subsection:rtp-background}

The aim of this section is to construct a relative tensor
product of suitably defined left and right modules over a
general $C^{*}$-algebra $B$ such that i) the construction
shares the main properties of the ordinary tensor product of
bimodules over rings like functoriality and associativity
and ii) the modules admit representations of
$C^{*}$-algebras that do not commute with the module
structures.  The latter condition will be needed to
construct fiber products of $C^{*}$-algebras; see Section
\ref{section:fiber}.

The internal tensor product of Hilbert $C^{*}$-modules meets
condition i) but not ii) because $C^{*}$-algebras represented
on such modules necessarily commute with the right module
structure.  An approach to quantum groupoids based on the
internal tensor product was developed in
\cite{timmermann:hopf} but remained
restricted to very special cases.

What we are looking for is an analogue of Connes' fusion of
correspondences. Here, $B$ is a von Neumann algebra, and
left and right modules are Hilbert spaces equipped with
suitable representation or antirepresentation of $B$,
respectively.  The relative tensor product of a right module
$H$ and a left module $K$ is then constructed as
follows. Choose a normal, semi-finite, faithful (n.s.f.)
weight $\mu$ on $B$, construct a $B$-valued inner product
$\langle \frei|\frei\rangle_{\mu}$ on the dense subspace
$H_{0}\subseteq H$ of all bounded vectors, and define $H
\underset{\mu}{\otimes} K$ to be the separated completion of
the algebraic tensor product $H_{0} \odot K$ with respect to
the sesquilinear form given by $\langle \xi \odot
\eta|\xi'\odot \eta'\rangle=\langle
\eta|\langle\xi|\xi'\rangle_{\mu}\eta'\rangle$. The
definition of bounded vectors involves the GNS-space
$\frakH:=H_{\mu}$ for $\mu$ which --- by Tomita-Takesaki
theory --- is bimodule over $B$, and each
bounded vector $\xi \in H_{0}$ gives rise to a map $L(\xi)
\in {\cal L}(\frakH_{B},H_{B})$ of right $B$-modules such
that $\langle \xi|\xi'\rangle_{\mu} = L(\xi)^{*}L(\xi') \in
B \subseteq {\cal L}(\frakH)$.

\begin{example*} Assume that $B=L^{\infty}(X,\mu)$ for some
  nice measure space $(X,\mu)$, and denote the weight on $B$
  given by integration by $\mu$ as well. Then
  $\frakH=L^{2}(X,\mu)$, and we can disintegrate $H$ and $K$
  into measurable fields $(H_{x})_{x}$ and $(K_{x})_{x}$ of
  Hilbert spaces over $X$ such that $H\cong
  \int^{\oplus}_{X} H_{x} d\mu(x)$ and $K\cong
  \int^{\oplus}_{X} K_{x} d\mu(x)$. Each vector $\xi$ of $H$
  or $K$ corresponds to a measurable section $ x \mapsto
  \xi(x)$ with square-integrable norm function $|\xi| \colon
  x \mapsto \|\xi_{x}\|$, and is bounded with respect to
  $\mu$ if and only if this norm function is essentially
  bounded. Then for all $\xi,\xi' \in H_{0}$, $x \in X$,
  $\eta,\eta'\in K$,
  \begin{align*}
    \langle \xi|\xi'\rangle_{\mu}(x) &=
    \langle\xi(x)|\xi'(x)\rangle_{H_{x}}, & \langle \xi
    \odot \eta|\xi'\odot \eta'\rangle &= \int_{X} \langle
    \xi(x)|\xi'(x)\rangle\langle \eta(x)|\eta'(x)\rangle
    d\mu(x),
\end{align*}
and $H \underset{\mu}{\otimes} K \cong \int_{X}^{\oplus}
H_{x} \otimes K_{x} d\mu(x)$.  Note that the sesquilinear
form above need not extend to $H \odot K$ because the
integrand need not be in $L^{1}(X,\mu)$ for arbitrary
$\xi,\xi' \in H$ and $\eta,\eta'\in K$.
\end{example*} 
For our purpose, the following algebraic description of $H
\underset{\mu}{\otimes} K$ is useful. This relative tensor
product can be identified with the separated completion of
algebraic tensor product
\begin{align} \label{eq:rtp-connes-algebraic} {\cal
L}(\frakH_{B},H_{B}) \odot \frakH \odot {\cal
L}({_{B}\frakH},{_{B}K})
\end{align} with respect to the sesquilinear form $\langle S
\odot \zeta \odot T|S' \odot \zeta' \odot T'\rangle =
\langle \zeta| S^{*}S' T^{*}T'\zeta'\rangle = \langle
\zeta|T^{*}T' S^{*}S' \zeta'\rangle$, where ${\cal
  L}(\frakH_{B},H_{B})$ and ${\cal L}({_{B}\frakH},{_{B}K})$
are all bounded maps of right or left $B$-modules,
respectively.  We adapt this definition to the setting of
$C^{*}$-algebras, making the following modifications:
\begin{itemize}
\item[(A)] The construction above depends on the choice of
  some n.s.f.\ weight $\mu$ or, more precisely,  the
  triple $(H_{\mu},\pi_{\mu}(B),\pi_{\mu}(B)')$, but any
  other $\mu$ yields a triple which is unitarily
  equivalent. In the setting of $C^{*}$-algebras, such a
  canonical triple does not exist but has to be chosen.
\item[(B)] The module structure of $H$ and $K$ can
  equivalently be described in terms of
  (anti)repre\-sentations of $B$ or in terms of the spaces
  ${\cal L}(\frakH_{B},H_{B})$ and ${\cal
    L}({_{B}\frakH},{_{B}K})$. In the setting of
  $C^{*}$-algebras, this equivalence breaks down, and we
  shall make suitable closed subspaces of intertwiners the
  primary object. In the commutative case, a representation
  corresponds to a measurable field of Hilbert spaces, and
  the subspaces fix a continuous structure.
\item[(C)] If $H$ and $K$ are bimodules, then so is $H
  \underset{\mu}{\otimes} K$. Here, a bimodule structure on
  $H$ is given by the additional choice of a representation
  of some von Neumann algebra $A$ that commutes with the
  antirepresentation of $B$ or, equivalently, satisfies
  $A{\cal L}(\frakH_{B},H_{B})={\cal
    L}(\frakH_{B},H_{B})$. If we pass to $C^{*}$-algebras,
  then commutation is too weak, and we shall adopt the second
  condition, where ${\cal L}(\frakH_{B},H_{B})$ is replaced
  by the subspace of intertwiners mentioned above.
\end{itemize}

\subsection{Modules and bimodules over $C^{*}$-bases}

\label{subsection:rtp-modules}

Observation (A) leads us to adopt the following terminology.
\begin{definition} \label{definition:cbase} A {\em
$C^{*}$-base} $\frakb=\cbasesb$ consists of a Hilbert space
$\frakH$ and commuting nondegenerate $C^{*}$-algebras
$\frakB,\frakBo \subseteq {\cal L}(\frakK)$, respectively.
The {\em opposite} of $\frakb$ is the
$C^{*}$-base $\frakbo:=\cbaseosb$.  A $C^{*}$-base
$\cbasesa$ is {\em equivalent} to $\frakb$ if
$\Ad_{V}(\frakA) = \frakB$ and $\Ad_{V}(\frakAo) = \frakBo$
for some unitary $V \in {\cal L}(\frakH, \frakK)$.
\end{definition} 
Clearly, the Hilbert space $\complex$ and
twice the algebra $\complex\equiv{\cal L}(\complex)$ form a
trivial $C^{*}$-base $\trivbase
=(\complex,\complex,\complex)$.
\begin{example}\label{example:rtp-cbase} 
  Let $\mu$ be a proper, faithful KMS-weight on a
  $C^{*}$-algebra $A$ \cite{kustermans:kms} with GNS-space
  $H_{\mu}$, GNS-representation $\pi_{\mu} \colon A \to
  {\cal L}(H_{\mu})$, modular conjugation $J_{\mu} \colon
  H_{\mu} \to H_{\mu}$, and opposite GNS-representation
  $\pi_{\mu^{op}} \colon A^{op} \to {\cal L}(H_{\mu})$, $a
  \mapsto J_{\mu}\pi_{\mu}(a^{*})J_{\mu}$.  Then
  $(H_{\mu},\pi_{\mu}(A),\pi_{\mu^{op}}(A^{op}))$ is a
  $C^{*}$-base.  Its opposite is equivalent to the
  $C^{*}$-base associated to the opposite weight $\mu^{op}$
  on $A^{op}$. Indeed, $H_{\mu}$ can be considered as the
  GNS-space for $\mu^{op}$ via the opposite GNS-map
  $\Lambda_{\mu^{op}} \colon \frakN_{\mu^{op}} \to H_{\mu}$,
  $a^{op} \mapsto J_{\mu}\Lambda_{\mu}(a^{*})$, and then
  $J_{\mu^{op}}\pi_{\mu^{op}}(A^{op})J_{\mu^{op}}=\pi_{\mu}(A)$.
\end{example} Let $\frakb=\cbasesb$ be a $C^{*}$-base. We
define $C^{*}$-modules over $\frakb$ as indicated in comment
(B).
\begin{definition} \label{definition:rtp-modules} A {\em
$C^{*}$-$\frakb$-module} $H_{\alpha}=(H,\alpha)$ is a
Hilbert space $H$ with a closed subspace $\alpha \subseteq
{\cal L}(\frakK,H)$ satisfying $[\alpha \frakK]=H$, $[\alpha
\frakB]=\alpha$, $[\alpha^{*}\alpha]=\frakB$.  A {\em
semi-morphism} between $C^{*}$-$\frakb$-modules $H_{\alpha}$
and $K_{\beta}$ is an operator $T \in \mathcal{L}(H,K)$
satisfying $T\alpha \subseteq \beta$. If additionally
$T^{*}\beta \subseteq \alpha$, we call $T$ a {\em
morphism}. We denote the set of all (semi-)morphisms by
$\mathcal{ L}_{(s)}(H_{\alpha},K_{\beta})$.
\end{definition} 
Evidently, the class of all
$C^{*}$-$\fraka$-modules forms a category with respect to
all semi-morphisms, and a $C^{*}$-category in the
sense of \cite{roberts} with
respect to all morphisms.

\begin{lemma} \label{lemma:rtp-morphism}
  \begin{enumerate}
  \item Let $H,K$ be Hilbert spaces and $I \subseteq {\cal
      L}(H,K)$ such that $[IH]=K$. Then there exists a
    unique normal, unital $*$-homomorphism $\rho_{I} \colon
    (I^{*}I)' \to (II^{*})'$ such that $\rho_{I}(x)S=Sx$ for
    all $x \in (I^{*}I)'$, $S \in I$.
  \item Let $H,K,L$ be Hilbert spaces and $I \subseteq {\cal
      L}(H,K)$, $J \subseteq {\cal L}(K,L)$ such that
    $[IH]=K$, $[JK]=L$, and $J^{*}JI\subseteq I$. Then
    $\rho_{I}((I^{*}I)') \subseteq (J^{*}J)'$ and $\rho_{J}
    \circ \rho_{I}=\rho_{JI}$.
  \end{enumerate}
 \end{lemma}
 \begin{proof}
   i) Uniqueness is evident. Let $x \in (I^{*}I)'$
   and $S_{1},\ldots,S_{n}\in I$,
   $\xi_{1},\ldots,\xi_{n} \in H$. Since $x^{*}x$
   commutes with each $S_{i}^{*}S_{j}$, the matrix
   $(S_{i}^{*}S_{j}x^{*}x)_{i,j} \in M_{n}({\cal L}(H))$ is
   dominated by $\|x^{*}x\|(S_{i}^{*}S_{j})_{i,j}$, and
   \begin{align*} \|\sum_{i} S_{i}x\xi_{i}\|^{2} =
\sum_{i,j} \langle
\xi_{i}|S_{i}^{*}S_{j}x^{*}x\xi_{j}\rangle \leq \|x\|^{2}
\sum_{i,j} \langle
\xi_{i}|S_{i}^{*}S_{j}\xi_{j}\rangle = \|x\|^{2}
\|\sum_{i}S_{i}\xi_{i}\|^{2}.
   \end{align*} 
   Hence, there exists an operator $\rho_{I}(x) \in {\cal
     L}(K)$ as claimed. One easily verifies that the
   assignment $x \mapsto \rho_{I}(x)$ is a $*$-homomorphism.
   It is normal because $[IH]=K$ and for all $S,T\in I$,
   $\xi,\eta \in K$, the functional $x \mapsto \langle S
   \xi|\rho_{I}(x)T\eta\rangle = \langle \xi|x
   S^{*}T\eta\rangle$ is normal.

   ii) Let $x \in (I^{*}I)'$. Then $\rho_{I}(x) \in J^{*}J$
   because $S^{*}T\rho_{I}(x)R =
   S^{*}TRx=\rho_{I}(x)S^{*}TR$ for all $S,T \in J$, $R \in
   I$, and $\rho_{JI}(x)=\rho_{J}(\rho_{I}(x))$ because
   $\rho_{JI}(x)TR = TRx = \rho_{J}(\rho_{I}(x))TR$ for all
   $T\in J$, $R \in I$.   
 \end{proof}
 \begin{lemma} \label{lemma:rtp-modules} Let $H_{\alpha}$ be
   a $C^{*}$-$\frakb$-module.
    \begin{enumerate}
    \item $\alpha$ is a Hilbert $C^{*}$-$\frakB$-module with
      inner product $(\xi,\xi') \mapsto \xi^{*}\xi'$.
    \item There exist isomorphisms $\alpha \tr \frakK \to
      H$, $\xi \tr \zeta \mapsto \xi \zeta$, and $\frakK \tl
      \alpha \to H$, $\zeta \tl \xi \mapsto \xi\zeta$.
    \item There exists a unique normal, unital and faithful
      representation $\rho_{\alpha} \colon \frakB' \to
      \mathcal{L}(H)$ such that $\rho_{\alpha}(x)(\xi
      \zeta)= \xi x \zeta$ for all $x \in \frakB'$, $\xi \in
      \alpha$, $\zeta \in \frakK$.
\item Let $K_{\beta}$ be a $C^{*}$-$\frakb$-module and $T
  \in \mathcal{ L}_{s}(H_{\alpha},K_{\beta})$. Then
  $T\rho_{\alpha}(x) = \rho_{\beta}(x)T$ for all $x \in
  \frakB'$. If additionally $T \in \mathcal{
    L}(H_{\alpha},K_{\beta})$, then left multiplication by
  $T$ defines an operator in $\mathcal{
    L}_{\frakB}(\alpha,\beta)$, again denoted by $T$.
\end{enumerate}
  \end{lemma}
  \begin{proof} 
    Assertions i) and ii) are obvious, and iii) follows from
    the preceding lemma.  To prove iv), let $x\in
    \frakB',\xi \in \alpha, \zeta \in \frakK$. Then $T\xi
    \in \beta$ and $T\rho_{\alpha}(x)\xi \zeta = T\xi x
    \zeta = \rho_{\beta}(x)T\xi \zeta$.
  \end{proof}
  \begin{example} \label{example:rtp-module-commutative} Let
    $Z$ be a locally compact Hausdorff space, $\mu$ a Radon
    measure on $Z$ of full support, and ${\cal
      H}=(H_{z})_{z}$ a continuous bundle of Hilbert spaces
    on $Z$ with full support. Then the Hilbert space
    $\frakK=L^{2}(Z,\mu)$ together with the $C^{*}$-algebras
    $\frakB=\frakBo=C_{0}(Z) \subseteq {\cal L}(\frakK)$
    forms a $C^{*}$-base. Let $H=\int^{\oplus}_{Z} H_{z}
    d\mu(z)$ and $\alpha=m(\Gamma_{0}({\cal H}))$, where for
    each section $\xi \in \Gamma_{0}({\cal H})$, the
    operator $m(\xi) \in {\cal L}(\frakK,H)$ is given by
    pointwise multiplication, $m(\xi)f = (\xi(z)f(z))_{z \in
      Z}$. Then $H_{\alpha}$ is a $C^{*}$-$\frakb$-module
    and $\rho_{\alpha} \colon \frakB' = L^{\infty}(Z,\mu)
    \to {\cal L}(H)$ is given by pointwise multiplication of
    sections by functions. Every $C^{*}$-$\frakb$-module
    arises in this way from a continuous bundle; see Section
    \ref{section:commutative}.
  \end{example}
  Let also $\fraka=\cbasesa$ be a $C^{*}$-base.  We define
  $C^{*}$-$(\frakao,\frakb)$-bimodules as indicated in (C).
  \begin{definition} \label{definition:rtp-bimodules} A {\em
      $C^{*}$-$(\frakao,\frakb)$-module} is a triple $\aHb=
    (H,\alpha,\beta)$, where $H$ is a Hilbert space,
    $(H,\alpha)$  a $C^{*}$-$\frakao$-module, $(H,\beta)$
    a $C^{*}$-$\frakb$-module, and
    $[\rho_{\alpha}(\frakA)\beta]=\beta$,
    $[\rho_{\beta}(\frakBo)\alpha]=\alpha$. The set of {\em
      (semi-)morphisms} between
    $C^{*}$-$(\frakao,\frakb)$-modules $\aHb$ and $\cKd$ is
    $\mathcal{ L}_{(s)}(\aHb,\cKd) := {\cal
      L}_{(s)}(H_{\alpha},K_{\gamma}) \cap {\cal
      L}_{(s)}(H_{\beta},K_{\delta})$.
 \end{definition}
 \begin{remark} By Lemma \ref{lemma:rtp-modules},
$[\rho_{\alpha}(\frakA),\rho_{\beta}(\frakBo)]=0$ for every
$C^{*}$-$(\frakao,\frakb)$-module $\aHb$.
\end{remark} Again, the class of all
$C^{*}$-$(\frakao,\frakb)$-modules forms a category with
respect to all semi-morphisms, and a
$C^{*}$-category with respect to all morphisms.
\begin{examples} \label{examples:rtp-bimodules}
  \begin{enumerate}
  \item $\frakH_{\frakA}$ is a $C^{*}$-$\fraka$-module,
    $\rho_{\frakA}(x)=x$ for all $x\in \frakA'$, and
    $\frakAoHA$ is a $C^{*}$-$(\frakao,\fraka)$-module
    because
    $[\rho_{\frakAo}(\frakA)\frakA]=[\frakA\frakA]=\frakA$
    and $[\rho_{\frakA}(\frakAo)\frakAo]=\frakAo$.
  \item Let $H_{\beta}$ be a $C^{*}$-$\frakb$-module, let
    $\trivbase=(\complex,\complex,\complex)$ be the trivial
    $C^{*}$-base, and let $\alpha={\cal
      L}(\complex,H)$. Then $\aHb$ is a
    $C^{*}$-$(\trivbase,\frakb)$-module.
  \item Let $({\cal H}_{i})_{i}$ be a family of
    $C^{*}$-$(\frakao,\frakb)$-modules, where ${\cal H}_{i}=
    (H_{i},\alpha_{i},\beta_{i})$ for each $i$.  Denote by
    $\boxplus_{i} \alpha_{i} \subseteq {\cal L}\big(\frakH,
    \oplus_{i} H_{i}\big)$ the norm-closed linear span of
    all operators of the form $\zeta \mapsto
    (\xi_{i}\zeta)_{i}$, where $(\xi_{i})_{i}$ is in the
    algebraic direct sum $ \bigoplus^{\mathrm{alg}}_{i}
    \alpha_{i}$, and similarly define $\boxplus_{i}
    \beta_{i} \subseteq {\cal L}\big(\frakK, \oplus_{i}
    H_{i}\big)$. Then the triple $\boxplus_{i} {\cal H}_{i}
    := \big(\oplus_{i} H_{i}, \boxplus_{i} \alpha_{i},
    \boxplus_{i} \beta_{i}\big)$ is a
    $C^{*}$-$(\frakao,\frakb)$-module, for each $j$, the
    canonical inclusions $\iota_{j} \colon H_{j} \to
    \oplus_{i} H_{i}$ and projection $\pi_{j} \colon
    \oplus_{i} H_{i} \to H_{j}$ are morphisms ${\cal H}_{j}
    \to \boxplus_{i} {\cal H}_{i}$ and $\boxplus_{i} {\cal
      H}_{i} \to {\cal H}_{j}$, and with respect to
    these maps, $\boxplus_{i} {\cal H}_{i}$ is the direct
    sum of the family $({\cal H}_{i})_{i}$.
 \end{enumerate}
 \end{examples}
 The following example shows how bimodules arise from
 conditional expectations.
 \begin{example} \label{example:rtp-bimodule} Let $B$ be a
   $C^{*}$-algebra with a KMS-state $\mu$ and associated
   $C^{*}$-base $\frakb$ (Example \ref{example:rtp-cbase}),
   let $A$ be a unital $C^{*}$-algebra containing $B$ such
   that $1_{A} \in B$, and let $\phi \colon A \to B$ be a
   faithful conditional expectation such that $\nu:=\mu
   \circ \phi$ is a KMS-state and $\phi \circ
   \sigma^{\nu}_{t} = \sigma^{\mu}_{t} \circ \phi$ for all
   $t \in \reals$.  Fix a GNS-construction $\pi_{\nu} \colon
   A \to {\cal L}(H_{\nu})$ for $\nu$ with modular
   conjugation $J_{\nu}\colon H_{\nu} \to H_{\nu}$, and
   define $\pi_{\nu}^{op} \colon A^{op} \to {\cal
     L}(H_{\nu})$ by $a\mapsto
   J_{\nu}\pi_{\nu}(a^{*})J_{\nu}$. Then the inclusion $B
   \hookrightarrow A$ extends to an isometry $\zeta \colon
   \frakK=H_{\mu} \hookrightarrow H_{\nu}=H$, and we obtain
   a $C^{*}$-$(\frakbo,\frakb)$-module $\aHb$, where
   $H=H_{\nu}$, $\alpha=[J_{\nu}\pi_{\nu}(A)\zeta]$, $\beta=
   [\pi_{\nu}(A)\zeta]$, and $\rho_{\alpha} \circ
   \pi_{\mu^{op}} = \pi_{\nu}^{op}$, $\rho_{\beta} \circ
   \pi_{\mu} = \pi_{\nu}$. Moreover, $\pi_{\nu}(A)+
   \pi_{\nu}^{op}((A \cap B')^{op}) \subseteq {\cal
     L}(H_{\alpha})$, $\pi_{\nu^{op}}(A^{op}) + \pi_{\nu}(A
   \cap B') \subseteq {\cal L}(H_{\beta})$.  For details,
   see \cite[\S 2--3]{timmermann:leiden}.
 \end{example}

\subsection{The relative tensor product}

\label{subsection:rtp-rtp}

The concepts introduced above allow us to adapt the
algebraic formulation of Connes' fusion to the setting of
$C^{*}$-algebras as follows.  Let $\frakb=\cbasesb$ be a
$C^{*}$-base, $H_{\beta}$ a $C^{*}$-$\frakb$-module, and
$K_{\gamma}$ a $C^{*}$-$\frakbo$-module. Then the {\em
  relative tensor product} of $H_{\beta}$ and $K_{\gamma}$
is the Hilbert space
\begin{align*} 
  \HfibreK := \beta \tr \frakK \tl \gamma,
\end{align*} 
which is spanned by elements $\xi \tr \zeta \tl \eta$, where
$\xi \in \beta$, $\zeta \in \frakK$, $\eta \in \gamma$, the
inner product being given by $\langle \xi \tr \zeta \tl
\eta|\xi' \tr \zeta' \tl \eta'\rangle = \langle \zeta |
\xi^{*}\xi' \eta^{*}\eta' \zeta'\rangle = \langle
\zeta|\eta^{*}\eta' \xi^{*}\xi' \zeta'\rangle$ for all
$\xi,\xi' \in \beta$, $\zeta,\zeta' \in \frakK$, $\eta,\eta'
\in \gamma$.
\begin{examples} \label{examples:rtp-rtp} 
  \begin{enumerate}\item 
    If $\frakb$ is the trivial $C^{*}$-base
    $\trivbase=(\complex,\complex,\complex)$, then
    $\beta={\cal L}(\complex,H)$, $\gamma={\cal
      L}(\complex,K)$, and $\HtensorK \cong H \otimes K$ via
    $\xi \tr \zeta \tl \eta \mapsto \xi\zeta \otimes \eta 1
    = \xi 1 \otimes \eta \zeta$.
  \item Let $Z$ be a locally compact Hausdorff space, $\mu$
    a Radon measure on $Z$ of full support, ${\cal
      H}=(H_{z})_{z}$ and ${\cal K}=(K_{z})_{z}$ continuous
    bundles of Hilbert spaces on $Z$ with full support, and
    $H_{\alpha},K_{\beta}$ the associated
    $C^{*}$-$\frakb$-modules as defined in Example
    \ref{example:rtp-module-commutative}. One easily checks
    that then we have an isomorphism
    \begin{align*}
      \HtensorK \to \int_{Z}^{\oplus} H_{z} \otimes K_{z} \,
      d\mu(z),  \quad
      m(\xi)\tr \zeta \tl m(\eta) \mapsto (\xi(z) \zeta(z)
      \otimes \eta(z))_{z \in Z}.
    \end{align*}
  \end{enumerate}
\end{examples}
Let us list some easy observations and a few definitions.
\begin{enumerate}
\item The isomorphisms in Lemma \ref{lemma:rtp-modules} ii),
applied to $H_{\beta}$ and $K_{\gamma}$, respectively, yield
the following identifications which we shall use without
further notice:
  \begin{gather*} \beta \tr_{\rho_{\gamma}} K \cong
\HtensorK \cong H {_{\rho_{\beta}} \tl} \gamma, \quad \xi
\tr \eta\zeta \equiv \xi \tr \zeta \tl \eta \equiv \xi \zeta
\tl \eta.
  \end{gather*}
\item For each $\xi \in \beta$ and $\eta \in \gamma$, there
exist bounded linear operators
  \begin{align*} |\xi\rangle_{{1}} \colon K &\to \beta
\tr_{\rho_{\gamma}} K = \HtensorK, \ \omega \mapsto \xi \tr
\omega, & |\eta\rangle_{{2}} \colon H &\to H
{_{\rho_{\beta}}\tl} \gamma= \HtensorK, \ \omega \mapsto
\omega \tl \eta,
  \end{align*} whose adjoints $ \langle
\xi|_{1}:=|\xi\rangle_{1}^{*}$ and $\langle\eta|_{{2}} :=
|\eta\rangle_{{2}}^{*}$ are given by
  \begin{align*} \langle \xi|_{1}\colon \xi' \tr \omega
&\mapsto \rho_{\gamma}(\xi^{*}\xi')\omega, &
\langle\eta|_{{2}} \colon \omega \tl\eta' &\mapsto
\rho_{\beta}( \eta^{*}\eta')\omega.
  \end{align*} We put $\kbeta{1} := \{ |\xi\rangle_{{1}}
\,|\, \xi \in \beta\} \subseteq {\cal L}(K,\HtensorK)$ and
similarly define $\bbeta{1}$, $\kgamma{2}$, $\bgamma{2}$.
\item For all $S \in \rho_{\beta}(\frakBo)'$ and $T \in
\rho_{\gamma}(\frakB)'$, we have operators
  \begin{align*} S \tl \Id &\in \mathcal{L}(H
{_{\rho_{\beta}}\tl} \gamma) = \mathcal{L}(\HtensorK), & \Id
\tr T &\in \mathcal{L}(\beta \tr_{\rho_{\gamma}} K) =
\mathcal{L}(\HtensorK).
  \end{align*} If these operators commute, we let $S
\btensor T:=(S\tl \Id)(\Id \tr T) = (\Id \tr T)(S \tl \Id)$.
The commutativity condition holds in each of the following
cases:
  \begin{enumerate}
  \item $S \in {\cal L}_{s}(H_{\beta})$; then $(S \btensor
T)(\xi \tr \omega)=S\xi \tr T\omega$ for each $\xi \in
\beta,\omega \in K$;
  \item $T \in {\cal L}_{s}(K_{\gamma})$; then $(S \btensor
T)(\omega \tl \eta) = S\omega \tl T\eta$ for each $\omega
\in H,\eta \in \gamma$;
\item $(\frakBo)'=\frakB''$; then for all $\xi,\xi' \in
  \beta$ and $\eta,\eta' \in \gamma$, the elements
  $\eta^{*}T\eta' \in \frakB'$ and $\xi^{*}S\xi' \in
  (\frakBo)'$ commute, and if $\zeta,\zeta'\in \frakK$ and
  $\omega=\xi \tr \zeta \tl \eta$, $\omega'=\xi' \tr \zeta'
  \tl \eta'$, then $\langle \omega|(\Id \tr T)(S \tl
  \Id)\omega'\rangle = \langle \zeta|
  (\eta^{*}T\eta')(\xi^{*} S \xi')\zeta'\rangle = \langle
  \zeta| (\xi^{*} S \xi') (\eta^{*}T\eta') \zeta'\rangle =
  \langle \omega|(S \tl \Id)(\Id \tr T)\omega'\rangle$.
\end{enumerate}
\end{enumerate} Let $\fraka=\cbasesa$ and $\frakc=\cbasesc$
be further $C^{*}$-bases.  Then the relative tensor product
of bimodules over $(\frakao,\frakb)$ and $(\frakbo,\frakc)$
is a bimodule over $(\frakao,\frakc)$:
\begin{proposition} \label{proposition:rtp-mod} Let ${\cal
    H} = \aHb$ be a $C^{*}$-$(\frakao,\frakb)$-module,
  ${\cal K}=\cKd$ a $C^{*}$-$(\frakbo,\frakc)$-module, and
  \begin{align} \label{eq:rtp-mod} \alpha \lt \gamma &:=
    [\kgamma{2} \alpha] \subseteq {\cal L}(\frakH,
    \HfibreK), & \beta \rt \delta &:= [\kbeta{1}\delta]
    \subseteq {\cal L}(\frakL,\HfibreK).
    \end{align} Then ${\cal H} \btensor {\cal K}:={_{(\alpha
\lt \gamma)}(\HfibreK)_{(\beta \rt \delta)}}$ is a
$C^{*}$-$(\frakao,\frakc)$-module and
\begin{align} \label{eq:rtp-rtp-rho} \rho_{(\alpha \lt
    \gamma)}(x) &= \rho_{\alpha}(x) \tl \Id \text{ for all }
  x \in (\frakAo)', & \rho_{(\beta \rt \delta)}(y) &= \Id
  \tr \rho_{\delta}(y) \text{ for all } y \in \frakC'.
    \end{align}
\end{proposition}
\begin{proof} 
  $(\HfibreK)_{(\alpha \lt \gamma)}$ is a
  $C^{*}$-$\frakao$-module because
  $[\alpha^{*}\bgamma{2}\kgamma{2}\alpha] = [\alpha^{*}
  \rho_{\beta}(\frakBo)\alpha] = \frakAo$,
  $[|\gamma\rangle_{2}\alpha\frakAo]
  =[|\gamma\rangle_{2}\alpha]$, and $ [|\gamma\rangle_{2}
  \alpha\frakH] = [|\gamma\rangle_{2} H] =
  \HfibreK$. Likewise, $(\HfibreK)_{(\beta \rt \delta)}$ is
  a $C^{*}$-$\frakc$-module.

  For all $x \in (\frakAo)'$, $\zeta \in \frakH$, $\theta
\in \alpha$, $\eta \in \gamma$, we have
$|\eta\rangle_{2}\theta \in \alpha \lt \gamma$ and hence
    \begin{align*} \rho_{(\alpha \lt \gamma)}(x)(\theta\zeta
\tl \eta) &= \rho_{(\alpha \lt
\gamma)}(x)|\eta\rangle_{2}\theta \zeta =
|\eta\rangle_{2}\theta x \zeta = \rho_{\alpha}(x)
\theta\zeta \tl \eta = (\rho_{\alpha}(x) \tl
\Id)(\theta\zeta \tl \eta).
    \end{align*} The first equation in
\eqref{eq:rtp-rtp-rho} follows, and a similar agument proves
the second one.

Finally, ${_{(\alpha \lt \gamma)}(\HfibreK)_{(\beta \rt
\delta)}}$ is a $C^{*}$-$(\frakao,\frakc)$-module because
$[\rho_{(\alpha \lt \gamma)}(\frakA)\kbeta{1} \delta]
=[|\rho_{\alpha}(\frakA)\beta\rangle_{1}\delta] =
[\kbeta{1}\delta]$ and $[\rho_{(\beta \rt
\delta)}(\frakCo)\kgamma{2}\alpha] = [\kgamma{2}\alpha]$.
\end{proof} In the situation above, we call ${\cal H}
\btensor {\cal K}$ the {\em relative tensor product} of
${\cal H}$ and ${\cal K}$. Note the following commutative
diagram of Hilbert spaces and closed spaces of operators
between them:
\begin{align*} \xymatrix@C=25pt@R=0pt{ {\frakH}
\ar[rd]^{\alpha} \ar@/_1pc/[rrdd]_(0.4){\alpha \lt \gamma} &
& {\frakK} \ar[ld]_{\beta} \ar[rd]^{\gamma} & & {\frakL}
\ar[ld]_{\delta} \ar@/^1pc/[lldd]^(0.4){\beta \rt \delta} \\
& H \ar[rd]^{\kgamma{2}} & & K \ar[ld]_{\kbeta{1}} & \\ & &
\HfibreK & & }
\end{align*}

Given a $C^{*}$-$\frakb$-module ${\cal H} = H_{\beta}$ and a
$C^{*}$-$(\frakbo,\frakc)$-module ${\cal K}=\cKd$, we
abbreviate $\HtensorK_{\delta}:=(\HtensorK)_{\beta \rt
\delta}$. Likewise, we write ${_{\alpha}\HtensorK}$ for
$(\HtensorK)_{\alpha \lt \gamma}$ and
${_{\alpha}\HtensorK_{\delta}}$ for ${_{\alpha \lt
\gamma}}(\HtensorK)_{\beta \rt \delta}$.

The relative tensor product is functorial, associative,
unital, and compatible with direct sums in the following
sense:
\begin{proposition} \label{proposition:rtp-properties} Let
${\cal H}=\aHb$, ${\cal H}^{1}=\aHbn{1},{\cal
H}^{2}=\aHbn{2}$ be $C^{*}$-$(\frakao,\frakb)$-modules,
${\cal K}=\cKd$, ${\cal K}^{1}=\cKdn{1}$, ${\cal
K}^{2}=\cKdn{2}$ $C^{*}$-$(\frakbo,\frakc)$-modules, and
${\cal L}=\eLf$ a $C^{*}$-$(\frakco,\frakd)$-module.
    \begin{enumerate}
    \item $S \rtensorh T \in {\cal L}\big({\cal H}^{1}
\btensor {\cal K}^{1}, {\cal H}^{2} \btensor {\cal
K}^{2}\big)$ for all $S \in {\cal L}({\cal H}^{1},{\cal
H}^{2})$, $T \in {\cal L}({\cal K}^{1},{\cal K}^{2})$.
\item The composition of the isomorphisms $(H_{\beta}
  \btensor {_{\gamma}K}_{\delta}) \ctensor {_{\epsilon}L}
  \cong (H_{\beta} \btensor {_{\gamma}K}) {_{\rho_{(\beta
        \rt \delta)}} \tl} \epsilon \cong \beta
  \tr_{\rho_{\gamma}} K {_{\rho_{\delta}}\tl} \epsilon$ and
  $\beta \tr_{\rho_{\gamma}} K {_{\rho_{\delta}}\tl}
  \epsilon \cong \beta \tr_{\rho_{(\gamma \lt \epsilon)}}
  (K_{\delta} \ctensor {_{\epsilon}L}) \cong H_{\beta}
  \btensor ({_{\gamma}K_{\delta}} \ctensor {_{\epsilon}L})$
is an isomorphism of $C^{*}$-$(\frakao,\frakc)$-modules $a_{\fraka,\frakb,\frakc,\frakd}({\cal
  L},{\cal K},{\cal H}) \colon ({\cal H} \btensor {\cal K})
\ctensor {\cal L} \to {\cal H} \btensor ({\cal K} \ctensor
{\cal L})$.
  \item Put ${\cal U}:=\frakBoKB$. Then there exist
isomorphisms
    \begin{align*} r_{\fraka,\frakb}({\cal H}) &\colon {\cal
H} \btensor {\cal U} \to {\cal H}, \ \xi \tr \zeta \tl
b^{\dag} \mapsto \xi b^{\dag} \zeta =
\rho_{\beta}(b^{\dag})\xi\zeta, \\ l_{\frakb,\frakc}({\cal
K}) &\colon {\cal U} \btensor {\cal K} \to {\cal K}, \ b \tr
\zeta \tl \eta \mapsto \eta b\zeta =
\rho_{\gamma}(b)\eta\zeta.
    \end{align*}
  \item Let $({\cal H}^{i})_{i}$ be a family of
$C^{*}$-$(\frakao,\frakb)$-modules and $({\cal K}^{j})_{j}$
a family of $C^{*}$-$(\frakbo,\frakc)$-modules. For each
$i,j$, denote by $\iota^{i}_{{\cal H}} \colon {\cal H}^{i}
\to \boxplus_{i'} {\cal H}^{i'}$, $\iota^{j}_{{\cal K}}
\colon {\cal K}^{j} \to \boxplus_{j'} {\cal K}^{j'}$ and
$\pi^{i}_{{\cal H}} \colon \boxplus_{i'} {\cal H}^{i'} \to
{\cal H}^{i}$, $\pi^{j}_{{\cal K}} \colon \boxplus_{j'}
{\cal K}^{j'} \to {\cal K}^{j}$ the canonical inclusions and
projections, respectively.  Then there exist inverse
isomorphisms $\boxplus_{i,j} ( {\cal H}^{i} \btensor {\cal
K}^{j}) \leftrightarrows (\boxplus_{i} {\cal H}^{i})
\btensor (\boxplus_{j} {\cal K}^{j})$, given by
$(\omega_{i,j})_{i,j} \mapsto \sum_{i,j} (\iota^{i}_{{\cal
H}} \btensor \iota^{j}_{{\cal K}})(\omega_{i,j})$ and
$\big((\pi^{i}_{{\cal H}} \btensor \pi^{j}_{{\cal
K}})(\omega)\big)_{i,j} \mapsfrom \omega$, respectively.
  \end{enumerate}
\end{proposition}
  \begin{proof} i) If $S,T$ are as above and ${\cal
H}^{i}=\aHbn{i}$, ${\cal K}^{j}=\cKdn{j}$ for $i,j=1,2$,
then $(S \rtensorh T)|\gamma_{1}\rangle_{2}\alpha_{1} =
|T\gamma_{1}\rangle_{2} S\alpha_{1} \subseteq
|\gamma_{2}\rangle_{2}\alpha_{2}$ and similarly $(S
\rtensorh T)|\beta_{1}\rangle_{1}\delta_{1} \subseteq
|\beta_{2}\rangle_{1} \delta_{2}$, $(S \rtensorh
T)^{*}|\gamma_{2}\rangle_{2}\alpha_{2} \subseteq
|\gamma_{1}\rangle_{2}\alpha_{1}$, $(S \rtensorh
T)^{*}|\beta_{2}\rangle_{1} \delta_{2} \subseteq
|\beta_{1}\rangle_{1} \delta_{1}$.
    
    ii) Straightforward.
    
    iii) $r_{\fraka,\frakb}({\cal H}) \cdot (\alpha \lt
\frakBo) = [\rho_{\beta}(\frakBo)\alpha]=\alpha$ and
$r_{\fraka,\frakb}({\cal H}) \cdot (\beta \rt \frakB ) =
[\beta \frakB]=\beta$. For $l_{\frakb,\frakc}({\cal K})$,
the arguments are similar.
    
    iv) Straightforward.
\end{proof}
\begin{remark} \label{remark:rtp-bicategory} The relative
  tensor product of modules and morphisms can be considered
  as the composition in a bicategory as follows. Recall that
  a bicategory $\bfB$ consists of a class of objects $\ob
  \bfB$, a category $\bfB(A,B)$ for each $A,B \in \ob \bfB$
  whose objects and morphisms are called {\em 1-cells} and
  {\em 2-cells}, respectively, a functor $c_{A,B,C} \colon
  \bfB(B,C) \times \bfB(A,B) \to \bfB(A,C)$
  (``composition'') for each $A,B,C \in \ob \bfB$, an object
  $1_{A} \in \bfB(A,A)$ (``identity'') for each $A \in \ob
  \bfB$, an isomorphism $a_{A,B,C,D}(f,g,h) \colon
  c_{A,B,D}(c_{B,C,D}(h,g),f) \to
  c_{A,C,D}(h,c_{A,B,C}(g,f))$ in $\bfB(A,D)$
  (``associativity'') for each triple of 1-cells $A
  \xrightarrow{f} B \xrightarrow{g} C \xrightarrow{h} D$ in
  $\bfB$, and isomorphisms $l_{A}(f) \colon
  c_{A,A,B}(f,1_{A}) \to f$ and $r_{B}(f) \colon
  c_{A,B,B}(1_{B},f) \to f$ in $\bfB(A,B)$ for each 1-cell
  $A \xrightarrow{f} B$ in $\bfB$, subject to several axioms
  \cite{leinster}.  Tedious but straightforward calculations
  show that there exists a bicategory $\bfcs$-$\bfbimod$ such that
\begin{enumerate}
\item the objects are all $C^{*}$-bases and
  $\bfcs$-$\bfbimod(\fraka,\frakb)$ is the category of all
  $C^{*}$-$(\frakao,\frakb)$-modules with morphisms (not
  semi-morphisms) for all $C^{*}$-bases $\fraka,\frakb$;
\item the functor $c_{\fraka,\frakb,\frakc}$ is given by
$(\cKd,\aHb) \mapsto \aHb \btensor \cKd$ and $(T,S) \mapsto
S \btensor T$, respectively, and the identity $1_{\fraka}$
is ${_{\frakAo}\frakH_{\frakA}}$ for all $C^{*}$-bases
$\fraka$, $\frakb$, $\frakc$, $\frakd$;
  \item $a, r, l$ are as in Proposition
\ref{proposition:rtp-properties}.
  \end{enumerate}
\end{remark}

\section{The spatial fiber product of $C^{*}$-algebras}

\label{section:fiber}

\subsection{Background}

\label{subsection:fiber-motivation}

We now use the relative tensor product to construct a fiber
product of $C^{*}$-algebras that are represented on
$C^{*}$-modules over $C^{*}$-bases.  To motivate our
approach, let us first review several related constructions.
In each case, the task is to construct a relative
tensor product or ``fiber product'' of two algebras $A$ and
$C$ with respect to a common subalgebra $B$.

First, assume that we are working in the category of unital
commutative rings.  Then the fiber product is just the
push-out of the diagram formed by $A,B,C$.  Explicitly, it
is the algebraic tensor product $A \underset{B}{\odot} C$,
where $A$ and $C$ are considered as modules over $B$, and
the multiplication is defined componentwise.  In the
category of commutative $C^{*}$-algebras, the push-out is
the maximal completion of the algebraic tensor product $A
\underset{B}{\odot} C$ and, as usual in the setting of
$C^{*}$-algebras, also other interesting completions exist
\cite{blanchard}.  For example, if $B=C_{0}(X)$ for some
locally compact Hausdorff space and if $A$ and $C$ are
represented on Hilbert spaces $H$ and $K$, respectively,
then $H$ and $K$ can be disintegrated over $X$ with respect
to some measure $\mu$ (see Subsection
\ref{subsection:rtp-background}), and the algebra $A
\underset{B}{\odot} C$ has a natural representation $\pi$ on
the relative tensor product $H \underset{\mu}{\otimes} K =
\int_{X}^{\oplus} H_{x} \otimes K_{x} d\mu(x)$, leading to a
minimal completion $\overline{\pi(A \underset{B}{\odot}
  C)}$.  In the setting of von Neumann algebras, $H$ and $K$
are intrinsic, and the desired fiber product is $\pi(A
\underset{B}{\odot} C)''\subseteq {\cal L}(H
\underset{\mu}{\otimes} K)$.  Note that all of these constructions
 do not depend on commutativity of $A$ and $C$
and make sense as long as $B$ is central in
$A$ and in $C$.

Next, consider the case where $A,B,C$ are non-commutative,
$B$ is a subalgebra of $A$, and the opposite $B^{op}$ is a
subalgebra of $C$.  Then one can consider $A$ and $C$ as
modules over $B$ via right multiplication, and form the
algebraic tensor product $A \underset{B}{\odot} C$, but
componentwise multiplication is well defined only on the
subspace $A \underset{B}{\times} C \subseteq A
\underset{B}{\odot} C$ which consists of all elements
$\sum_{i} a_{i} \odot c_{i}$ satisfying $\sum_{i} ba_{i}
\odot c_{i} = \sum_{i} a_{i} \odot b^{op}c_{i}$ for all $b
\in B$. This subspace was first considered by Takeuchi and
provides the right notion of a fiber product for the
algebraic theory of quantum groupoids
\cite{boehm:hopf,xu}.  In the setting of
$C^{*}$-algebras, the Takeuchi product $A
\underset{B}{\times} C$ may be $0$ even when we expect a
nontrivial fiber product on the level of $C^{*}$-algebras;
therefore, the latter can not be obtained as the completion
of the former. In the setting of von Neumann algebras, a
fiber product can be constructed as follows
\cite{sauvageot:2}. If $A$ and $C$ act on Hilbert spaces $H$
and $K$, respectively, one can form the Connes fusion $H
\underset{\mu}{\otimes} K$ with respect to some weight $\mu$
on $B$ and the actions of $B$ on $H$ and $B^{op}$ on $K$
which --- by functoriality --- carries a representation $\pi
\colon A' \odot C' \to {\cal L}(H \underset{\mu}{\otimes}
K)$, and the desired fiber product is $A
\underset{\mu}{\ast} C=\pi(A' \odot C')'$. A categorical
interpretation of this construction is given in
\ref{subsection:vn-category}.

We now modify the last construction to define a fiber product for
$C^{*}$-algebras $A$ and $C$ as follows.
\begin{itemize}
\item[(A)] We assume that $A$ and $C$ are represented on a
$C^{*}$-$\frakb$-module $H_{\beta}$ and a
$C^{*}$-$\frakbo$-module $K_{\gamma}$, respectively, where
$\frakb=\cbasesb$ is a $C^{*}$-base, such that
$\rho_{\beta}(\frakB)$ and $\rho_{\gamma}(\frakBo)$ take
the places of $B$ and $B^{op}$, respectively.
\item[(B)] On the relative tensor product $\HtensorK$, we
define $C^{*}$-algebras $\Ind_{\kgamma{2}}(A)$ and
$\Ind_{\kbeta{1}}(C)$ which, roughly, take the places of
$\pi(A' \odot \Id_{K})'$ and $\pi(\Id_{H} \odot C')'$.
\item[(C)] The fiber product is then the intersection
$\AfibreB=\Ind_{\kgamma{2}}(A) \cap \Ind_{\kbeta{1}}(C)
\subseteq {\cal L}(\HtensorK)$.
\end{itemize}

\subsection{$C^{*}$-algebras represented on $C^{*}$-modules}

\label{subsection:fiber-algebras}

Let $\frakb=\cbasesb$ be a $C^{*}$-base.  As indicated in
step (A), we adopt the following terminology.
\begin{definition} \label{definition:fiber-algebra} A {\em
    $C^{*}$-$\frakBo$-algebra} $(A,\rho)$, briefly written
  $A_{\rho}$, is a $C^{*}$-algebra $A$ with a
  $*$-homo\-morphism $\rho \colon \frakBo \to M(A)$.  A {\em
    morphism} of $C^{*}$-$\frakBo$-algebras $A_{\rho}$ and
  $B_{\sigma}$ is a $*$-homomorphism $\pi \colon A \to B$
  satisfying $\sigma(x)\pi(a) = \pi(\rho(x)a)$ for all $x
  \in \frakBo, a \in A$.  We denote the category of all
  $C^{*}$-$\frakBo$-algebras by $\calg{\frakBo}$.

  A {\em (nondegenerate) $C^{*}$-$\frakb$-algebra} is a pair
  $A_{H}^{\alpha}=(H_{\alpha},A)$, where $H_{\alpha}$ is a
  $C^{*}$-$\frakb$-module, $A \subseteq {\cal L}(H)$ a
  (nondegenerate) $C^{*}$-algebra, and
  $\rho_{\alpha}(\frakBo)A \subseteq A$.  A (semi-)morphism
  between $C^{*}$-$\frakb$-algebras $A_{H}^{\alpha}$,
  $B_{K}^{\beta}$ is a $*$-homomorphism $\pi \colon A \to B$
  satisfying $\beta=[{\cal
    L}_{(s)}^{\pi}(H_{\alpha},K_{\beta})\alpha]$, where
  ${\cal L}_{(s)}^{\pi}(H_{\alpha}, K_{\beta}) := \{ T \in
  {\cal L}_{(s)}(H_{\alpha},K_{\beta}) \mid \forall a \in A:
  Ta=\pi(a)T\}$. We denote the category of all
  $C^{*}$-$\frakb$-algebras together with all
  (semi-)morphisms by $\calg{\frakb}{}^{(s)}$.
  \end{definition}
We first give some examples of $C^{*}$-$\frakb$-algebras
  and then study the relation between $\calg{\frakBo}$ and
  $\calg{\frakb}$.
  \begin{examples} \label{examples:fiber-algebras}
    \begin{enumerate}
    \item If $H$ is a Hilbert space and $A \subseteq {\cal
        L}(H)$  a
      $C^{*}$-algebra, then $A^{\alpha}_{H}$ is a
      $C^{*}$-$\trivbase$-algebra, where
      $\trivbase=(\complex,\complex,\complex)$ denotes the
      trivial $C^{*}$-base and $\alpha={\cal L}(\complex,H)$.
    \item Let $A^{\alpha}_{H}$ be a nondegenerate
      $C^{*}$-$\frakb$-algebra. If we identify $M(A)$ with a
      $C^{*}$-subalgebra of ${\cal L}(H)$ in the canonical
      way, $M(A)^{\alpha}_{H}$ becomes a
      $C^{*}$-$\frakb$-algebra.
    \item Let $({\cal A}_{i})_{i}$ be a family of
      $C^{*}$-$\frakb$-algebras, where ${\cal A}_{i}=({\cal
        H}_{i},A_{i})$ for each $i$. Then the $c_{0}$-sum
      $\bigoplus_{i} A_{i}$ and the $l^{\infty}$-product
      $\prod_{i} A_{i}$ are naturally represented on the
      underlying Hilbert space of $\boxplus_{i} {\cal
        H}_{i}$, and we obtain $C^{*}$-$\frakb$-algebras
      $\boxplus_{i} {\cal A}_{i} := \big(\boxplus_{i} {\cal
        H}_{i}, \bigoplus_{i} A_{i}\big)$
      and $\prod_{i} {\cal A}_{i}:=\big(\boxplus_{i} {\cal
        H}_{i}, \prod_{i} A_{i}\big)$.
      For each $j$, the canonical maps $A_{j} \to
      \bigoplus_{i} A_{i} \to \prod_{i} A_{i} \to A_{j}$ are
      evidently morphisms of $C^{*}$-$\frakb$-algebras
      ${\cal A}_{j} \to \boxplus_{i} {\cal A}_{i} \to
      \prod_{i}{\cal A}_{i} \to {\cal A}_{j}$.
\end{enumerate}
\end{examples}
The following example is a continuation of Example
\ref{example:rtp-bimodule}.
\begin{example}
  Let $B$ be a $C^{*}$-algebra with a KMS-state $\mu$ and
  associated $C^{*}$-base $\frakb$, and let $A$ be a
  $C^{*}$-algebra containing $B$ with a conditional
  expectation $\phi\colon A \to B$ as in Example
  \ref{example:rtp-bimodule}. With the notation introduced
  before, $\pi_{\nu}(A)^{\beta}_{H}$ is a nondegenerate
  $C^{*}$-$\frakb$-algebra because
  $\rho_{\beta}(\frakB)\pi_{\nu}(A) =
  \pi_{\nu}(B)\pi_{\nu}(A) \subseteq \pi_{\nu}(A)$, and
  similarly, $(\pi_{\nu}^{op}(A^{op}))^{\alpha}_{H}$ is a
  nondegenerate $C^{*}$-$\frakbo$-algebra \cite[\S
  2--3]{timmermann:leiden}.
\end{example}
The categories $\calgs{\frakb}$ and $\calg{\frakBo}$ are
related by a pair of adjoint functors, as we shall see now.
\begin{lemma} \label{lemma:fiber-morphism} Let $\pi$ be a
  semi-morphism of $C^{*}$-$\frakb$-algebras
  $A^{\alpha}_{H}$ and $B^{\beta}_{K}$. Then $\pi$ is normal
  and $\pi(a\rho_{\alpha}(x)) = \pi(a)\rho_{\beta}(x)$ for
  all $x \in \frakBo$, $a \in A$.
\end{lemma}
\begin{proof}
  Let $T,T' \in {\cal L}^{\pi}_{s}(H_{\alpha},K_{\beta})$,
  $\xi,\xi' \in \alpha$, $\zeta,\zeta' \in \frakK$, $a \in
  A$, $x \in \frakBo$. Then $\langle T\xi\zeta|\pi(a)
  T'\xi'\zeta'\rangle =\langle \xi\zeta| a
  T^{*}T'\xi'\zeta'\rangle$ and
  $\pi(a\rho_{\alpha}(x))T\xi\zeta =
  Ta\rho_{\alpha}(x)\xi\zeta = \pi(a)T\xi x \zeta =
  \pi(a)\rho_{\beta}(x)T\xi\zeta$ because $T\xi \in \beta$.
  Now, the assertions follow since $K=[{\cal
    L}^{\pi}_{s}(H_{\alpha},K_{\beta})\alpha\frakK]$.
\end{proof}
The preceding lemma shows that there exists a forgetful
functor
\begin{align*}
 \bfU_{\frakb} \colon \calgs{\frakb} \to
\calg{\frakBo}, \quad 
\begin{cases}
  A^{\alpha}_{H} \mapsto A_{\rho_{\alpha}} &\text{for each
    object } A^{\alpha}_{H}, \\
  \pi \mapsto \pi &\text{for each morphism } \pi.
\end{cases}
\end{align*}
We shall see that this functor has a partial adjoint that
associates to a $C^{*}$-$\frakBo$-algebra a universal
representation on a $C^{*}$-$\frakb$-module.  For the
discussion, we fix a $C^{*}$-$\frakBo$-algebra $C_{\sigma}$.
\begin{definition}
  A {\em representation of $C_{\sigma}$ in $\calgs{\frakb}$}
  is a pair $({\cal A},\phi)$, where ${\cal
    A}=A^{\alpha}_{H} \in \calgs{\frakb}$ and $\phi \in
  \calg{\frakBo}(C_{\sigma},\bfU{\cal A})$.  Denote by
  $\bfrep_{\frakb}(C_{_{\sigma}})$ the category of all such
  representations, where the {\em morphisms} between objects
  $({\cal A},\phi)$ and $({\cal B},\psi)$ are all $\pi \in
  \calgs{\frakb}({\cal A},{\cal B})$ satisfying
  $\psi=\bfU\pi\circ \phi$. 
\end{definition}
Note that $\bfrep_{\frakb}(C_{_{\sigma}})$ is just the comma
category $(C_{\sigma} \downarrow \bfU_{\frakb})$
\cite{maclane}.  Unfortunately, we have no general method
like the GNS-construction to produce representations of
$C_{\sigma}$ in  in $\calgs{\frakb}$. In
particular, we have no good criteria to decide whether there
are any and, if so, whether there exists a faithful
one. However, we now show that if there are any
representations, then there also is a universal one. The
proof involves the following direct product construction.
\begin{example}
  Let $({\cal A}_{i},\phi_{i}) \in
  \bfrep_{\frakb}(C_{\sigma})$ for all $i$, where ${\cal
    A}_{i}=({\cal H}_{i},A_{i})$, and define $\phi \colon C
  \to \prod_{i} A_{i}$ by $c \mapsto
  (\phi_{i}(c))_{i}$. Then 
  $\prod_{i}({\cal A}_{i},\phi_{i}) :=( \prod_{i}{\cal
    A}_{i}, \phi) \in \bfrep_{\frakb}(C_{\sigma})$, and the
  canonical maps ${\cal A}_{j} \to \prod_{i} {\cal A}_{i}
  \to {\cal A}_{j}$ are morphisms between $({\cal
    A}_{j},\phi_{j})$ and $(\prod_{i}{\cal A}_{i},\phi)$ for
  each $j$.
\end{example}
\begin{proposition} \label{proposition:fiber-universal} If
  the category $\bfrep_{\frakb}(C_{\sigma})$ is non-empty,
  then it has an initial object.
\end{proposition}
\begin{proof}
  Assume that $\bfrep_{\frakb}(C_{\sigma})$ is non-empty.
  We first use a cardinality argument to show that
  $\bfrep_{\frakb}(C_{\sigma})$ has an initial set of
  objects, and then apply the direct product construction to
  this set to obtain an initial object.

  Given a topological vector space $X$ and a cardinal number
  $c$, let us call $X$ {\em $c$-separable} if $X$ has a
  linearly dense subset of cardinality $c$.  Choose a
  cardinal number $d$ such that $\frakB$ and $C \times
  \frakK$ are $d$-separable, and let $e:=|\naturals|\sum_{n}
  d^{n}$.  Then the isomorphism classes of
  $e$-separable Hilbert $C^{*}$-$\frakB$-modules form a set,
  and hence there exists a set ${\cal R}$ of objects in
  $\bfrep_{\frakb}(C_{\sigma})$ such that each
  $(A^{\alpha}_{H},\phi) \in \bfrep_{\frakb}(C_{\sigma})$
  with $e$-separable $\alpha$ is isomorphic to some element
  of ${\cal R}$. Let $(A^{\alpha}_{H},\phi) = \boxplus_{R
    \in {\cal R}} R$.  We show that
  $(\phi(C)^{\alpha}_{H},\phi)$ is initial in
  $\bfrep_{\frakb}(C_{\sigma})$.

  Let $(B^{\beta}_{K},\psi)
  \in\bfrep_{\frakb}(C_{\sigma})$. We show that there exists
  a morphism $\pi \in \calgs{\frakb}(\phi(C)^{\alpha}_{H},
  B^{\beta}_{K})$ such that $\psi=\pi \circ \phi$, and
  uniqueness of such a $\pi$ is evident.  Let $\xi \in
  \beta$ be given. Since $\frakB$ and $C \times \frakK$ are
  $d$-separable, we can inductively choose subspaces
  $\beta_{0} \subseteq \beta_{1} \subseteq \cdots \subseteq
  \beta$ and cardinal numbers $d_{0}, d_{1}, \ldots$
  such that $\xi \in \beta_{0}$,
  $[\beta_{0}^{*}\beta_{0}]=\frakB$, $d_{0} \leq 2d+1$,
  $\beta_{0}$ is $d_{0}$-separable and
  for all $n \geq 0$,
  \begin{align*}
      \beta_{n}\frakB &\subseteq \beta_{n+1}, &
      \psi(C)\beta_{n}\frakK &\subseteq [\beta_{n+1}\frakK],
      & d_{n+1} &\leq |\naturals| d d_{n}, & \beta_{n+1}
      &\text{ is $d_{n+1}$-separable}.
    \end{align*}
    Let $\tilde \beta:=[\bigcup_{n} \beta_{n}] \subseteq
    \beta$ and $ \tilde K:=[\tilde\beta\frakK] \subseteq K$.
    By construction, $[\tilde \beta^{*}\tilde\beta]=\frakB$,
    $ \tilde\beta\frakB \subseteq \tilde\beta$, $
    \psi(C)\tilde K \subseteq \tilde K$, so that
    $(\psi(C)|_{\tilde K})^{\tilde \beta}_{\tilde K}$ is in
    $\calg{\frakb}$. Define $\tilde \psi \colon C \to
    \psi(C)|_{\tilde K}$ by $c \mapsto \psi(c)|_{\tilde
      K}$. Then $(\tilde \psi(C)^{\tilde \beta}_{\tilde
      K},\tilde\psi)$ is in
    $\bfrep_{\frakb}(C_{\sigma})$. Since $\tilde \beta$ is
    $e$-separable, $(\tilde \psi(C)^{\tilde \beta}_{\tilde
      K},\tilde\psi)$ is isomorphic to some element of
    ${\cal R}$. Hence, there exists a surjection $\tilde T
    \colon H \to \tilde K$ such that $\tilde T\alpha=\tilde
    \beta$, and the composition with the inclusion $\tilde K
    \to K$ gives an operator $T \in {\cal
      L}_{s}(H_{\alpha},K_{\beta})$ such that
    $\psi(c)T=T\phi(c)$ for all $c \in C$. Since $\xi \in
    \tilde \beta = T\alpha$ and $\xi \in \beta$ was
    arbitrary, we can conclude the existence of $\pi$ as
    desired.
\end{proof}

Evidently, every morphism $\Phi$ between
$C^{*}$-$\frakBo$-algebras $C_{\sigma}$ and $D_{\tau}$
yields a functor
\begin{align*}
  \Phi^{*} \colon \bfrep_{\frakb}(D_{\tau}) \to
  \bfrep_{\frakb}(C_{\sigma}), \quad
  \begin{cases}
   (A^{\alpha}_{H},\phi) \mapsto (A^{\alpha}_{H}, \phi \circ
   \Phi) &\text{for each object } (A^{\alpha}_{H},\phi), \\
   \pi \mapsto \pi &\text{for each morphism } \pi.
 \end{cases}
\end{align*}
Denote by $\calgr{\frakBo}$ the full subcategory of
$\calg{\frakBo}$  consisting of all objects $C_{\sigma}$
for which  $\bfrep(C_{\sigma})$ is non-empty.
\begin{theorem} \label{theorem:fiber-algebras-adjoint} There
  exist a functor $\bfR_{\frakb} \colon \calgr{\frakBo} \to
  \calgs{\frakb}$ and natural transformations $\eta \colon
  \Id_{\calgr{\frakBo}} \to \bfU_{\frakb}\bfR_{\frakb}$ and
  $\epsilon \colon \bfR_{\frakb}\bfU_{\frakb} \to
  \Id_{\calgs{\frakb}}$ such that for every $C_{\sigma},
  D_{\tau} \in \calgr{\frakBo}$,  $\Phi \in
  \calgr{\frakBo}(C_{\sigma},D_{\tau})$, 
  $A^{\alpha}_{H} \in \calgs{\frakb}$,
  \begin{itemize}
  \item $\bfR_{\frakb}(C_{\sigma})\in
    \bfrep_{\frakb}(C_{\sigma})$ is an initial object and
    $\bfR_{\frakb}(\Phi)$ is the unique morphism from
    $\bfR_{\frakb}(C_{\sigma})$ to
    $\Phi^{*}(\bfR_{\frakb}(D_{\tau}))$,
  \item $\eta_{C_{\sigma}}=\phi$ if
    $\bfR_{\frakb}(C_{\sigma}) = (B^{\beta}_{K},\phi)$, and
    $\epsilon_{A^{\alpha}_{H}}$ is the unique morphism from
    $\bfR_{\frakb}\bfU_{\frakb} (A^{\alpha}_{H})$ to
    $(A^{\alpha}_{H},\Id_{A})$.
  \end{itemize}
  Moreover, $\bfR_{\frakb}$ is left adjoint to
  $\bfU_{\frakb}$ and $\eta$, $\epsilon$ are the unit and
  counit of the adjunction, respectively.
\end{theorem}
\begin{proof}
  This follows from Proposition
  \ref{proposition:fiber-universal} and \cite[\S IV Theorem 2]{maclane}.
\end{proof}

We next consider $C^{*}$-algebras represented on
$C^{*}$-bimodules.  Let $\fraka=\cbasesa$ be a $C^{*}$-base.
\begin{definition}
  A {\em $C^{*}$-$(\frakA,\frakBo)$-algebra} is a triple
  $(A,\rho,\sigma)$, briefly written $A_{\rho,\sigma}$,
  where $A_{\rho}$ is a $C^{*}$-$\frakA$-algebra,
  $A_{\sigma}$ a $C^{*}$-$\frakBo$-algebra, and
  $[\rho(\frakA),\sigma(\frakBo)]=0$.  A {\em morphism} of
  $C^{*}$-$(\frakA,\frakBo)$-algebras is a morphism of the
  underlying $C^{*}$-$\frakA$-algebras and
  $C^{*}$-$\frakBo$-algebras.  We denote the category of all
  $C^{*}$-$(\frakA,\frakBo)$-algebras by
  $\calg{(\frakA,\frakBo)}$.

  A {\em (nondegenerate) $C^{*}$-$(\frakao,\frakb)$-algebra}
  is a pair $A_{H}^{\alpha,\beta}=(\aHb,A)$, where $\aHb$ is
  a $C^{*}$-$(\frakao,\frakb)$-module, $A^{\alpha}_{H}$ a
  (nondegenerate) $C^{*}$-$\frakao$-algebra, and
  $A^{\beta}_{H}$ a $C^{*}$-$\frakb$-algebra.  A
  (semi-)morphism of
  $C^{*}$-$(\frakao,\frakb)$-algebras $A^{\alpha,\beta}_{H}$
  and $B^{\gamma,\delta}_{K}$ is a $*$-homomorphism $\pi
  \colon A \to B$ satisfying $\gamma=[{\cal
    L}_{(s)}^{\pi}(\aHb,\cKd)\alpha]$ and $\delta=[{\cal
    L}_{(s)}^{\pi}(\aHb,\cKd)\beta]$, where ${\cal
    L}_{(s)}^{\pi}(\aHb,\cKd) := \{ T \in {\cal
    L}_{(s)}(\aHb,\cKd) \mid \forall a \in A:
  Ta=\pi(a)T\}$. We denote the category of all
  $C^{*}$-$(\frakao,\frakb)$-algebras together with all
  (semi-)morphisms by $\calgss{(\frakao,\frakb)}$.
  \end{definition}
  \begin{remark}
    Note that the condition on a (semi-)morphism between
    $C^{*}$-$(\frakao,\frakb)$-algebras above is stronger
    than just being a (semi-)morphism of the underlying
    $C^{*}$-$\frakao$-algebras and $C^{*}$-$\frakb$-algebras.
  \end{remark}
  Examples \ref{examples:fiber-algebras} ii) and iii)
  naturally extend to $C^{*}$-$(\frakao,\frakb)$-algebras, and  the
  categories $\calg{(\frakA,\frakBo)}$ and
  $\calgs{(\frakao,\frakb)}$ are again related by a pair of
  adjoint functors.
  \begin{theorem} \label{theorem:fiber-algebras-adjoint-2}
    There exists a functor $\bfU_{(\frakao,\frakb)} \colon
    \calgs{(\frakao,\frakb)} \to \calg{(\frakA,\frakBo)}$,
    given by $A^{\alpha,\beta}_{H} \mapsto
    A_{\rho_{\alpha},\rho_{\beta}}$ on objects and $\pi
    \mapsto \pi$ on morphisms.  Denote by
    $\calgr{(\frakA,\frakBo)}$ the full subcategory of
    $\calg{(\frakA,\frakBo)}$ consisting of all objects
    $C_{\sigma,\rho}$ for which the comma category
    $(C_{\sigma,\rho} \downarrow \bfU_{(\frakao,\frakb)})$ is
    non-empty. Then the corestriction of
    $\bfU_{(\frakao,\frakb)}$ to $\calgr{(\frakA,\frakBo)}$
    has a left adjoint $\bfR_{(\frakao,\frakb)} \colon
    \calgr{(\frakA,\frakBo)} \to \calgs{(\frakao,\frakb)}$.
\end{theorem}
\begin{proof}
  The proof proceeds as in the case of
  $C^{*}$-$\frakb$-algebras with straightforward
  modifications, so we only indicate the necessary changes
  for the second half of the proof of Proposition
  \ref{proposition:fiber-universal}. Given a
  $C^{*}$-$(\frakA,\frakBo)$-algebra $C_{\sigma,\tau}$ and
  a $C^{*}$-$(\frakao,\frakb)$-algebra
  $B^{\gamma,\delta}_{K}$ with a morphism $\psi \colon
  C_{\sigma,\tau} \to B_{\rho_{\gamma},\rho_{\delta}}$, one
  constructs $\tilde \gamma \subseteq \gamma $ and $\tilde
  \delta \subseteq \delta$ for given $\xi \in \gamma$, $\eta
  \in \delta$ as follows. One first fixes a cardinal number $d$
  such that $\frakA,\frakAo,\frakH,\frakB,\frakBo,\frakH$
  are $d$-separable, and then inductively chooses cardinal
  numbers $d_{0},d_{1},\ldots$ and closed subspaces
  $\gamma_{0} \subseteq \gamma_{1} \subseteq \cdots
  \subseteq \gamma$ and $\delta_{0} \subseteq \delta_{1}
  \subseteq \cdots \subseteq \delta$ such that
  \begin{gather*}
    \begin{aligned}
      \xi &\in \gamma_{0}, & \eta &\in \delta_{0}, &
      [\gamma_{0}^{*}\gamma_{0}] &= \frakAo, &
      [\delta_{0}^{*}\delta_{0}]&=\frakB, & d_{0} &\leq
      2d+1, & \gamma_{0}, \delta_{0} &\text{ are
        $d_{0}$-separable},
    \end{aligned} \\
    \begin{aligned}
      \rho_{\delta}(\frakBo)\gamma_{n}+ \gamma_{n} \frakAo
      &\subseteq \gamma_{n+1}, &
      \rho_{\gamma}(\frakA)\gamma_{n} + \delta_{n} \frakB
      &\subseteq \delta_{n+1}, &
      \psi(C)\gamma_{n}\frakH+\psi(C)\delta_{n} \frakK
      &\subseteq [\gamma_{n+1} \frakH]\cap
      [\delta_{n+1}\frakK],
    \end{aligned} \\
    \begin{aligned}
      d_{n+1} &\leq |\naturals|d^{2} d_{n}, &
      \gamma_{n+1}, \delta_{n+1} &\text{ are $d_{n+1}$-separable}    \end{aligned}
  \end{gather*}
  for all $n \geq 0$, and finally lets $\tilde
  \gamma:=[\bigcup_{n} \gamma_{n}], \tilde
  \delta:=[\bigcup_{n} \delta_{n}], \tilde K:=[\tilde\gamma
  \frakH]=[\tilde \delta \frakK]$.
\end{proof}
\begin{remark} \label{remark:fiber-adjoints}
  Let $C_{\rho,\sigma}$ be a
  $C^{*}$-$(\frakA,\frakBo)$-algebra, 
  $A^{\alpha,\beta}_{H} =
  \bfR_{(\frakao,\frakb)}(C_{\rho,\sigma})$, and  $\phi =
  \eta_{C_{\rho,\sigma}} \colon C_{\rho,\sigma}\to
  A_{\rho_{\alpha},\rho_{\beta}}$  the morphism given by
  the unit of the adjunction above. Then $(A^{\alpha},\phi)
  \in \bfrep_{\frakao}(C_{\rho})$ and $(A^{\beta},\phi) \in
  \bfrep_{\frakb}(C_{\sigma})$, whence we have
  semi-morphisms $\bfR_{\frakao}(C_{\sigma}) \to
  A^{\alpha}_{H}$ and $\bfR_{\frakb} (C_{\rho}) \to
  A^{\beta}_{H}$. 
\end{remark}

\subsection{The spatial fiber product for $C^{*}$-algebras
  represented on $C^{*}$-modules}

\label{subsection:fiber-fiber}

Our definition of the fiber product of $C^{*}$-algebras
represented on $C^{*}$-modules --- more precisely, step (B) in
the introduction --- involves the following construction.

Let $H$ and $K$ be Hilbert spaces, $I \subseteq {\cal
  L}(H,K)$ a subspace and $A \subseteq {\cal L}(H)$ a
$C^{*}$-algebra such that $ [IH]=K$, $[I^{*}K]=H$,
$[II^{*}I]=I$, $I^{*}IA \subseteq A$.  We define a new
$C^{*}$-algebra
\begin{align*}
  \Ind_{I}(A):= \{ T \in {\cal L}(K) \mid TI + T^{*}I
  \subseteq [I A]\} \subseteq {\cal L}(K).
\end{align*}
\begin{definition} \label{definition:fiber-topologies} The {\em
    $I$-strong-$*$, $I$-strong, and $I$-weak topology} on
  ${\cal L}(K)$ are the topologies induced by the families
  of semi-norms $T \mapsto \|T\xi\| + \|T^{*}\xi\|$ ($\xi \in
  I$), $T \mapsto \| T\xi\|$ $(\xi \in I)$, and $T \mapsto
  \| \xi^{*}T\xi'\|$ $(\xi,\xi' \in I)$, respectively.
  Given a subset $X \subseteq {\cal L}(K)$, denote by
  $[X]_{I}$ the closure of $\lspan X$ with respect to the
  $I$-strong-$*$ topology.
\end{definition}
Evidently, the multiplication in ${\cal L}(K)$ is separately
continuous with respect to the topologies introduced above,
and the involution $T \mapsto T^{*}$ is continuous with
respect to the $I$-strong-$*$ and the $I$-weak topology.
Define $\rho_{I} \colon (I^{*}I)' \to {\cal L}(K)$ as in Lemma
\ref{lemma:rtp-morphism}.
\begin{lemma}
  \label{lemma:fiber-ind}
  \begin{enumerate}
  \item $[I^{*}\Ind_{I}(A)I] \subseteq A$ and $\Ind_{I}(A) = [I
    AI^{*}]_{I}$.
  \item $\Ind_{I}(M(A)) \subseteq M(\Ind_{I}(A))$.
  \item $\Ind_{I}(A) \subseteq {\cal L}(K)$ is nondegenerate
    if and only if $A \subseteq {\cal L}(H)$ is
    nondegenerate.
  \item If $A \subseteq {\cal L}(H)$ is nondegenerate, then $A' \subseteq
    (I^{*}I)'$ and $\Ind_{I}(A) \subseteq \rho_{I}(A')'$.
  \end{enumerate}
\end{lemma}
\begin{proof}
  i) We have $[I^{*}\Ind_{I}(A)I] \subseteq [I^{*}IA]
  \subseteq A$ by definition and $[I AI^{*}]_{I} \subseteq
  \Ind_{I}(A)$ because $[I A I^{*}]_{I} I \subseteq [I A
  I^{*}I] \subseteq [I A]$.  To see that $[I AI^{*}]_{I}
  \supseteq \Ind_{I}(A)$, choose a bounded approximate unit
  $(u_{\nu})_{\nu}$ for the $C^{*}$-algebra $[II^{*}]$ and
  observe that for each $T \in \Ind_{I}(A)$, the net
  $(u_{\nu} T u_{\nu})_{\nu}$ lies in the space $[II^{*}
  \Ind_{I}(A) II^{*}] \subseteq [I A I^{*}]$ and converges
  to $T$ in the $I$-strong-$*$ topology because $\lim_{\nu}
  T^{(*)}u_{\nu}\xi =T^{(*)}\xi \in [IA]$ for all $\xi \in I$ and
  $\lim_{\nu} u_{\nu}\omega=\omega$ for all $\omega \in
  [IA]$.

  ii) If $S \in \Ind_{I}(M(A))$, $T \in \Ind_{I}(A)$, then
  $ST \in \Ind_{I}(A)$ because $STI \subseteq [S I A]
  \subseteq [I M(A)A]=[I A]$ and $T^{*}S^{*} I \subseteq [TI
  M(A)] \subseteq [I AM(A)]=[I A]$.

  iii) If $\Ind_{I}(A) \subseteq {\cal L}(K)$ is
  nondegenerate, then $[A H] \supseteq [I^{*}\Ind_{I}(A)IH]
  = [I^{*} \Ind_{I}(A)K] = [I^{*}K] = H$.  Conversely, if
  $A$ is nondegenerate, then $[IAI^{*}]$ and hence also
  $\Ind_{I}(A)$ is nondegenerate.

  iv) Assume that $A$ is nondegenerate. Then $I^{*}I
  \subseteq M(A) \subseteq {\cal L}(H)$ and hence $A'
  \subseteq (I^{*}I)'$.  For all $x \in \Ind_{I}(A)$, $y \in
  A'$, $S,T \in I$, we have $S^{*}x\rho_{I}(y)T = S^{*}xTy =
  yS^{*}xT= S^{*}\rho_{I}(y)xT$ because $S^{*}xT\in A$, and
  since $[IH]=K$, we can conclude that
  $x\rho_{I}(y)=\rho_{I}(y)x$.
\end{proof}

Let $ \frakb=\cbasesb$ be a $C^{*}$-base, $A^{\beta}_{H}$ a
$C^{*}$-$\frakb$-algebra, and $B^{\gamma}_{K}$ a
$C^{*}$-$\frakbo$-algebra. We apply the construction above
to $A$, $B$ and $\kgamma{2} \subseteq {\cal L}(H,\HfibreK)$,
$\kbeta{1} \subseteq {\cal L}(K,\HfibreK)$, respectively,
and define the {\em fiber product} of $A^{\beta}_{H}$ and
$B^{\gamma}_{K}$ to be the $C^{*}$-algebra
\begin{align*}
  \AfibreB &:= \ind{\kgamma{2}}{A} \cap \ind{\kbeta{1}}{B} \\
  & = \{ T \in {\cal L}(\HfibreK) \mid T\kgamma{2} +
  T^{*}\kgamma{2} \subseteq [\kgamma{2}A], T\kbeta{1} +
  T^{*}\kbeta{1} \subseteq [\kbeta{1}B] \}.
\end{align*}
The spaces of operators involved are visualized as arrows
in the following diagram:
\begin{align*}
  \xymatrix@C=40pt@R=12pt{ H \ar[d]_{A}
    \ar[r]^(0.4){\kgamma{2}} &
    {\HfibreK} \ar[d]^{\AfibreB}& K \ar[l]_(0.4){\kbeta{1}} \ar[d]^{B} \\
    H \ar[r]^(0.4){\kgamma{2}} & {\HfibreK} & K
    \ar[l]_(0.4){\kbeta{1}} }
\end{align*}
Even in very special situations, it seems to be difficult to
give a more explicit description of the fiber product.  The
main drawback of the definition above is that apart from
special situations, we do not know how to produce elements
of the fiber product.

Let $\fraka=\cbasesa$ and $\frakc=\cbasesc$ be further
$C^{*}$-bases.
\begin{proposition} \label{proposition:fiber-fiber} Let
  ${\cal A}=A^{\alpha,\beta}_{H}$ be a
$C^{*}$-$(\frakao,\frakb)$-algebra and ${\cal
B}=B^{\gamma,\delta}_{K}$ a
$C^{*}$-$(\frakbo,\frakc)$-algebra. Then ${\cal A} \bfibre
  {\cal B}:=(\aHb \btensor \cKd, \AfibreB)$ is a
  $C^{*}$-$(\frakao,\frakc)$-algebra.
\end{proposition}
\begin{proof}
  The product $X:= \rho_{(\alpha \lt
    \gamma)}(\frakAo)(\AfibreB)$ is contained in $\AfibreB$
  because
  \begin{align*}
    X \kbeta{1} &\subseteq
    [|\rho_{\alpha}(\frakA)\beta\rangle_{1}B] =
    [\kbeta{1}B], & X^{*} \kbeta{1} & =
    (\AfibreB)|\rho_{\alpha}(\frakA)\beta\rangle_{1}
    \subseteq
    [\kbeta{1}B], \\
    X\kgamma{2} &\subseteq [
    \kgamma{2}\rho_{\alpha}(\frakA)A] \subseteq
    [\kgamma{2}A], & X^{*}\kgamma{2} &= (\AfibreB)
    \kgamma{2}\rho_{\alpha}(\frakA) \subseteq [\kgamma{2}A]
  \end{align*}
  by equation \eqref{eq:rtp-rtp-rho}.  A similar argument
  shows that $\rho_{(\beta \rt \delta)}(\frakCo)(\AfibreB)
  \subseteq \AfibreB$.
\end{proof}
In the situation above, we call ${\cal A} \bfibre {\cal B}$
the {\em fiber product} of ${\cal A}$ and ${\cal B}$.
Forgetting $\alpha$ or $\delta$, we obtain a
$C^{*}$-$\frakc$-algebra $\AfibreB_{\delta}:=A^{\beta}_{H} \bfibre B^{\gamma,\delta}_{H}:=
(\HtensorK_{\delta},\AfibreB)$ and a $C^{*}$-$\frakao$-algebra
${_{\alpha}\AfibreB}=A^{\alpha,\beta}_{H} \bfibre B^{\gamma}_{K}$.

Denote by $A' \subseteq {\cal L}(H)$ and $B' \subseteq {\cal
  L}(K)$ the commutants of $A$ and $B$, respectively, and let
\begin{align*}
  \begin{aligned}
    A^{(\beta)} &:=A \cap {\cal L}(H_{\beta}), &
    B^{(\gamma)} &:=B \cap {\cal L}(K_{\gamma}), & X
    &:=(A^{(\beta)}\btensor \Id)+ (\Id \btensor
    B^{(\gamma)}),
  \end{aligned} \\
  M_{s}(A^{(\beta)} \btensor B^{(\gamma)}) := \{ T \in {\cal
    L}(\HfibreK) \mid TX, XT\subseteq A^{(\beta)} \btensor
  B^{(\gamma)}\}.
\end{align*}
\begin{lemma} \label{lemma:fiber-properties}
  \begin{enumerate}
  \item $\bbeta{1}(\AfibreB)\kbeta{1} \subseteq B$,
    $\bgamma{2}(\AfibreB)\kgamma{2} \subseteq A$, and $M(A)
    \fibrebc M(B) \subseteq M(\AfibreB)$.
  \item $A^{(\beta)} \btensor B^{(\gamma)} \subseteq
    \AfibreB$.
  \item If $[A^{(\beta)}\beta]=\beta$ and
    $[B^{(\gamma)}\gamma]=\gamma$, then $\AfibreB$ is
    nondegenerate and $M_{s}(A^{(\beta)} \btensor
    B^{(\gamma)}) \subseteq \AfibreB$.
  \item If $\rho_{\beta}(\frakBo) \subseteq A$, then
    $\Id_{H} \btensor B^{(\gamma)} \subseteq \AfibreB$. If
    $\rho_{\gamma}(\frakB) \subseteq B$, then $A^{(\beta)}
    \btensor \Id_{K} \subseteq \AfibreB$.
  \item   $\Id_{(\HfibreK)}\in \AfibreB$ if and only if
    $\rho_{\beta}(\frakBo) \subseteq A$ and
    $\rho_{\gamma}(\frakB) \subseteq B$.
  \item If $A^{\alpha,\beta}_{H}$ is a
    $C^{*}$-$(\frakao,\frakb)$-algebra and
    $B^{\gamma,\delta}_{K}$ a
    $C^{*}$-$(\frakbo,\frakc)$-algebra such that
    $\rho_{\alpha}(\frakA)+\rho_{\beta}(\frakBo) \subseteq
    A$ and $\rho_{\gamma}(\frakB)+\rho_{\delta}(\frakCo)
    \subseteq B$, then $\rho_{(\alpha \lt
      \gamma)}(\frakA)+\rho_{(\beta \rt \delta)}(\frakCo)
    \subseteq \AfibreB$.
  \item If $\AfibreB$ is nondegenerate, then the
  $C^{*}$-algebra $[\beta^{*}A\beta] \cap
  [\gamma^{*}B\gamma] \subseteq {\cal L}(\frakK)$ is
  nondegenerate.
\item If $A$ and $B$ are nondegenerate, then $A' \subseteq
  \rho_{\beta}(\frakBo)'$, $B' \subseteq
  \rho_{\gamma}(\frakB)'$, and $\AfibreB \subseteq
  \rho_{\kgamma{2}}(A') \cap \rho_{\kbeta{1}}(B') = (A'
  \btensor \Id_{K})' \cap (\Id_{H} \btensor B')'$.
  \end{enumerate}
\end{lemma}
\begin{proof}
  i) Immediate from Lemma \ref{lemma:fiber-ind}.

  ii) Use $(A^{(\beta)} \btensor B^{(\gamma)})\kbeta{1}
  \subseteq [|A^{(\beta)}\beta\rangle_{1}B^{(\gamma)}]
  \subseteq [\kbeta{1}B]$, $(A^{(\beta)} \btensor
  B^{(\gamma)})\kgamma{2} \subseteq
  [|B^{(\gamma)}\gamma\rangle_{1}A^{(\beta)}] \subseteq
  [\kgamma{2}A]$.

  iii) Assume $[A^{(\beta)}\beta]=\beta$ and
  $[B^{(\gamma)}\gamma]=\gamma$. Then $A^{(\beta)} \btensor
  B^{(\gamma)} \subseteq \AfibreB$ is nondegenerate and for
  each $T
  \in M_{s}(A^{(\beta)} \btensor B^{(\gamma)})$, we have
  $T\kbeta{1} \subseteq [T(A^{(\beta)} \btensor \Id)
  \kbeta{1}] \subseteq [(A^{(\beta)} \btensor
  B^{(\gamma)})\kbeta{1}] \subseteq [\kbeta{1}B]$ and
  similarly $T^{*}\kbeta{1} \subseteq [\kbeta{1}B]$,
  $T\kgamma{2}+T^{*}\kgamma{2} \subseteq [\kgamma{2}A]$.

  iv) If $\rho_{\gamma}(\frakB) \subseteq B$, then
  $(A^{(\beta)} \btensor \Id_{K})\kgamma{2} =
  \kgamma{2}A^{(\beta)}$ and $[(A^{(\beta)} \btensor
  \Id_{K})\kbeta{1}]\subseteq \kbeta{1} = [|\beta
  \frakB\rangle_{1}] = [\kbeta{1}\rho_{\gamma}(\frakB)]
  \subseteq [\kbeta{1}B]$.  The second assertion follows
  similarly.

  v) If $\Id_{(\HfibreK)} \in \AfibreB$, then
  $\rho_{\beta}(\frakBo) = [\bgamma{2}\kgamma{2}] \subseteq
  A$, $\rho_{\gamma}(\frakB) =[\bbeta{1}\kbeta{1}]
  \subseteq B$ by i). Conversely, if the last two inclusions
  hold, then $\kgamma{2} = [|\gamma \frakBo\rangle_{2}] =
  [\kgamma{2}\rho_{\beta}(\frakBo)] \subseteq [\kgamma{2}A]$
  and similarly $\kbeta{1} \subseteq [\kbeta{1}B]$, whence
  $\Id_{(\HfibreK)} \in \AfibreB$.

  vi) Immediate from iv).

  vii) The $C^{*}$-algebra $C:= [\beta^{*}A\beta] \cap
  [\gamma^{*}B\gamma]$ contains  $ \beta^{*}
  \bgamma{2} (\AfibreB) \kgamma{2}\beta = \gamma^{*}
  \bbeta{1}(\AfibreB) \kbeta{1}\gamma$. If $\AfibreB$ is
  nondegenerate, we therefore must have $[C \frakK]
  \supseteq [\beta^{*} \bgamma{2}
  (\AfibreB)(\HfibreK)]=\frakK$.

  viii) Immediate  from Lemma \ref{lemma:fiber-ind}.
\end{proof}
Even in the case of a trivial $C^{*}$-base, we have no 
explicit description of the fiber product.
\begin{examples} \label{examples:fiber} Let $H$ and $K$ be
  Hilbert spaces, $\beta={\cal L}(\complex,H)$,
  $\gamma={\cal L}(\complex,K)$, $\frakb=\trivbase$ the
  trivial $C^{*}$-base $(\complex,\complex,\complex)$, and
  identify $\HtensorK$ with $H \otimes K$ as in Example
  \ref{examples:rtp-rtp}.
  \begin{enumerate}
  \item Let $A \subseteq {\cal L}(H)$ and $B \subseteq {\cal
      L}(K)$ be nondegenerate $C^{*}$-algebras. Then
    $A^{(\beta)} = A$, $B^{(\gamma)}=B$, and by Lemma
    \ref{lemma:fiber-properties}, $\AfibreB$ contains the
    minimal tensor product $A \otimes B \subseteq {\cal L}(H
    \otimes K)$ and $M_{s}(A \otimes B)=\{ T \in {\cal L}(H
    \otimes K) \mid T^{(*)}(1 \otimes B), T^{(*)}(A \otimes
    1) \subseteq A \otimes B\}$. If $A$ or $B$ is
    non-unital, then $\Id_{H\otimes K} \not\in \AfibreB$ by
    Lemma \ref{lemma:fiber-properties} and so $M(A \otimes
    B) \not \subseteq \AfibreB$. In Example
    \ref{examples:fp-commutative} iii), we shall see that also
    $\AfibreB \nsubseteq M(A \otimes B)$ is possible.
  \item Assume that $H=K=l^{2}(\naturals)$ and identify
    $\beta=\gamma={\cal L}(\complex,H)$ with $H$. Then the
    flip $\Sigma \colon H \otimes H \to H \otimes H$, $\xi
    \otimes \eta \mapsto \eta \otimes \xi$, is not contained
    in ${\cal L}(H) \fibre{\beta}{\frakb}{\gamma} {\cal
      L}(H)$.  Indeed, let $(\xi_{\nu})_{\nu}$ be an
    orthonormal basis for $H$ and let $\eta \in H$ be
    non-zero. Then $\langle
    \xi_{\nu}|_{1}\Sigma|\eta\rangle_{1}=|\eta\rangle\langle
    \xi_{\nu}|$ for each $\nu$ and hence $\sum_{\nu} \langle
    \xi_{\nu}|_{1}\Sigma |\eta\rangle_{1}$ does not converge
    in norm. On the other hand, one easily verifies that
    $\sum_{\nu} \langle \xi_{\nu}|_{1}S$  converges in
    norm for each $S \in [|H\rangle_{1}{\cal L}(H)]$.
    Hence, $\Sigma|\eta\rangle_{1} \not\in
    [|H\rangle_{1}{\cal L}(H)]$. 
  \end{enumerate}
\end{examples}

\subsection{Functoriality
  and slice maps}

\label{subsection:fiber-functorial}

We first show that the fiber product constructed above is
functorial, and then consider various slice maps. The
results concerning functoriality were stated in slightly
different form in
\cite{timmermann:leiden,timmermann:coactions,timmermann:cpmu-hopf}
with proofs referring to unpublished material. We use the
opportunity to rectify this situation.  As before, let
$\fraka=\cbasesa,\frakb=\cbasesb,\frakc=\cbasesc$ be
$C^{*}$-bases.

\begin{lemma}
 \label{lemma:fiber-morphism-ind}
 Let $\pi$ be a (semi-)morphism of $C^{*}$-$\frakb$-algebras
 $A_{H}^{\beta}$ and $C^{\lambda}_{L}$,  let $\cKd$ be a
 $C^{*}$-$(\frakbo,\frakc)$-module, and let $ I:=\mathcal{
   L}_{(s)}^{\pi}(H_{\beta},L_{\lambda}) \btensor \Id
 \subseteq {\cal L}(\HtensorK, L
 \rtensor{\lambda}{\frakb}{\gamma} K)$.
   \begin{enumerate}
   \item ${\cal X}:=(\HtensorK_{\delta}, (I^{*}I)')$ and
     ${\cal Y}:=(L \rtensor{\lambda}{\frakb}{\gamma}
     K_{\delta},(II^{*})')$ are nondegenerate
     $C^{*}$-$\frakc$-algebras.
   \item There exists a unique $\rho_{I} \in
     \Mor_{(s)}({\cal X},{\cal Y})$ such that $\rho_{I}(x)S
     =Sx$ for all $x\in (I^{*}I)', S\in I$.
   \item There exists a unique linear contraction
     $j_{\pi}\colon [\kgamma{2}A] \to [\kgamma{2}C]$ given
     by $|\eta\rangle_{2}a \mapsto |\eta\rangle_{2}\pi(a)$.
   \item $\Ind_{\kgamma{2}}(A) \subseteq (I^{*}I)'$ and
     $\rho_{I}(x)|\eta\rangle_{2}=j_{\pi}(x|\eta\rangle_{2})$
     for all $x \in \Ind_{\kgamma{2}}(A)$, $\eta \in
     \gamma$.
 \item Let $B^{\gamma}_{K}$ be a $C^{*}$-$\frakbo$-algebra.
   Then $\AfibreB \subseteq (I^{*}I)'$ and $\rho_{I}(\AfibreB)
   \subseteq \Cl \bfibre \gB$.
 \end{enumerate}
\end{lemma}
\begin{proof} 
  i) Clearly, $(I^{*}I)'$ and $(II^{*})'$ are nondegenerate
  $C^{*}$-algebras, and ${\cal X}$ and ${\cal Y}$ are
  $C^{*}$-$\frakc$-algebras because $\rho_{(\beta \rt
    \delta)}(\frakCo) = \Id \rtensor{\beta}{\frakb}{\gamma}
  \rho_{\delta}(\frakCo) \subseteq (I^{*}I)'$ and
  $\rho_{(\lambda \rt \delta)}(\frakCo) = \Id
  \rtensor{\lambda}{\frakb}{\gamma} \rho_{\delta}(\frakCo)
  \subseteq (II^{*})'$.

  ii) There exists a unique $*$-homomorphism $\rho_{I}\colon
  (I^{*}I)' \to (II^{*})'$ satisfying the formula above by
  Lemma \ref{lemma:rtp-morphism}, and this is a
  (semi-)morphism because $[I(\beta \rt \delta)]= [\lambda
  \rt \delta]$ by assumption on $\pi$.
  
  iii) Let $\eta_{1},\ldots,\eta_{n} \in \gamma$ and
  $a_{1},\ldots, a_{n} \in A$. Then $ \| \sum_{j}
  |\eta_{j}\rangle_{2}\pi(a_{j})\|^{2} = \|\sum_{i,j}
  \pi(a_{i}^{*})\rho_{\lambda}(\eta_{i}^{*}\eta_{j})\pi(a_{j})\|
  \leq \|\sum_{i,j}
  a_{i}^{*}\rho_{\beta}(\eta_{i}^{*}\eta_{j})a_{j}\|=
  \|\sum_{j} |\eta_{j}\rangle_{2}a_{j}\|^{2}$ by Lemma
  \ref{lemma:fiber-morphism}. The claim follows.

  iv) The first assertion follows  from Lemma
  \ref{lemma:fiber-ind} and the relation $I^{*}I
  \subseteq A' \btensor \Id =
  \rho_{|\gamma\rangle_{2}}(A')$, and the second one from
  the fact that  for all $x \in \Ind_{\kgamma{2}}(A),\eta\in
  \gamma,S \in {\cal L}^{\pi}_{(s)}(H_{\beta},L_{\lambda})$,
  we have $\rho_{I}(x)|\eta\rangle_{2}S = \rho_{I}(x)
  (S \btensor \Id)|\eta\rangle_{2} =   (S \btensor \Id) x
  |\eta\rangle_{2} =  j_{\pi}(x|\eta\rangle_{2})S$.

  v) First, $\AfibreB \subseteq (I^{*}I)'$ by Lemma
  \ref{lemma:fiber-properties}.  The second assertion
  follows from the relations
  \begin{align*}
    \rho_{I}(\AfibreB)\kgamma{2} &\subseteq
    \rho_{I}(\Ind_{\kgamma{2}}(A))\kgamma{2} \subseteq
    j_{\pi}([\kgamma{2}A]) = [\kgamma{2}C], \\
    \rho_{I}(\AfibreB)|\lambda\rangle_{1} &=
    \rho_{I}(\AfibreB)[I\kbeta{1}] \subseteq [I (\AfibreB)
    \kbeta{1}] \subseteq[I\kbeta{1}B] =
    [|\lambda\rangle_{1}B]. \qedhere
  \end{align*}
\end{proof}
\begin{theorem} \label{theorem:fiber-functorial}
  Let $\phi$ be a (semi-)morphism of
  $C^{*}$-$(\fraka,\frakb)$-algebras $\mathcal{
    A}=A^{\alpha,\beta}_{H}$ and $\mathcal{
    C}=C^{\kappa,\lambda}_{L}$, and $\psi$ a (semi-)morphism
  of $C^{*}$-$(\frakbo,\frakc)$-algebras $\mathcal{
    B}=B_{K}^{\gamma,\delta}$ and $\mathcal{
    D}=D^{\mu,\nu}_{M}$. Then there exists a unique
  (semi-)morphism of $C^{*}$-$(\fraka,\frakc)$-algebras
  $\phi \ast \psi$ from ${\cal A} \bfibre {\cal B}$ to
  ${\cal C} \bfibre {\cal D}$ such that
  \begin{align*}
    (\phi \ast
    \psi)(x)R = R x \quad \text{for all } x \in \AfibreB
    \text{ and }R \in I_{M}J_{H} + J_{L}I_{K},
  \end{align*}
  where $I_{X}=\mathcal{
    L}_{(s)}^{\phi}(H_{\beta},L_{\lambda}) \btensor \Id_{X}$
  and $J_{Y}= \Id_{Y} \btensor \mathcal{
    L}_{(s)}^{\psi}(K_{\gamma},M_{\mu})$ for
  $X\in\{K,M\},Y\in\{H,L\}$.
\end{theorem}
\begin{proof}
  By Lemma \ref{lemma:fiber-morphism-ind}, we can define $\phi
  \ast \psi$ to be the restriction of $\rho_{I_{M}} \circ
  \rho_{J_{H}}$ or of $\rho_{J_{L}} \circ \rho_{I_{K}}$ to
  $\AfibreB$.  Uniqueness follows from the fact that
  $[I_{M}J_{H}(\HtensorK)]=[J_{L}I_{K}(\HtensorK)]=L
  \rtensor{\lambda}{\frakb}{\mu} M$.
\end{proof}
\begin{remark}
  Let $A_{H}^{\beta}$, $C^{\lambda}_{L}$ be
  $C^{*}$-$\frakb$-algebras, $B_{K}^{\gamma}$, $D^{\mu}_{M}$
  $C^{*}$-$\frakbo$-algebras, and $\phi \in
  \Mor(A^{\beta}_{H},M(C)^{\lambda}_{L})$, $\psi \in
  \Mor(B^{\gamma}_{K},M(D)^{\mu}_{M})$ such that
  $[\phi(A)C]=C$, $[\psi(B)D]=D$. Then there exists a
  $*$-homomorphism $\phi \bfibre \psi \colon \AfibreB \to
  M(C) \fibre{\lambda}{\frakb}{\mu} M(D) \hookrightarrow
  M(\Cl \bfibre \mD)$, but  in general, we do not know
  whether this is nondegenerate.
\end{remark}

Next, we briefly discuss two kinds of slice maps on fiber
products.  For applications and further details, see
\cite{timmermann:cpmu-hopf}.  The first class of slice maps
arises from a completely positive map on one factor and
takes values in operators on a certain KSGNS-construction,
that is, an internal tensor product with respect to a
completely positive linear map  \cite[\S 4--\S
5]{lance}.
\begin{proposition} \label{proposition:fiber-slice-cp} Let
  $A^{\beta}_{H}$ be a $C^{*}$-$\frakb$-algebra,
  $K_{\gamma}$ a $C^{*}$-$\frakbo$-module, $L$ a Hilbert
  space, $\phi \colon [A + \rho_{\beta}(\frakBo)] \to {\cal
    L}(L)$ a c.p.\ map, and $\theta = \phi \circ
  \rho_{\beta} \colon \frakBo\to {\cal L}(L)$.  Then there
  exists a unique c.p.\ map $\phi \ast \Id \colon
  \Ind_{\kgamma{2}}(A) \to {\cal L}(L {_{\theta}\tl}
  \gamma)$ such that for all $ \zeta,\zeta' \in L, \eta,\eta' \in
  \gamma, x \in \Ind_{\kgamma{2}}(A)$,
  \begin{align} \label{eq:slice-two}
  \langle \zeta \tl \eta|(\phi \ast \Id)(x)(\zeta' \tl
  \eta')\rangle = \langle
  \zeta|\phi(\langle\eta|_{2}x|\eta'\rangle_{2})\zeta'\rangle.
\end{align}
  If $B^{\gamma}_{K}$ is a $C^{*}$-$\frakbo$-algebra, then
  $(\phi \ast \Id)(\AfibreB) \subseteq (\phi(A)' {_{\theta}\tl}
  (B' \cap {\cal L}(K_{\gamma}))' \subseteq {\cal L}(L
  {_{\theta}\tl} \gamma)$.
\end{proposition}
\begin{proof}
  Let $x=(x_{ij})_{i,j} \in M_{n}( \Ind_{\kgamma{2}}(A))$ be
  positive, $\zeta_{1},\ldots,\zeta_{n} \in L$,
  $\eta_{1},\ldots,\eta_{n} \in \gamma$, where $n \in
  \naturals$, and
  $d=\mathrm{diag}(|\eta_{1}\rangle_{2},\ldots,|\eta_{n}\rangle_{2})$.
  Then $0 \leq (\langle
  \eta_{i}|_{2}x_{ij}|\eta_{j}\rangle_{2})_{i,j} = d^{*} x d
  \leq \|x\| d^{*}d $ and hence $0 \leq
  (\phi(\langle\eta_{i}|_{2}x_{ij}|\eta_{j}\rangle_{2}))_{i,j}
  \leq \|x\|  \phi(d^{*}d)$ and
\begin{align*}
  0 \leq \sum_{i,j} \langle \zeta_{i}|\phi(\langle\eta_{i}|_{2} x_{ij} |
  \eta_{j}\rangle_{2}) \zeta_{j}\rangle \leq \|x\| \sum_{i,j} \langle
  \zeta_{i} \tl \eta_{i}| \zeta_{j} \, \tl \eta_{j}\rangle.
\end{align*}
Hence, there exists a map $\phi \ast \Id$ as claimed.  The
verification of the assertion concerning $B^{\gamma}_{K}$ is
straightforward.
\end{proof} 
\begin{remark}
  If $C^{\lambda}_{L}$ is a $C^{*}$-$\frakbo$-algebra and
  $\phi|_{A}$ is a semi-morphism of
  $C^{*}$-$\frakbo$-algebras, then the map $\phi
  \ast \Id$ extends the fiber product $\phi \ast \Id$
  defined in Theorem \ref{theorem:fiber-functorial}.
\end{remark}
Second, we show that the fiber product is functorial with
respect to the following class of maps.  A {\em spatially
  implemented} map of $C^{*}$-$\frakb$-algebras
$A^{\beta}_{H}$ and $C^{\lambda}_{L}$ is a map $\phi \colon
A \to C$ admitting sequences $(S_{n})_{n}$ and $(T_{n})_{n}$
in ${\cal L}(L_{\lambda},H_{\beta})$ such that
\begin{align} \label{eq:slice-spatial} \text{i) } & \
  \sum_{n} S_{n}^{*}S_{n} \text{ and } \sum_{n}
  T_{n}^{*}T_{n} \text{ converge in norm}, & \text{ii) } & \
  \phi(a)=\sum_{n} S_{n}^{*}aT_{n} \text{ for all } a \in A.
\end{align}
Note that condition i) implies norm-convergence of the sum
in ii). Evidently, such a map is linear, extends to a normal
map $\bar \phi \colon A' \to C'$,  its norm is
bounded by $\|\sum_{n}S_{n}^{*}S_{n}\|^{1/2}\|\sum_{n}
T_{n}^{*}T_{n}\|^{1/2}$, and the composition of spatially
implemented maps is spatially implemented again.

\begin{proposition} \label{proposition:fiber-slice-spatial}
  Let $\phi$ be a spatially implemented map of
  $C^{*}$-$\frakb$-algebras  $A^{\beta}_{H}$ and
  $C^{\lambda}_{L}$, and let $B^{\gamma,\delta}_{K}$ be a
  $C^{*}$-$(\frakbo,\frakc)$-algebra. Then there exists a
  spatially implemented map
from $A^{\beta}_{H} \bfibre
  B^{\gamma,\delta}_{K}$ to $C^{\lambda}_{H} \bfibre
  B^{\gamma,\delta}_{K}$ such that $\langle
  \eta|_{2}(\phi \ast
  \Id)(x)|\eta'\rangle_{2}=\phi(\langle\eta|_{2}
  x|\eta'\rangle_{2})$ for all $x \in \AfibreB$,
  $\eta,\eta'\in \gamma$.
\end{proposition}
\begin{proof}
  Uniqueness is clear.  Fix sequences $(S_{n})_{n}$,
  $(T_{n})_{n}$ as in \eqref{eq:slice-spatial} and let
  $\tilde S_{n}:= S_{n} \btensor \Id_{K}$, $\tilde T_{n} :=
  T_{n} \btensor \Id_{K}$ for all $n$. Then $\tilde
  S_{n},\tilde T_{n} \in {\cal L}(L
  \rtensor{\lambda}{\frakb}{\gamma}
  K_{\delta},\HtensorK_{\delta})$ for all $n$, we have
  $\|\sum_{n} \tilde S_{n}^{*}\tilde S_{n}\|=\|\sum_{n}
  S_{n}^{*}S_{n}\|$, $\|\sum_{n} \tilde T_{n}^{*}\tilde
  T_{n}\|=\|\sum_{n} T_{n}^{*}T_{n}\|$, and the map $\phi
  \ast \Id \colon \AfibreB \to {\cal L}(L
  \rtensor{\lambda}{\frakb}{\gamma} K)$ given by $x \mapsto
  \sum_{n} \tilde T_{n}^{*}x \tilde S_{n}$ has the desired
  properties.  Indeed, let $x \in \AfibreB$, $\eta,\eta' \in
  \gamma$. Then $\tilde S_{n}|\eta\rangle_{2} =
  |\eta\rangle_{2} S_{n}$ and $\tilde
  T_{n}|\eta'\rangle_{2}=|\eta'\rangle_{2}T_{n}$ for all
  $n$, and hence $ \langle \eta|_{2}(\phi \ast
  \Id)(x)|\eta'\rangle_{2} = \phi(\langle
  \eta|_{2}x|\eta'\rangle_{2})$. It remains to show that
  $(\phi \ast \Id)(x) \in C \fibre{\lambda}{\frakb}{\gamma}
  B$. Consider the expression $(\phi \ast
  \Id)(x)|\eta'\rangle_{2} = \sum_{n} \tilde S_{n}^{*}
  x|\eta'\rangle_{2} T_{n}$.  This sum converges in norm and
  each summand lies in $[\kgamma{2}{\cal L}(H)]$ because
  $x|\eta'\rangle_{2} \in [\kgamma{2}A]$ and $[\tilde
  S^{*}_{n}\kgamma{2}]=[\kgamma{2}S_{n}^{*}]$. Since
  $\langle\eta''|_{2}(\phi \ast \Id)(x)|\eta'\rangle_{2} \in
  C$ for each $\eta'' \in \gamma$, we
  can conclude that the sum lies in $[\kgamma{2}C]$.
  Finally, consider the expression $(\phi \ast
  \Id)(x)|\xi\rangle_{1} = \sum_{n} \tilde S_{n} x \tilde
  T_{n}|\xi\rangle_{1}$, where $\xi \in \lambda$. Again, the
  sum converges in norm and each summand lies in
  $[|\lambda\rangle_{1}B]$ because $\tilde S_{n}^{*}x\tilde
  T_{n}|\xi\rangle_{1} =\tilde S_{n}^{*}x |T_{n}\xi\rangle_{1}
  \in \tilde S_{n}^{*}(\AfibreB)\kbeta{1} \subseteq [\tilde
  S_{n}^{*}\kbeta{1}B] \subseteq [|\lambda\rangle_{1}B]$.
\end{proof}
\begin{remarks}
  \begin{enumerate}
  \item The map $\phi \ast \Id$ constructed above is a
    ``slice map'' in the case where $C^{\lambda}_{L} = {\cal
      L}(\frakK)_{\frakK}^{\frakB}$ and $S_{n},T_{n} \in
    \beta \subseteq {\cal L}(\frakK_{\frakB},H_{\beta})$ for
    all $n$. Then, we can identify $C
    \fibre{\lambda}{\frakb}{\gamma} B$ with a
    $C^{*}$-subalgebra of ${\cal L}(K)$, and $\phi \ast
    \Id$ is just the map $\AfibreB \to B$ given by $x
    \mapsto \sum_{n} \langle S_{n}|_{1}X|T_{n}\rangle_{1}$.
  \item Assume that the extension $\tilde \phi \colon
    [A+\rho_{\beta}(\frakBo)] \to C$ given by $x \mapsto
    \sum_{n} S_{n}^{*}xT_{n}$ is completely positive. Here,
    we use the notation of the proof above. Then the map
    $\tilde \phi \ast \Id$ constructed in Proposition
    \ref{proposition:fiber-slice-cp} extends the map $\phi \ast
    \Id$ of Proposition
    \ref{proposition:fiber-slice-spatial} because then
    $\theta = \rho_{\lambda}$ and
    hence $\langle \eta|_{2} (\tilde \phi \ast
    \Id)(x)|\eta'\rangle_{2} = \tilde \phi(\langle
    \eta|_{2}x|\eta'\rangle_{2})$ for all $x \in \AfibreB$
    and $\eta,\eta' \in \gamma$.
  \end{enumerate}
\end{remarks}
 Of course, slice maps of the form $\Id \ast \phi$ can
  be constructed in a similar way.

\subsection{Further categorical properties}
\label{subsection:fiber-categorical}

The fiber product of $C^{*}$-algebras is neither
associative, unital, nor compatible with infinite sums.

\paragraph{Non-associativity}
Let ${\cal A}=A^{\alpha,\beta}_{H}$ be a
$C^{*}$-$(\frakao,\frakb)$-algebra, ${\cal
  B}=B^{\gamma,\delta}_{K}$ a
$C^{*}$-$(\frakbo,\frakc)$-algebra, and ${\cal
  C}=C^{\epsilon,\phi}_{L}$ a
$C^{*}$-$(\frakco,\frakd)$-algebra. Then we can form the
fiber products $({\cal A} \bfibre {\cal B}) \cfibre{\cal C}$
and ${\cal A} \bfibre ({\cal B} \cfibre {\cal C})$. The
following example shows that these $C^{*}$-algebras need
{\em not} be identified by the canonical isomorphism
$a_{\fraka,\frakb,\frakc,\frakd}(\eLf,\cKd,\aHb)$ of
Proposition \ref{proposition:rtp-properties}.  A similar
phenomenon occurs in the purely algebraic setting with the
Takeuchi $\times_{R}$-product \cite{takeuchi}.
\begin{example} \label{example:fp-nonassociative} Let
  $\fraka=\frakb= \frakc=\frakd$ be the trivial
  $C^{*}$-base, $H=l^{2}(\naturals)$, $\alpha={\cal
    L}(\complex,H)$,  ${\cal A}={\cal B}={\cal C} = {\cal
    L}(H)^{\alpha,\alpha}_{H}$. Identify $H
  \rtensor{\alpha}{\frakb}{\alpha} K
  \rtensor{\alpha}{\frakc}{\alpha} L\cong \alpha \otimes H
  \otimes \alpha$ with $H \otimes H \otimes H$ via
  $|\xi\rangle \tr \zeta \tl |\eta\rangle \equiv \xi
  \otimes \zeta \otimes \eta$, fix an orthonormal basis
  $(e_{n})_{n \in \naturals}$ of $H$, and define $T \in
  {\cal L}(H^{\otimes 3})$ by
  \begin{align*}
    T(e_{k} \otimes e_{l} \otimes e_{m}) =
    \begin{cases}
      e_{k} \otimes e_{l} \otimes e_{m} &
      \text{for all } k,l,m \in \naturals \text{ s.t. } m \leq k
      +l, \\
      e_{l} \otimes e_{k} \otimes e_{m} &
      \text{for all } k,l,m \in \naturals \text{ s.t. } m > k+l.
    \end{cases} 
  \end{align*}
  We show that $T$ belongs to the underlying $C^{*}$-algebra
  of $ ({\cal A} \bfibre {\cal B}) \cfibre {\cal C}$, but not
  of ${\cal A} \bfibre ({\cal B} \cfibre {\cal C})$.

  For each $\xi \in H$ and $\omega \in H^{\otimes 2}$,
  define $|\xi\rangle_{1},|\xi\rangle_{3} \in {\cal
    L}(H^{\otimes 2}, H^{\otimes 3})$ and
  $|\omega\rangle_{12} \in {\cal L}(H,H^{\otimes 3})$ by
  $\upsilon \mapsto \xi \otimes \upsilon$, $\upsilon \mapsto
  \upsilon \otimes \xi$, and $\zeta \mapsto \omega \otimes
  \zeta$, respectively.  Then for all $k,l,m \in \naturals$,
  \begin{align*}
    T|e_{k} \otimes e_{l}\rangle_{12} &=
    |e_{k} \otimes e_{l}\rangle_{12} P_{l+k} +
    |e_{l} \otimes e_{k}\rangle_{12} (\Id -
    P_{l+k}), \text{ where } P_{l+k} := \sum_{m \leq k+l}
    |e_{m}\rangle\langle e_{m}|, \\
    T|e_{m}\rangle_{3} &=|e_{m}\rangle_{3} (\Id
    +
    \Sigma_{m}), \ \text{ where } \Sigma_{m} := \sum_{\substack{k,l \\
        k+l < m}} |e_{l}\otimes e_{k} -
    e_{k}\otimes e_{l}\rangle\langle
    e_{k}\otimes e_{l}|,
  \end{align*}
and therefore,
\begin{align*}
 T|H^{\otimes 2}\rangle_{12} &\in [|H^{\otimes
    2}\rangle_{12}{\cal L}(H)], & T|\alpha\rangle_{3} &\in
  [|\alpha_{3}\rangle(\Id + {\cal K}(H) \otimes {\cal
    K}(H))] \subseteq [|\alpha\rangle_{3}({\cal L}(H)
  \fibre{\alpha}{\frakb}{\alpha} {\cal L}(H))].
\end{align*}
Since $T=T^{*}$, we can conclude that $T$ belongs to $({\cal
  L}(H) \fibre{\alpha}{\frakb}{\alpha} {\cal L}(H)_{\alpha})
\fibre{}{\frakb}{\alpha} {\cal L}(H)$.  However,
  \begin{align*}
    T|e_{0}\rangle_{1} =|e_{0}\rangle_{1} Q +
    \sum_{l} |e_{l}\rangle_{1} Q_{l}, \ \text{ where
    } Q &= \sum_{m \leq l} |e_{l} \otimes
    e_{m}\rangle\langle e_{l} \otimes
    e_{m}| \\ \text{ and } Q_{l} &= \sum_{m > l}
    |e_{0} \otimes e_{m}\rangle\langle
    e_{l} \otimes e_{m}|,
  \end{align*}
  and $|e_{0}\rangle_{1}Q \in [|\alpha\rangle_{1}
  {\cal L}(H \otimes H)]$,  but $\sum_{l}
  |e_{l}\rangle_{1} Q_{l} \not\in  [|\alpha\rangle_{1}
  {\cal L}(H \otimes H)]$ because the sum
  \begin{align*}
    \sum_{l} Q_{l}^{*}Q_{l} = \sum_{l} \sum_{m > l}
    |e_{l} \otimes e_{m}\rangle\langle
    e_{l} \otimes e_{m}|
  \end{align*}
   does not converge in norm. Hence,
  $T|e_{0}\rangle_{1} \not \in [|\alpha\rangle_{1}{\cal
    L}(H \otimes H)]$ and  $T \not\in {\cal L}(H)
  \fibre{\alpha}{\frakb}{} (_{\alpha}{\cal L}(H)
  \fibre{\alpha}{\frakb}{\alpha} {\cal L}(H)) $.
  \end{example}

\paragraph{Unitality}
A unit for the fiber product relative to $\frakb$ would be
a $C^{*}$-$(\frakbo,\frakb)$-algebra ${\cal
  U}=\frakU^{\frakBo,\frakB}_{\frakK}$ such that for all
$C^{*}$-$(\frakao,\frakb)$-algebras ${\cal
  A}=A^{\alpha,\beta}_{H}$ and all
$C^{*}$-$(\frakbo,\frakc)$-algebras ${\cal
  B}=B^{\gamma,\delta}_{K}$, we have ${\cal A}=\Ad_{r}({\cal A}
\bfibre {\cal U})$ and ${\cal B}=\Ad_{l}({\cal U}\bfibre
{\cal B})$, where $r=r_{\fraka,\frakb}(\aHb)$ and
$l=l_{\frakb,\frakc}(\cKd)$ (see Proposition
\ref{proposition:rtp-properties}). The relations $r
\kbeta{1}=\beta$,
$r|\frakBo\rangle_{2}=\rho_{\beta}(\frakBo)$,
$l\kgamma{2}=\gamma$,
$l|\frakB\rangle_{1}=\rho_{\gamma}(\frakB)$ imply
  \begin{align} \label{eq:fp-unital}
    \begin{aligned}
      \Ad_{r}(\Ab \bfibre {_{\frakBo} \frakU}) &=
      \Ind_{\beta}(\frakU) \cap
      \Ind_{\rho_{\beta}(\frakBo)}(A), &
      \Ad_{l}(\frakU_{\frakB} \bfibre \gB) &=
      \Ind_{\rho_{\gamma}(\frakB)}(B) \cap
      \Ind_{\gamma}(\frakU).
    \end{aligned}
  \end{align}
 If $\frakBo$ and $\frakB$ are unital, then
  $\Ind_{\rho_{\beta}(\frakBo)}(A)=A$ and
  $\Ind_{\rho_{\gamma}(\frakB)}(B)=B$, and then the
  $C^{*}$-$(\frakbo,\frakb)$-algebra ${\cal
    L}(\frakK)^{\frakBo,\frakB}_{\frakK}$ is a unit for the
  fiber product on the full subcategories of all
  $A^{\alpha,\beta}_{H}$ and $B^{\gamma,\delta}_{K}$
  satisfying $A \subseteq \Ind_{\beta}({\cal L}(\frakK))$
  and $B \subseteq \Ind_{\gamma}({\cal L}(\frakK))$.  
\begin{remarks}
  \begin{enumerate}
  \item If $A \subseteq \Ind_{\alpha}({\cal L}(\frakH))$ and
  $B \subseteq \Ind_{\gamma}({\cal L}(\frakL))$, then 
  $\AfibreB \subseteq \Ind_{(\alpha \lt \gamma)}({\cal
    L}(\frakH)) \cap \Ind_{(\beta \rt \delta)}({\cal
    L}(\frakK))$.
  \item $\Ind_{\beta}(\frakBo)={\cal L}(H_{\beta})$, and if
    $\frakBo$ is unital, then $\Ad_{r}(\Ab \bfibre
    {_{\frakBo}} \frakBo) = A \cap {\cal L}(H_{\beta}) =
    A^{(\beta)}$.
    \item $\Ad_{r}(\frakB \fibre{\frakB}{\frakb}{\frakBo}
      \frakBo) = {\cal L}(\frakK_{\frakB}) \cap {\cal
        L}(\frakK_{\frakBo}) = M(\frakB) \cap M(\frakBo)$.
    \end{enumerate}
\end{remarks}

\paragraph{Compatibility with sums and products}
The fiber product is compatible with finite sums in the
following sense. Let $({\cal A}^{i})_{i}$ be a finite family
of $C^{*}$-$(\frakao,\frakb)$-algebras and $({\cal
  B}^{j})_{j}$ a finite family of
$C^{*}$-$(\frakbo,\frakc)$-algebras. For each $i,j$, denote
by $\iota^{i}_{{\cal A}} \colon {\cal A}^{i} \to
\boxplus_{i'} {\cal A}^{i'} $, $\iota^{j}_{{\cal B}} \colon
{\cal B}^{j} \to \boxplus_{j'} {\cal B}^{j'} $ and
$\pi^{i}_{{\cal A}} \colon \boxplus_{i'} {\cal A}^{i'} \to
{\cal A}^{i}$, $\pi^{j}_{{\cal B}} \colon \boxplus_{j'}
{\cal B}^{j'} \to {\cal B}^{j}$ the canonical inclusions and
projections, respectively.  One easily verifies that 
there exist inverse isomorphisms $\boxplus_{i,j} {\cal
  A}^{i} \bfibre {\cal B}^{j} \leftrightarrows (\boxplus_{i}
{\cal A}^{i}) \bfibre (\boxplus_{j} {\cal B}^{j})$, given by
$(x_{i,j})_{i,j} \mapsto \sum_{i,j} (\iota^{i}_{{\cal A}}
\bfibre \iota^{j}_{{\cal B}})(x_{i,j})$ and
$\big((\pi^{i}_{{\cal A}} \bfibre \pi^{j}_{{\cal
    B}})(y)\big)_{i,j} \mapsfrom y$, respectively.  However,
the fiber product is neither compatible with infinite sums
nor infinite products:
\begin{examples}
  Let $\trivbase=(\complex,\complex,\complex)$ be the
  trivial $C^{*}$-base.
\begin{enumerate}
\item For each $i,j \in \naturals$, let ${\cal A}^{i}$ and
  ${\cal B}^{j}$ be the $C^{*}$-$\trivbase$-algebra
  $\complex^{\complex}_{\complex}$. Identify the Hilbert
  space $\bigoplus_{i,j} \complex
  \rtensor{\complex}{\trivbase}{\complex} \complex$ with
  $l^{2}(\naturals \times \naturals)$ in the canonical
  way. Then $\bigoplus_{i,j} {\cal A}^{i} \sfibre{\trivbase}
  {\cal B}^{j}$ corresponds to $C_{0}(\naturals \times
  \naturals)$, represented on $l^{2}(\naturals \times
  \naturals)$ by multiplication operators, but
  $(\bigoplus_{i} {\cal A}^{i}) \sfibre{\trivbase}
  (\bigoplus_{j} {\cal B}^{j}) \cong C_{0}(\naturals)
  \sfibre{\trivbase} C_{0}(\naturals)$ is strictly larger
  and contains, for example, the characteristic function of
  the diagonal $\{(x,x) \mid x \in \naturals\}$ (see Example
  \ref{examples:fp-commutative}).
\item Let $H=l^{2}(\naturals)$, $\alpha={\cal
    L}(\complex,H)$, and let ${\cal A}$ and ${\cal B}^{j}$
  be the $C^{*}$-$\trivbase$-algebra ${\cal
    K}(H)_{H}^{\alpha}$ for all $j$. Identify $H
  \rtensor{\alpha}{\trivbase}{\alpha} H$ with $H \otimes H$
  as in Example \ref{examples:rtp-rtp} i),  choose an
  orthonormal basis $(e_{k})_{k \in \naturals}$ of $H$, and
  put $y_{j}:=|e_{j} \otimes e_{0}\rangle\langle e_{0}
  \otimes e_{0}| \in {\cal K}(H \otimes H)$ for each $j \in
  \naturals$.  Then $y:=(y_{j})_{j} \in \prod_{j} {\cal A}
  \sfibre{\trivbase} {\cal B}^{j}$ because $y_{j} \in {\cal
    K}(H) \otimes {\cal K}(H) \subset {\cal A}
  \sfibre{\trivbase} {\cal B}^{j}$ for all $j \in
  \naturals$, but with respect to the canonical
  identification $\bigoplus_{j} H \otimes H \cong H \otimes
  \left(\bigoplus_{j} \otimes H\right)$, we have $y \not\in
  {\cal A} \sfibre{\trivbase} (\prod_{j} {\cal B}^{j})$
  because $y |e_{0}\rangle_{1}$ corresponds to the family
  $(|e_{j}\rangle_{1}|e_{0}\rangle\langle e_{0}|)_{j} \in
  \prod_{j} {\cal L}(H, H \otimes H ) \subseteq {\cal
    L}(\bigoplus_{j} H, \bigoplus_{j} H \otimes H )$ which
  is not contained in the space $[|\alpha\rangle_{1}{\cal
    L}(\bigoplus_{j} H)]$.
\end{enumerate}
  \end{examples}

  \subsection{A fiber product of non-represented
    $C^{*}$-algebras}

\label{subsection:fiber-universal}

The spatial fiber product of $C^{*}$-algebras represented
on $C^{*}$-modules yields a fiber product of
non-represented $C^{*}$-algebras as follows.  

Let $\frakb=\cbasesb$ be a $C^{*}$-base. In Subsection
\ref{subsection:fiber-algebras}, we constructed a functor
$\bfR_{\frakb} \colon\calgr{\frakBo} \to \calgs{\frakb}$
that associates to each $C^{*}$-$\frakBo$-algebra a
universal representation in form of a
$C^{*}$-$\frakb$-algebra. Replacing $\frakb$ by $\frakbo$,
we  obtain a functor $\bfR_{\frakbo} \colon
\calgr{\frakB} \to \calgs{\frakb}$, and composition of these
with the spatial fiber product gives a fiber product of
non-represented $C^{*}$-algebras in form of a functor
\begin{align*}
  \calgr{\frakBo} \times \calgr{\frakB}
  \xrightarrow{\bfR_{\frakb} \times \bfR_{\frakbo}}
  \calgs{\frakb} \times \calgs{\frakbo} \to \calg, \quad
  (C_{\sigma},D_{\tau}) \mapsto \bfR_{\frakb}(C_{\sigma})
  \bfibre \bfR_{\frakbo}(D_{\tau}),
\end{align*}
where $\bfcs$ denotes the category of $C^{*}$-algebras and
$*$-homomorphisms.  In categorical terms, this is the right
Kan extension of the spatial fiber product on
$\calgs{\frakb} \times \calgs{\frakbo}$ along the product of
the forgetful functors $\bfU_{\frakb} \times \bfU_{\frakbo}
\colon \calgs{\frakb} \times \calgs{\frakbo} \to
\calgr{\frakBo} \times \calgr{\frakB}$ \cite[\S X]{maclane}.

Given further $C^{*}$-bases $\fraka=\cbasesa$ and
$\frakc=\cbasesc$, we similarly obtain a functor
\begin{align*}
  \calgr{(\frakA,\frakBo)} \times \calgr{(\frakB,\frakCo)}
  \xrightarrow{\bfR_{(\frakao,\frakb)} \times
    \bfR_{(\frakbo,\frakc)}} \calgs{(\frakao,\frakb)} \times
  \calgs{(\frakbo,\frakc)} \to \calgs{(\frakao,\frakc)}
  \xrightarrow{\bfU_{(\frakao,\frakc)}}
  \calgr{(\frakA,\frakCo)},
\end{align*}
and, using Remark \ref{remark:fiber-adjoints}, a natural
transformation between the compositions in the square
\begin{align*}
  \xymatrix@R=15pt{ {\calgr{(\frakA,\frakBo)} \times
      \calgr{(\frakB,\frakCo)}} \ar[r] \ar[d]
    & {\calgr{(\frakA,\frakCo)}}   \ar[d]   \ar@{=>}[ld]  \\
    { \calgr{\frakBo} \times \calgr{\frakB}} \ar[r] &
    {\bfcs,}},
\end{align*}
where the vertical maps are the forgetful functors.

\section{Relation to the setting of von Neumann algebras}

\label{section:vn}

Throughout this section, let $N$ be a von Neumann algebra
with a n.s.f.\ weight $\mu$, denote by
$\frakN_{\mu},H_{\mu},\pi_{\mu},J_{\mu}$  the usual
objects of Tomita-Takesaki theory \cite{takesaki:2}, and
define the antirepresentation $\pi_{\mu}^{op} \colon N  \to
{\cal L}(H_{\mu})$ by $x \mapsto J_{\mu}\pi_{\mu}(x^{*})J_{\mu}$.

\subsection{Adaptation to von Neumann algebras}

\label{subsection:vn-adapt}

The definitions and constructions presented in Sections
\ref{section:rtp} and \ref{section:fiber} can be adapted to
a variety of other settings.  We now briefly explain what
happens when we pass to the setting of von Neumann algebras.
Instead of a $C^{*}$-base, we start with the triple
$\frakb=(\frakK,\frakB,\frakBo)$, where $\frakK=H_{\mu}$,
$\frakB=\pi_{\mu}(N)$, and $\frakBo=J_{\mu}
\pi_{\mu}(N)J_{\mu}$.  Next, we define
$W^{*}$-$\frakb$-modules,
$W^{*}$-$(\frakbo,\frakb)$-modules, their relative tensor
product, $W^{*}$-$\frakb$-algebras, and the fiber product by
just replacing the norm closure $[\, \cdot\, ]$ by the
closure with respect to the weak operator topology $[\,
\cdot\, ]_{w}$ everywhere in Sections \ref{section:rtp} and
\ref{section:fiber}. We then recover Connes' fusion of
Hilbert bimodules over $N$ and Sauvageot's fiber product:
\begin{description}
\item [Modules] Let $H$ be some Hilbert space.  If
  $(H,\rho)$ is a right $N$-module, then the space
  \begin{align*}
   \alpha={\cal L}((\frakK,\pi_{\mu}^{op}), (H,\rho)) := \{ T
  \in {\cal L}(\frakK,H) : T\pi_{\mu}^{op}(x)=\rho(x)T \text{ for
  all  } x \in N\} 
  \end{align*}
  satisfies $[\alpha \frakK]=H,
  [\alpha^{*}\alpha]_{w}=\frakB, \alpha\frakB \subseteq
  \alpha$, and $\rho_{\alpha} \circ \pi_{\mu}^{op}$ (see
  Lemma \ref{lemma:rtp-morphism}) coincides with
  $\rho$. Conversely, if $\alpha \subseteq {\cal
    L}(\frakK,H)$ is a weakly closed subspace satisfying the
  three preceding equations, then $(H,\rho_{\alpha} \circ
  \pi_{\mu}^{op})$ is a right $N$-module and $\alpha = {\cal
    L}((\frakK,\pi_{\mu}^{op}),(H,\rho_{\alpha}\circ\pi_{\mu}^{op}))$
  \cite{skeide}.  We thus obtain a bijective correspondence
  between right $N$-modules and
  $W^{*}$-$\frakb$-modules. This correspondence is an
  isomorphism of categories since for every other right
  $N$-module $(K,\sigma)$, an operator $T \in {\cal L}(H,K)$
  intertwines $\rho$ and $\sigma$ if and only if $T\alpha$
  is contained in $ \beta:= {\cal
    L}((\frakK,\pi_{\mu}^{op}),(K,\sigma))$.  For
  $W^{*}$-$\frakb$-modules, the notions of morphisms and
  semi-morphisms coincide.
\item [Algebras] Let $H,\rho,\alpha$ be as above and let $A
  \subseteq {\cal L}(H)$ be a von Neumann algebra. Then
  $\rho(N)\subseteq A$ if and only if
  $\rho_{\alpha}(\frakB)A \subseteq A$.  Thus,
  $W^{*}$-$\frakb$-algebras correspond with von Neumann
  algebras equipped with a normal unital embedding of
  $N$. Moreover, let $K,\sigma,\beta$ be as above, let $B
  \subseteq {\cal L}(K)$ be a von Neumann algebra, assume
  $\rho(N)\subseteq A$ and $\sigma(N) \subseteq B$, and let
  $\pi \colon A \to B$ be a $*$-homomorphism satisfying $\pi
  \circ \rho = \sigma$. Then $\pi$ is normal if and only if
  $[{\cal
    L}^{\pi}(H_{\alpha},K_{\beta})\alpha]_{w}=\beta$. Indeed,
  the ``if'' part is straightforward (see Lemma
  \ref{lemma:fiber-morphism}), and the ``only if'' part
  follows easily from the fact that every normal
  $*$-homo\-morphism is the composition of an amplification,
  reduction, and unitary transformation \cite[\S
  4.4]{dixmier:2}.
\item [Bimodules] Let $(H,\rho)$ be a left $N$-module,
  $(H,\sigma)$ a right $N$-module, $\alpha={\cal
    L}((\frakK,\pi_{\mu}),(H,\rho))$ and $\beta={\cal
    L}((\frakK,\pi_{\mu}^{op}),(H,\sigma))$. Then
  $(H,\rho,\sigma)$ is an $N$-bimodule if and only if
  $\rho(N) \beta = \beta$ and $\sigma(N)\alpha=\alpha$, and
  thus we obtain an isomorphism between the category of
  $N$-bimodules and the category of
  $W^{*}$-$(\frakbo,\frakb)$-modules.
\item [Fusion] The preceding considerations and formula
  \eqref{eq:rtp-connes-algebraic} show that the relative
  tensor product of $W^{*}$-$(\frakbo,\frakb)$-modules
  corresponds to Connes' fusion of $N$-bimodules.
\item [Fiber product] Let $(H,\rho)$ be a right $N$-module,
  $(K,\sigma)$ a left $N$-module, $\alpha={\cal
    L}((\frakK,\pi_{\mu}^{op}),(H,\rho))$, $\beta={\cal
    L}((\frakK,\pi_{\mu}),(K,\sigma))$, and let $A \subseteq
  {\cal L}(H)$ and $B \subseteq {\cal L}(K)$ be von Neumann
  algebras satisfying $\rho(N) \subseteq H$ and $\sigma(N)
  \subseteq K$. One easily verifies the equivalence of the
  following conditions for each $x \in {\cal L}(\HtensorK)$:
  i) $x \kalpha{1} \subseteq [\kalpha{1}B]_{w}$, ii)
  $\balpha{1}x\kalpha{1} \subseteq B$, iii) $x \in (\Id_{H}
  \btensor B')'$.  Consequently, the fiber product of $A$
  and $B$, considered as a $W^{*}$-$\frakb$-algebra and a
  $W^{*}$-$\frakbo$-algebra, coincides with the fiber
  product $(\Id_{H} \btensor B')' \cap (A' \btensor
  \Id_{K})' = (A' \btensor B')'$ of Sauvageot.
\end{description}

\subsection{Relation to Connes' fusion and Sauvageot's fiber
  product}

\label{subsection:vn-relation}

Let $\frakb=\cbasesb$ be a $C^{*}$-base such that
$\frakK=H_{\mu}$, $\frakB'' = \pi_{\mu}(N)$, $(\frakBo)'' =
\pi_{\mu}^{op}(N) = \frakB'$.

Denote by $\cbimodr{\frakb}$ the category of all
$C^{*}$-$(\frakbo,\frakb)$-modules with all semi-morphisms,
and by $\wbimodr{N}$ the category of all $N$-bimodules,
respectively.  Lemmas \ref{lemma:rtp-morphism} and
\ref{lemma:rtp-modules} imply:
\begin{proposition}
  There exists a faithful functor $\bfF \colon
  \cbimodr{\frakb} \to \wbimodr{N}$, given by
  ${_{\alpha}H_{\beta}} \mapsto (H,\rho_{\alpha} \circ
  \pi_{\mu}, \rho_{\beta} \circ \pi_{\mu}^{op})$ on objects
  and $ T \mapsto T$ on morphisms. \qed
\end{proposition}
The categories $ \cbimodr{\frakb}$ and $\wbimodr{N}$ carry
the structure of a monoidal category \cite{maclane}, and we now
show that the functor $\bfF$ above is monoidal.
Let $H_{\beta}$ be a $C^{*}$-$\frakb$-module, $K_{\gamma}$ a
$C^{*}$-$\frakbo$-module, and let
\begin{align*}
  \rho&=\rho_{\beta} \circ \pi_{\mu}^{op}, & X &= {\cal
    L}((\frakK,\pi_{\mu}^{op}),(H,\rho)), & \sigma
  &=\rho_{\gamma} \circ \pi_{\mu}, & Y&= {\cal
    L}((\frakK,\pi_{\mu}),(K,\sigma)).
\end{align*}
Given subspaces $X_{0} \subseteq X$ and $Y_{0} \subseteq Y$,
we define a sesquilinear form $\langle \frei|\frei\rangle$
on the algebraic tensor product $X_{0} \odot \frakK \odot
Y_{0}$ such that for all $\xi,\xi' \in X_{0}, \zeta,\zeta' \in \frakK,
\eta,\eta' \in Y_{0}$,
\begin{align*}
  \langle \xi \odot \zeta \odot \eta|\xi' \odot \zeta' \odot
  \eta'\rangle = \langle
  \zeta|(\xi^{*}\xi')(\eta^{*}\eta')\eta'\rangle = \langle
  \zeta|(\eta^{*}\eta')(\xi^{*}\xi')\eta'\rangle 
\end{align*}
Denote by $X_{0} \tr \frakK \tl Y_{0}$ the Hilbert space
obtained by forming the separated completion.
\begin{lemma} \label{lemma:vn-rtp-embedding} Let $X_{0}
  \subseteq X$ and $Y_{0} \subseteq Y$ be subspaces
  satisfying $[X_{0}\frakK]=H$ and $[Y_{0}\frakK]=K$. Then
  the natural map $X_{0} \tr \frakK \tl Y_{0} \to X \tr
  \frakK \tl Y$ is an isomorphism.
\end{lemma}
\begin{proof}
  Injectivity is clear. The natural map $X_{0} \tr \frakK \tl Y_{0}
  \to X \tr \frakK \tl Y_{0}$ is surjective because both
  spaces coincide with the separated completion of the
  algebraic tensor product
  $H \odot Y_{0}$ with respect to the sesquilinear inner
  form given by $\langle \omega \odot \eta|\omega' \odot
  \eta'\rangle = \langle \omega|
  \rho_{\beta}(\eta^{*}\eta')\omega'\rangle$, and a
  similar argument shows that the natural map
  $X \tr \frakK \tl Y_{0} \to X \tr \frakK \tl Y$ is surjective.
\end{proof}
We conclude that Connes' original definition of the relative
tensor product $H \rtensor{\rho}{\mu}{\sigma} K$ via bounded
vectors coincides with the algebraic one given in
\eqref{eq:rtp-connes-algebraic} and with the relative tensor
product $\HtensorK$.
\begin{theorem} \label{theorem:vn-fusion}
  There exists a natural isomorphism between the
  compositions in the square
  \begin{align*}
    \xymatrix@C=35pt@R=15pt{ {\cbimodr{\frakb} \times
        \cbimodr{\frakb}} \ar[r]^(0.62){-\btensor -} \ar[d]_{\bfF
        \times \bfF} &
      {\cbimodr{\frakb}} \ar@{=>}[ld] \ar[d]^{\bfF} \\
      {\wbimodr{N} \times \wbimodr{N}}
      \ar[r]_(0.62){-\rtensor{}{\mu}{}-} &{\wbimodr{N},} }
  \end{align*}
  given for each object $(\aHb,\cKd) \in \cbimodr{\frakb}
  \times \cbimodr{\frakb}$ by the natural map
  \begin{align} \label{eq:vn-rtp-embedding}
    \HtensorK = \beta \tr \frakK \tl \gamma \to X \tr \frakK
    \tl Y = H \rtensor{\rho}{\mu}{\sigma} K.
  \end{align}
 With respect to this isomorphism, the
  functor $\bfF \colon \cbimodr{\frakb} \to \wbimodr{N}$ is
  monoidal.
\end{theorem}
\begin{proof}
  Lemma \ref{lemma:vn-rtp-embedding} implies that the map
  \eqref{eq:vn-rtp-embedding} is an isomorphism. Evidently,
  this map is natural with respect to $\aHb$ and $\cKd$. The
  verification of the assertion concerning $\bfF$ is now
  tedious but straightforward.
\end{proof}
Denote by $\bfcsnd_{(\frakbo,\frakb)}$ the category of all
$C^{*}$-$(\frakbo,\frakb)$-algebras $A^{\alpha,\beta}_{H}$
satisfying $\rho_{\alpha}(\frakB)+\rho_{\beta}(\frakBo)
\subseteq A$ together with all semi-morphisms, and by
$\bfws_{(N,N^{op})}$ the category of all von Neumann
algebras $A$ equipped with a normal, unital embedding and
anti-embedding $\iota^{(op)}_{A} \colon N \to A$ such
that $[\iota_{A}(N),\iota_{A}^{op}(N)]=0$, together with all
morphisms preserving these (anti-)embeddings. Lemma
\ref{lemma:fiber-morphism} implies:
\begin{proposition} 
  There exists a faithful functor $\bfG \colon
  \bfcsnd_{(\frakbo,\frakb)} \to \bfws_{(N,N^{op})}$, given
  by $(\aHb,A) \mapsto (A'',\rho_{\alpha} \circ \pi_{\mu},
  \rho_{\beta} \circ \pi_{\mu}^{op})$ on objects and $\phi
  \mapsto \phi''$ on morphisms, where $\phi''$ denotes the
  normal extension of $\phi$. \qed
\end{proposition}
By Lemma \ref{lemma:fiber-properties}, ${\cal A} \bfibre
{\cal B} \in \bfcsnd_{(\frakbo,\frakb)}$ for all ${\cal
  A},{\cal B} \in \bfcsnd_{(\frakbo,\frakb)}$, but
$ \bfcsnd_{(\frakbo,\frakb)}$ is not a monoidal category
with respect to the fiber product because the latter is not
associative (see Subsection \ref{subsection:fiber-categorical}).
\begin{proposition} \label{proposition:vn-fiber} There exists a natural transformation
  \begin{align*}
    \xymatrix@C=35pt@R=15pt{ {\bfcsnd_{(\frakbo,\frakb)} \times
        \bfcsnd_{(\frakbo,\frakb)}} \ar[r]^(0.6){-\fibre{}{\frakb}{}-} \ar[d]_{\bfG\times\bfG} &
      {\bfcsnd_{(\frakbo,\frakb)}} \ar[d]^{\bfG} \ar@{=>}[ld] \\
      {\bfws_{(N,N^{op})} \times \bfws_{(N,N^{op})}} \ar[r]_(0.6){-\fibre{}{\mu}{}-}
      & {\bfws_{(N,N^{op})}},  }
  \end{align*}
  given for each object $A^{\alpha,\beta}_{H}$ and
  $B^{\gamma,\delta}_{K}$ by conjugation with the
  isomorphism \eqref{eq:vn-rtp-embedding}.
\end{proposition}
\begin{proof}
  Immediate from Theorem \ref{theorem:vn-fusion} and Lemma
  \ref{lemma:fiber-properties}.
\end{proof}

\subsection{A categorical interpretation of the fiber
  product of von Neumann algebras}

\label{subsection:vn-category}

We keep the notation introduced above, denote by $\hilb$ the
category of Hilbert spaces and bounded linear operators, and
call  a subcategory of $\wmodrn$  a {\em $*$-subcategory}
if it is closed with respect to the involution $T \mapsto
T^{*}$ of morphisms. 
\begin{definition} \label{definition:fiber-cat-catn} A {\em
    category over $\wmodrn$} is a category $\bfC$ equipped
  with a functor $\bfU_{\bfC} \colon \bfC \to
  \wmodr{(N,N^{op})}$ such that $\bfU_{\bfC} \bfC$ is a
  $*$-subcategory of $\wmodrn$.  Let $(\bfC,\bfU_{\bfC})$ be
  such a category. We loosely refer to $\bfC$ as a category
  over $\wmodrn$ without mentioning $\bfU_{\bfC}$
  explicitly, and denote by $\bfH_{\bfC}$ the composition of
  $\bfU_{\bfC}$ with the forgetful functor $\wmodrn \to
  \hilb$.  We call an object $G \in \bfC$ {\em separating}
  if $[\bfH_{\bfC}\bfC(G,X)(\bfH_{\bfC} G)]=\bfH_{\bfC} X$
  for each $X \in \bfC$.

  We denote by $\bfCatn$ the category of all categories over
  $\wmodr{(N,N^{op})}$ having a separating object, where the
  morphisms between objects $(\bfC,\bfU_{\bfC})$ and
  $(\bfD,\bfU_{\bfD})$ are all functors $\bfF \colon \bfC
  \to \bfD$ satisfying $\bfU_{\bfD} \bfF = \bfU_{\bfC}$.
\end{definition}
\begin{example}
  For each $A \in \bfwsn$, denote by $\wmodr{A}$ the
  category of all normal, unital representations $\pi \colon
  A \to {\cal L}(H)$ for which $\pi \circ \iota_{A}$ and
  $\pi\circ \iota_{A}^{op}$ are faithful, and all
  intertwiners.  This is a category over $\wmodrn$, where
  $\bfU_{A}\colon \wmodr{A} \to \wmodrn$ is given by
  $(L,\pi) \mapsto (L, \pi \circ \iota_{A},\pi \circ
  \iota^{op}_{A})$ on objects and $T \mapsto T$ on
  morphisms.  The only non-trivial thing to check is that
 $\wmodr{A}$ has a separating object; by
  \cite[Lemma 2.10]{connes:2} or \cite[IX Theorem 1.2
  iv)]{takesaki:2}, one can take the GNS-representation for
  a n.s.f.\ weight on $A$.

    For each morphism $\phi \colon A \to B$ in $\bfwsn$, we
    obtain a functor $\phi^{*} \colon \wmodr{B} \to
    \wmodr{A}$, given by $(L,\pi) \mapsto (L,\pi \circ
    \phi)$ on objects and $T \mapsto T$ on morphisms.
\end{example}
\begin{remark}
  In the definition above, $\bfCatn(\bfC,\bfD)$ need not be
  a set, and this may cause problems. There are several
  possible solutions: we can  fix a ``universe'' to
  work in, or replace the category $\wmodrn$ by a small
  subcategory and require categories over $\wmodrn$ to be
  small, too. It is clear how to modify the preceding example
  in that case.
\end{remark}
\begin{proposition} \label{proposition:fiber-cat-mod} There
  exists a contravariant functor $\bfMod \colon \bfwsn \to
  \bfCatn$ given by $A \mapsto \bfMod (A):=
  (\wmodr{A},\bfU_{A})$ on objects and $\phi \mapsto
  \bfMod(\phi):= \phi^{*}$ on morphisms. \qed
\end{proposition}
For each category $\bfC \in \bfCatn$, choose a separating
object $G_{\bfC}$. Fix $\bfC\in \bfCatn$, let
$\bfU=\bfU_{\bfC}$, $\bfH=\bfH_{\bfC}$ $G=G_{\bfC}$,
$(H,\rho,\sigma) = \bfU G$, and define
$\bfEnd(\bfC):=\bfH(\bfC(G,G))' \subseteq {\cal L}(H)$. Then
$\rho(N)+\sigma(N) \subseteq \bfEnd(\bfC)$ because
$\bfH(\bfC(G,G)) \subseteq (\rho(N)+ \sigma(N))'$, and we
can consider $\bfEnd(\bfC)$ as an element of $\bfwsn$ with
respect to $\rho$ and $\sigma$.
\begin{lemma} \label{lemma:fiber-cat-end} 
  There exists a morphism $\eta_{\bfC} \colon \bfC \to
    \bfMod ( \bfEnd(\bfC))$ in $\bfCatn$, given by $X \mapsto
    (\bfU X, \rho^{X})$ on objects and $T \mapsto \bfH T$ on
    morphisms, where $\rho^{X} = \rho_{\bfH \bfC(G,X)}$ for
    each $X \in \bfC$. In particular, $\rho^{X}(\bfEnd(\bfC))
    \subseteq \bfH(\bfC(X,X))'$ for each $X\in \bfC$.
\end{lemma}
\begin{proof}
 Let $X \in \bfC$ and $(K,\phi,\psi)=\bfU X$. Lemma
  \ref{lemma:rtp-morphism}, applied to $I:=\bfH\bfC(G,X)
  \subseteq {\cal L}(\bfH G,\bfH X)$, gives a normal
  representation $\rho_{I} \colon (I^{*}I)' \to {\cal
    L}(K)$. Since $I^{*}I \subseteq \bfH \bfC(G,X)$ by
  assumption on $\bfC$, we have $\bfEnd(\bfC) \subseteq
  (I^{*}I)'$ and can define
  $\rho^{X}=\rho_{I}|_{\bfEnd(\bfC)}$.  Each element of $I$
  intertwines $\rho$ with $\phi$ and $\sigma$ with $\psi$,
  whence $\bfU X = (K,\rho_{I} \circ \rho,\rho_{I} \circ
  \sigma) =\bfU_{\bfEnd(\bfC)}(\eta_{\bfC}X)$.

  Next, let $Y \in \bfC$, $T\in \bfC(X,Y)$,
  $J:=\bfH\bfC(G,Y)$. Then
  $\bfH(T)\rho_{I}(S)=\rho_{J}(S)\bfH(T)$ for all $S \in
  \bfEnd(G)$ because $\bfH(T) I \in J$, and therefore
  $\bfH(T)$ is a morphism from $(\bfH X,\rho^{X})$ to $(\bfH
  Y,\rho^{Y})$. By definition,
  $\bfH_{\bfEnd(\bfC)}(\eta_{\bfC}(T))= \bfH T$.
\end{proof}
\begin{remark} \label{remark:fiber-cat-rho}
    If $G' \in \bfC$ is another separating object, then
    $\rho^{G'} \colon \bfH(\bfC(G,G))' \to
    \bfH(\bfC(G',G'))'$ is an isomorphism with inverse
    $\rho_{\bfH\bfC(G',G)}$.
\end{remark}

We eventually show that the assignment $\bfC \to
\bfEnd(\bfC)$ extends to a functor $\bfEnd \colon \bfCatn
\to \bfwsn$ that is adjoint to $\bfMod$. The key is a
more careful analysis of functors from a category $\bfC \in
\bfCatn$ to categories of the form $\bfMod(A)$, where $A \in
\bfwsn$.  Such functors themselves can be considered as
objects of a category as follows.

For all $\bfC,\bfD \in \bfCatn$, the elements of
$\bfCatn(\bfC,\bfD)$ are the objects of a category, where
the morphisms are all natural transformations with the usual
composition.

Similarly, for all $A,B \in \bfCatn$, the morphisms in
$\bfwsn(A,B)$ can be considered as objects of a category,
where the morphisms between $\phi,\psi$ are all $b \in B$
satisfying $b\phi(a)=\psi(a)b$ for all $a \in A$, and where
composition is given by multiplication.

\begin{proposition} \label{proposition:fiber-cat-adjoint}
  Let $A \in \bfwsn$ and $\bfC \in \bfCatn$. Then there exists an
  isomorphism $\Phi_{\bfC,A} \colon \bfCatn(\bfC,\bfMod(A))
  \to \bfwsn(A,\bfEnd(\bfC))$ with inverse
  $\Psi_{\bfC,A}:=\Phi_{\bfC,A}^{-1}$ such that
  \begin{enumerate}
  \item $\Phi_{\bfC,A}(\bfF)$ is defined by $\bfF G_{\bfC} =
    (\bfH_{\bfC} G_{\bfC}, \Phi_{\bfC,A}(\bfF))$ for each functor
    $\bfF\colon \bfC \to \bfMod(A)$ and
    $\Phi_{\bfC,A}(\alpha)=\alpha_{G_{\bfC}}$ for each natural
    transformation $\alpha$ in $\bfCatn(\bfC,\bfMod(A))$,
  \item $\Psi_{\bfC,A}(\pi)=\bfMod(\pi)\circ \eta_{\bfC} \colon \bfC
    \to \bfMod(\bfEnd(\bfC)) \to \bfMod(A)$ for each object
    $\pi$ and $\Psi_{\bfC,A}(S) = (\rho^{X}(S))_{X \in \bfC} $ for
    each morphism $S$ in $\bfwsn(A,\bfEnd(\bfC))$.
  \end{enumerate}
\end{proposition}
Explicitly, $\Psi_{\bfC,A}(\pi)$ is given by $X \mapsto
(\bfH_{\bfC}X, \rho^{X} \circ \pi)$ on objects and $T
\mapsto \bfH_{\bfC} T$ on morphisms.

The proof of Proposition \ref{proposition:fiber-cat-adjoint}
involves the following result.
\begin{lemma} \label{lemma:fiber-cat-transformation} Write
  $\bfU_{\bfC}G_{\bfC}=(\bfH_{\bfC} G_{\bfC},\rho,\sigma)$.
  Then the assignments $\alpha \mapsto \alpha_{G_{\bfC}}$
  and $(\rho^{X}(S))_{X \in \bfC} \mapsfrom S$ are inverse
  bijections between all natural transformations $\alpha$ of
  $\bfH_{\bfC}$ (or $\eta_{\bfC}$) and all elements $S \in
  \bfEnd(G_{\bfC})$ (or $S \in \bfEnd(G_{\bfC}) \cap
  (\rho(N)+\sigma(N))'$, respectively).
\end{lemma}
\begin{proof}
  A family of morphisms $(\alpha_{X}\colon
  \bfH_{\bfC}X\to\bfH_{\bfC}X)_{X \in \bfC}$ is a natural
  transformation of $\bfH_{\bfC}$ if and only if
  $\alpha_{X}T=T\alpha_{X}$ for all $X \in \bfC$ and $T \in
  \bfH_{\bfC}(G_{\bfC},X)$, that is, if
  $\alpha_{X}=\rho^{X}(\alpha_{G_{\bfC}})$ and
  $\alpha_{G_{\bfC}} \in \bfEnd(\bfC)$.  Such a family is a
  natural transformation of $\eta_{\bfC}$ if and only if
  additionally, $\alpha_{X}=\rho^{X}(\alpha_{G_{\bfC}})$ is
  a morphism of $\bfU_{\bfC}X$ for each $X \in \bfC$ or,
  equivalently, if $\alpha_{G_{\bfC}} \in
  (\rho(N)+\sigma(N))'$.
\end{proof}
\begin{proof}[Proof of Proposition
  \ref{proposition:fiber-cat-adjoint}]
  Lemma \ref{lemma:fiber-cat-transformation} implies that
  $\Psi:=\Psi_{\bfC,A}$ is well defined by ii).
  Let us show that $\Phi:=\Phi_{\bfC,A}$ is well defined by
  i). For each $\bfF$ as above, the image
  $\bfH_{\bfMod(A)}(\bfF(\bfC(G_{\bfC},G_{\bfC}))) =
  \bfH_{\bfC}(\bfC(G_{\bfC},G_{\bfC}))$ consists of
  intertwiners for $\Phi(\bfF)$ and hence $(\Phi(\bfF))(A)
  \subseteq \bfH_{\bfC}(\bfC(G_{\bfC},G_{\bfC}))' =
  \bfEnd(\bfC)$.  Likewise, for each $\alpha$ as above,
  $\alpha_{G_{\bfC}}$ intertwines
  $\bfH_{\bfC}(\bfC(G_{\bfC},G_{\bfC}))$ and hence
  $\alpha_{G_{\bfC}} \in \bfEnd(\bfC)$. Finally, $\Phi(\alpha
  \circ \beta) = \alpha_{G_{\bfC}} \circ \beta_{G_{\bfC}} =
  \Phi(\alpha)\Phi(\beta)$ for all composable
  $\alpha,\beta$.

  Next, $\Phi \circ \Psi=\Id$ because for each $\pi$ as
  above, $\Psi(\pi)(G_{\bfC}) = (\bfH_{\bfC} G_{\bfC},
  \rho^{G_{\bfC}} \circ \pi)$ so that $\Phi(\Psi(\pi)) =
  \rho^{G_{\bfC}} \circ \pi = \pi$, and for each $S$ as
  above, the component of $(\rho^{X}(S))_{X \in \bfC}$ at
  $X=G_{\bfC}$ is $\rho^{G_{\bfC}}(S)=S$.

  Finally, we prove $\Psi \circ \Phi=\Id$.
  Let $\bfF$ be as above and define $\phi^{X}$ by $\bfF X
  =(\bfH_{C} X, \phi^{X})$ for each $X \in \bfC$. Then
  $\Phi(\bfF)= \phi^{G_{\bfC}}$, and for each $a \in A$, the
  family $(\phi^{X}(a))_{X \in \bfC}$ is a natural
  transformation of
  $\bfH_{\bfMod(A)}\circ\bfF=\bfH_{\bfC}$ and
  coincides by Lemma \ref{lemma:fiber-cat-transformation}
  with $(\rho^{X}(\phi^{G_{\bfC}}(a)))_{X \in
    \bfC}$. Therefore, $\bfF X = (\bfH_{\bfC} X,\phi^{X}) =
  (\bfH_{\bfC}X,\rho^{X} \circ \Phi(\bfF)) =
  \Psi(\Phi(\bfF))(X)$ for each $X \in \bfC$.  On morphisms,
  $\Psi(\Phi(\bfF))$ and $\bfF$ coincide anyway.  For each
  $\alpha$ as above, $\Psi(\Phi(\alpha)) =
  (\rho^{X}(\alpha_{G_{\bfC}}))_{X \in \bfC} = \alpha$ by
  Lemma \ref{lemma:fiber-cat-transformation}.
\end{proof}

\begin{corollary} \label{corollary:fiber-cat}
  \begin{enumerate}
  \item Let $A \in \bfwsn$ and consider $\Id_{A}$ as an
    object of $\bfC:=\bfMod(A)$. Then
    $\Phi_{\bfC,A}(\Id_{\bfC}) \colon A \to
    \bfEnd(\bfMod(A))$ is an isomorphism in $\bfwsn$ with
    inverse $\epsilon_{A}:=\rho^{\Id_{A}}$.
  \item Let $A,B \in \bfwsn$. The the isomorphism
    $\bfMod_{(A,B)} := \Psi_{\bfMod(B),A} \circ
    (\epsilon_{B}^{-1})_{*}  \colon \bfws_{(N,N^{op})}(A,B)
    \to \bfws_{(N,N^{op})}(A,\bfEnd(\bfMod(B))) 
    \to \bfCat_{(N,N^{op})}(\bfMod(B),\bfMod(A))$ is given by
    $\phi \mapsto \bfMod(\phi)$ on objects and $b \mapsto
    (\pi(b))_{(L,\pi)}$ on morphisms.
  \item Let $\bfC,\bfD \in \bfCatn$.  Then the functor
    $\bfEnd_{(\bfC,\bfD)}:= \Phi_{\bfC,\bfEnd(\bfD)} \circ
    (\eta_{\bfD})_{*} \colon \bfCatn(\bfC,\bfD) \to
    \bfCatn(\bfC,\bfMod(\bfEnd(\bfD))) \to
    \bfwsn(\bfEnd(\bfD),\bfEnd(\bfC))$ is given by $\bfF
    \mapsto \rho^{\bfF G_{\bfC}}$ on objects and $\alpha
    \mapsto \bfH_{\bfD}(\alpha_{G_{\bfC}})$ on
    morphisms. \qed
  \end{enumerate}
\end{corollary}
\begin{proof}
Assertions   i) and iii) follow  immediately from the
definitions and Proposition \ref{proposition:fiber-cat-adjoint}. 
Let us prove ii). For each object $\phi$, we have $G_{\bfMod(B)} =
  (\bfH_{\bfMod(B)},\epsilon_{B}^{-1})$ and
  $\Phi_{\bfMod(B),A}(\bfMod(\phi)) = \epsilon_{B}^{-1}
  \circ \phi$, whence $\Psi_{\bfMod(B),A}(\epsilon_{B}^{-1}
  \circ \phi) = \bfMod(\phi)$, and for each morphism $b$,
  the family $\alpha:=(\pi(b))_{(L,\pi)}$ is a natural
  transformation and $\Phi_{\bfMod(B),A}(\alpha) =
  \alpha_{G_{\bfMod(B)}} = \epsilon_{B}^{-1}(b)$.
\end{proof}
The relative tensor product on $\wmodrn$ induces a product
on $\bfCatn$ as follows. Let $\bfC,\bfD \in \bfCatn$. Then
$\bfC \times \bfD$ and the functor
\begin{align*}
  \bfU_{\bfC \times \bfD}=(- \rtensor{}{\mu}{} -) \circ
  (\bfU_{\bfC} \times \bfU_{\bfD}) \colon \bfC \times \bfD
  \to \wmodrn,
\end{align*}
form a category over $\wmodrn$ with
separating object $(G_{\bfC},G_{\bfD})$.  Thus, we
obtain a monoidal structure on $\bfCatn$, given by
$(\bfC,\bfD) \mapsto \bfC \times \bfD$ on objects and
$(\bfF,\bfG) \mapsto \bfF \times \bfG$ on morphisms.
\begin{corollary}
  For all $A,B,C \in \bfwsn$, there exists an isomorphism
  \begin{align*}
    \Xi\colon \bfws_{(N,N^{op})}(A, B \fibre{}{\mu}{} C) \to
    \bfCat_{(N,N^{op})}(\bfMod(B)\times\bfMod(C),\bfMod(A))
  \end{align*}
  such that for each object $\pi$, the functor $\Xi(\pi)$ is
  given by $((L,\tau),(M,\upsilon)) \mapsto (L
  \rtensor{}{\mu}{} M, (\tau \fibre{}{\mu}{} \upsilon) \circ
  \pi)$ and $(S,T) \mapsto S \rtensor{}{\mu}{} T$, and for
  each morphism $x \colon \pi_{1} \to \pi_{2}$, the
  transformation $\Xi(b) \colon \Xi(\pi_{1})
  \to\Xi(\pi_{2})$ is given by
  $\Xi(b)_{((L,\tau),(M,\upsilon))} = (\tau \fibre{}{\mu}{}
  \upsilon)(x)$.
\end{corollary}
\begin{proof}
  Let $\bfB:=\bfMod(B)$, $\bfC:=\bfMod(C)$,
  $\bfD:=\bfB\times \bfC$. Then
  $G:=(G_{\bfB},G_{\bfC})$ is 
  separating  and
  \begin{align*}
    \rho^{G} \colon \bfEnd(\bfD) \to \bfH_{\bfD}(\bfD(G,G))'
    = (\bfEnd(\bfB)' \rtensor{}{\mu}{} \bfEnd(\bfC)')' = 
    \bfEnd(\bfB) \fibre{}{\mu}{} \bfEnd(\bfC) \cong B
    \fibre{}{\mu}{} C
  \end{align*}
  is an isomorphism by Remark \ref{remark:fiber-cat-rho}.
  Moreover, if $X=(L,\tau) \in \bfB$, $Y=(M,\upsilon)\in
  \bfC$, then $\rho^{(X,Y)} = (\tau \fibre{}{\mu}{}\upsilon)
  \circ \rho^{G}$ by Lemma \ref{lemma:rtp-morphism} 
  because $\tau \fibre{}{\mu}{}\upsilon = \rho_{J}$, where
  $J=\bfH_{\bfB}(\bfB(G_{\bfB},X)) \rtensor{}{\mu}{}
  \bfH_{\bfC}(\bfC(G_{\bfC},Y))$, and $J \cdot
  \bfH_{\bfD}(\bfD(G_{\bfD},G)) \subseteq
  \bfH_{\bfD}(\bfD(G_{\bfD},(X,Y)))$.  Now, the assertion
  follows from Proposition
  \ref{proposition:fiber-cat-adjoint}.
\end{proof}

The categories $\bfwsn$ and $\bfCatn$ are enriched over the
monoidal category $\bfCat$ of small categories
\cite{kelly:enriched}, or, equivalently, are 2-categories,
meaning that the morphisms between fixed objects are
themselves objects of a small category, as explained before
Proposition \ref{proposition:fiber-cat-adjoint}, and that
the composition of morphisms between fixed objects extends
to a functor, where
\begin{align} \label{eq:fiber-cat-vertical-vn}
  \xymatrix{ B
 \ar@/^/[r]^{\psi_{1}} \ar@/_/[r]_{\psi_{2}}
    \ar@{}[r]|{\Downarrow\, c} & C} \circ
  \xymatrix{ A \ar@/^/[r]^{\phi_{1}} \ar@/_/[r]_{\phi_{2}}
    \ar@{}[r]|{\Downarrow\, b} & B }
 &= 
  \xymatrix@C=100pt{ A \ar@/^0.75pc/[r]^{\psi_{1}\circ\phi_{1}}
    \ar@/_0.75pc/[r]_{\psi_{2} \circ \phi_{2}}
    \ar@{}[r]|{\Downarrow\, \psi_{2}(b)c} & C }
  \text{ in } \bfwsn, \\ \label{eq:fiber-cat-vertical-cat}
 \xymatrix{ \bfC
 \ar@/^/[r]^{\bfG_{1}} \ar@/_/[r]_{\bfG_{2}}
    \ar@{}[r]|{\Downarrow\, \beta} & \bfD} \circ
  \xymatrix{ \bfB \ar@/^/[r]^{\bfF_{1}} \ar@/_/[r]_{\bfF_{2}}
    \ar@{}[r]|{\Downarrow\, \alpha} & \bfC }
 &= 
  \xymatrix@C=100pt{ \bfB \ar@/^0.75pc/[r]^{\bfG_{1}\circ\bfF_{1}}
    \ar@/_0.75pc/[r]_{\bfG_{2} \circ \bfF_{2}}
    \ar@{}[r]|{\Downarrow\, \beta_{\bfF_{2}} \circ \bfG_{1}\alpha } & \bfD }
  \text{ in } \bfCatn.
\end{align}
Recall that a contravariant functor between enriched
categories $\bfC,\bfD$ consists of an assignment $\bfF
\colon \ob \bfC \to \ob \bfD$ and, for each pair of objects
$X,Y\in \bfC$, a functor $\bfF_{(X,Y)} \colon \bfC(X,Y) \to
\bfD(\bfF Y,\bfF X)$ that is compatible with composition in
a natural sense.  We now show that the assignments
$\bfMod,\bfEnd$ defined above are functors in this sense and
that the isomorphisms in Proposition
\ref{proposition:fiber-cat-adjoint} form part of an
adjunction between $\bfMod$ and $\bfEnd$.  For background on
enriched categories, see \cite{kelly:enriched}.

\begin{theorem} \label{theorem:fiber-cat-adjunction} The
  assignments $\bfMod$, $\bfEnd$ define contravariant
  functors $\bfMod \colon \bfwsn \to \bfCatn$ and $\bfEnd
  \colon \bfCatn \to \bfwsn$ of enriched categories, and the
  isomorphisms $(\Phi_{\bfC,A})_{\bfC,A}$ define an
  adjunction whose unit is $(\eta_{\bfC})_{\bfC \in
    \bfCatn}$ and counit is $(\epsilon_{A})_{A \in \bfwsn}$.
\end{theorem}

\begin{proof}
  We first show that $\bfMod$ and $\bfEnd$ are functors of
  enriched categories. By Corollary
  \ref{corollary:fiber-cat}, it suffices to prove this for
  $\bfEnd$.  Consider a diagram as in
  \eqref{eq:fiber-cat-vertical-cat} and let
  $a=\bfEnd_{(\bfB,\bfC)}(\alpha)$,
  $b=\bfEnd_{(\bfC,\bfD)}(\beta)$,
  $c=\bfEnd_{(\bfB,\bfD)}(\beta_{\bfF_{2}}\circ
  \bfG_{1}\alpha)$. We have to show that then the cells
  \begin{align*}
    \xymatrix@C=40pt{ \bfEnd(\bfC)
      \ar@/^0.75pc/[r]^{\bfEnd_{(\bfB,\bfC)}(\bfF_{1})}
      \ar@{}[r]|{\Downarrow a}
      \ar@/_0.75pc/[r]_{\bfEnd_{(\bfB,\bfC)}(\bfF_{2})} &
      \bfEnd(\bfB) } \circ \xymatrix@C=40pt{ \bfEnd(\bfD)
      \ar@/^0.75pc/[r]^{\bfEnd_{(\bfC,\bfD)}(\bfG_{1})}
      \ar@{}[r]|{\Downarrow b}
      \ar@/_0.75pc/[r]_{\bfEnd_{(\bfC,\bfD)}(\bfG_{2})} &
      \bfEnd(\bfC) } \ \text{ and } \ \xymatrix@C=50pt{
      \bfEnd(D)
      \ar@/^0.75pc/[r]^{\bfEnd_{(\bfB,\bfD)}(\bfG_{1}\bfF_{1})}
      \ar@/_0.75pc/[r]_{\bfEnd_{(\bfB,\bfD)}(\bfG_{2}\bfF_{2})}
      \ar@{}[r]|{\Downarrow c} &\bfEnd(B) }
  \end{align*}
  are equal.  By definition,
  $a=\bfH_{\bfC}(\alpha_{G_{\bfB}})$,
  $b=\bfH_{\bfD}(\beta_{G_{\bfC}})$, and by Lemma
  \ref{lemma:fiber-cat-transformation},
  \begin{align*}
    c = \bfH_{\bfD}(\beta_{\bfF_{2} G_{\bfB}} \cdot
    \bfG_{1}(\alpha_{G_{\bfB}})) =
    \rho^{\bfF_{2}G_{\bfB}}(\bfH_{\bfD}(\beta_{G_{\bfC}})) \cdot
    \bfH_{\bfC}(\alpha_{G_{\bfB}}) = \bfEnd(\bfF_{2})(b)
    \cdot a.
  \end{align*}

  It remains to show that for all morphisms $\phi \colon A
  \to B$ in $\bfwsn$ and $\bfF \colon \bfC\to \bfD$ in
  $\bfCatn$, the diagram
  \begin{align*}
    \xymatrix@R=10pt{
      \bfCatn(\bfD,\bfMod(B)) \ar[r]^{\Phi_{\bfD,B}}  \ar[d]
      &
      \bfwsn(B,\bfEnd(\bfD)) \ar[d] \\
      \bfCatn(\bfC,\bfMod(A)) \ar[r]^{\Phi_{\bfC,A}} 
      &
      \bfwsn(A,\bfEnd(\bfC))
    }
  \end{align*}
  commutes, where the vertical maps are induced by $\bfF$
  and $\bfMod_{(A,B)}(\phi)$ on the left and $\phi$ and
  $\bfEnd_{(\bfC,\bfD)}(\bfF)$ on the right, respectively,
  or, more precisely, that for each object $\bfG$ and each
  morphism $\alpha$ in $\bfCatn(\bfD,\bfMod(B)) $,
  \begin{align*}
    \bfEnd_{(\bfC,\bfD)}(\bfF) \circ \Phi_{\bfD,B}(\bfG) \circ
    \phi &= \Phi_{\bfC,A}(\bfMod_{(A,B)}(\phi) \circ \bfG\circ
    \bfF), & \bfEnd_{(\bfC,\bfD)}(\bfF)(\alpha) &=
    \bfMod_{(A,B)}(\phi)(\alpha_{\bfF}).
  \end{align*}
  The second equation holds because of Lemma
  \ref{lemma:fiber-cat-transformation} and the relation
  \begin{align*}
   \bfEnd_{(\bfC,\bfD)}(\bfF)(\alpha_{G_{\bfC}}) =
  \rho^{\bfF G_{\bfC}}(\alpha_{G_{\bfD}}) = \alpha_{\bfF
    G_{\bfC}} = \bfMod_{(A,B)}(\phi)(\alpha_{\bfF
    G_{\bfC}}) 
  \end{align*}
  first one holds because by Corollary
  \ref{corollary:fiber-cat},
  \begin{align*}
    \bfEnd_{(\bfC,\bfD)}(\bfF) \circ \Phi_{\bfD,B}(\bfG)
    \circ \phi &= \rho^{\bfF G_{\bfC}} \circ
    \Phi_{\bfD,B}(\bfG)
    \circ \phi, \\
    (\bfMod_{(A,B)}(\phi) \circ \bfG \circ \bfF)(G_{\bfC})
    &= (\bfH_{\bfC}G_{\bfC},\rho^{\bfF G_{\bfC}} \circ
    \Phi_{\bfD,B}(\bfG) \circ \phi). \qedhere
  \end{align*}
\end{proof}

\section{The special case of a commutative base}

\label{section:commutative}

Let $Z$ be a locally compact Hausdorff
space with a Radon measure $\mu$ of full support,
and identify $C_{0}(Z)$ with multiplication operators on
${\cal L}(L^{2}(Z,\mu))$. Then the relative tensor product
and the fiber product over the $C^{*}$-base
$\frakb=(L^{2}(Z,\mu),C_{0}(Z),C_{0}(Z))$ can be related to
the fiberwise product of bundles as follows.

\paragraph{Modules and their relative tensor product}

Denote by $\bfMod_{\frakb}$, $\bfMod_{C_{0}(Z)}$,
$\bfBund_{Z}$ the categories of all $C^{*}$-$\frakb$-modules
with all morphisms, of all Hilbert $C^{*}$-modules over
$C_{0}(Z)$, and of all continuous Hilbert bundles over $Z$;
for the precise definition of the latter, see \cite{dupre}.
Each of these categories carries a monoidal structure, where
the product
\begin{itemize}
\item of $E,F \in \bfMod_{C_{0}(Z)}$ is the separated
  completion of $E \odot F$ with respect to the inner
  product $\langle \xi \odot \eta|\xi'\odot \eta'\rangle =
  \langle \xi|\xi'\rangle\langle\eta|\eta'\rangle$, denoted
  by $E \rtensor{}{C_{0}(Z)}{} F$,
\item of ${\cal E},{\cal F} \in \bfBund_{Z}$ is the
  fibrewise tensor product of ${\cal E}$ and ${\cal F}$,
\item of $H_{\beta},K_{\gamma}\in \bfMod_{\frakb}$ is $(H
  \rtensor{\beta}{\frakb}{\gamma} K, \beta \bowtie \gamma)$,
  where $\beta \bowtie \gamma := [\kgamma{2}\beta] =
  [\kbeta{1}\gamma]$; here, note that
  ${_{\beta}H_{\beta}},{_{\gamma}K_{\gamma}}$ are
  $C^{*}$-$(\frakb,\frakb)$-modules.
\end{itemize}
There exist equivalences of monoidal categories
$\bfMod_{C_{0}(Z)}
\underset{\Gamma_{0}}{\stackrel{\bfB}{\rightleftarrows}}
\bfBund_{Z}$ and $\bfMod_{C_{0}(Z)}
\underset{\bfU}{\stackrel{\bfF}{\rightleftarrows}}
\bfMod_{\frakb}$ such that for each $E \in
\bfMod_{C_{0}(Z)}$, ${\cal F} \in \bfBund_{Z}$, $H_{\beta}
\in \bfMod_{\frakb}$,
\begin{itemize}
\item $\bfB E=\bigsqcup_{z \in Z} E_{z}$ is and
  $\Gamma_{0}(\bfB E) =\{ (\xi_{z})_{z} \mid \xi \in E\}$,
  where $E_{z}$ is the completion of $E$ with respect to the
  inner product $(\xi,\eta) \mapsto \langle \xi|\eta\rangle
  (z)$, and $\xi \mapsto \xi_{z}$ denotes the quotient map
  $E \to E_{z}$,
\item the operations on the space of sections
  $\Gamma_{0}({\cal F}) \in \bfMod_{C_{0}(Z)}$ are defined
  fiberwise,
\item $\bfF E = (E \otimes_{C_{0}(Z)} L^{2}(Z,\mu), l(E))$,
  where $l(\xi) \eta = \xi \otimes_{C_{0}(Z)} \eta$ for each
  $\xi \in E, \eta \in L^{2}(Z,\mu)$,
\item $\bfU H_{\beta} = \beta \in \bfMod_{C_{0}(Z)}$.
\end{itemize}
 The first equivalence is explained in \cite{dupre},
and the second one is easily verified. Compare also Examples
\ref{example:rtp-module-commutative} and
\ref{examples:rtp-rtp} ii).

\paragraph{Algebras} 
Denote by $\czalg$ the category of all continuous
$C_{0}(Z)$-algebras with full support \cite{}, where the
morphisms between $A,B \in \czalg$ are all $C_{0}(Z)$-linear
nondegenerate $*$-homomorphisms $\pi \colon A \to M(B)$, and
by $\cbalg$ the category of all $C^{*}$-$\frakb$-algebras
$A^{\beta}_{H}$ satisfying $[\rho_{\beta}(C_{0}(Z))A]=A$ and
$[A\beta]= \beta$, where the morphisms between
$A^{\beta}_{H}$, $B^{\gamma}_{K} \in \cbalg$ are all $\pi
\in \bfcs_{\frakb
}(A^{\beta}_{H},M(B)^{\gamma}_{K})$
satisfying $[\pi(A)B]=B$.  Then there exists a functor
$\cbalg \to \czalg$, given by $A^{\beta}_{H} \mapsto
(A,\rho_{\alpha})$ and $\pi \mapsto \pi$, and this functor has a
full and faithful left adjoint which embeds $\czalg$ into $
\cbalg$ \cite[Theorem 6.6]{timmermann:coactions}.

\paragraph{The fiber product of commutative
  $C^{*}$-$\frakb$-algebras}
We finally discuss the fiber product of commutative
$C^{*}$-$\frakb$-algebras and start with preliminaries. Let
$Z$ be a locally compact space, $E$ a Hilbert $C^{*}$-module
over $C_{0}(Z)$, and $\bfB E= \bigsqcup_{z \in Z} E_{z}$ the
corresponding Hilbert bundle. The topology on $\bfB E$ is
generated by all open sets of the form
$U_{V,\eta,\epsilon}=\{ \zeta | z \in V, \zeta \in E_{z},
\|\eta_{z} - \zeta\|_{E_{z}}<\epsilon\}$, where $V \subseteq
Z$ is open, $\eta \in E$, $\epsilon>0$.  Denote by $q \colon
\bigsqcup_{z \in Z} {\cal L}(E_{z}) \to Z$ the natural
projection and define for each $\eta,\eta' \in E$ maps
\begin{align*}
 \omega_{\eta,\eta'}& \colon \bigsqcup_{z \in Z}
{\cal L}(E_{z}) \to \complex, \ T \mapsto \langle
\eta_{q(T)}|T\eta'_{q(T)}\rangle, & \upsilon^{(*)}_{\eta}
\colon
\colon \bigsqcup_{z \in Z} {\cal L}(E_{z}) \to \bigsqcup_{z
  \in Z} E_{z}, \ T \mapsto T^{(*)}\eta_{q(T)}.
\end{align*}
The {\em weak topology (strong-*-topology)} on $\bigsqcup_{z
  \in Z} {\cal L}(E_{z})$ is the weakest one that makes $q$
and all maps of the form $\omega_{\eta,\eta'}$ (of the form
$\upsilon^{(*)}_{\eta}$) continuous.  

Let $A$ be a commutative $C^{*}$-algebra, $\pi \colon
C_{0}(Z) \to M(A)$ a $*$-homomorphism, and $\chi \in
\hA$. Then we identify $E \otimes_{\phi^{*}} A
\otimes_{\chi} \complex$ with $E_{z}$, where $z\in Z$
corresponds to $\chi\circ \pi \in \widehat{C_{0}(Z)}$, via
$\eta \otimes_{\pi}  a\otimes_{\chi} \lambda \mapsto \lambda
\chi(a)\eta_{z}$. A map $T \colon \hA \to \bigsqcup_{z \in Z}
{\cal L}(E_{z})$ is {\em weakly
  vanishing (strong-$*$-vanishing) at infinity} if for all
$\eta,\eta' \in E$, the map $\omega_{\eta,\eta'} \circ T$ (the
maps $\chi \mapsto \|\upsilon^{(*)}_{\eta}(T(\chi))\|$) vanish
at infinity.
\begin{lemma} \label{lemma:fp-commutative} Let
  $A^{\beta}_{H}$ be a $C^{*}$-$\frakb$-algebra,
  $K_{\gamma}$ a $C^{*}$-$\frakbo$-module, $x \in {\cal
    L}(\HfibreK)$. Assume that $A$ is commutative,
  $[\rho_{\beta}(C_{0}(Z))A]=A$, and $\langle
  \gamma|_{2}x|\gamma\rangle_{2} \subseteq A$. Define $F_{x}
  \colon \hA \to \bigsqcup_{z \in Z} {\cal L}(\gamma_{z})$
  by $\chi \mapsto (\chi \ast \Id)(x)$. Then:
  \begin{enumerate}
  \item $F_{x}$ is weakly  continuous, weakly vanishing at infinity.
  \item $x \in \Ind_{|\gamma\rangle_{2}}(A)$ if and only if
    $F_{x}$ is strong-$*$ continuous, strong-$*$-vanishing at infinity.
  \end{enumerate}
\end{lemma}
\begin{proof}
  First, note that for all $\eta,\eta' \in \gamma$ and $\chi
  \in \hat A$,
  \begin{align*} \chi(\langle \eta|_{2}x|\eta'\rangle_{2}) =
    \langle 1 {_{(\chi \circ \rho_{\beta})}\tl} \eta | (\chi
    \ast \Id ) (x) (1 {_{(\chi \circ \rho_{\beta})}\tl}
    \eta')\rangle = \langle \eta_{(\chi \circ \rho_{\beta})} |
    F_{x}(\chi) \eta'_{(\chi \circ \rho_{\beta})}\rangle.
  \end{align*}
  i) For each $\eta',\eta \in \gamma$, the map $\chi \mapsto
  \langle \eta_{(\chi \circ \rho_{\beta})} | F_{x}(\chi)
  \eta'_{(\chi \circ \rho_{\beta})}\rangle$  equals
  $\langle\eta|_{2}x|\eta'\rangle_{2} \in A$.
  
  ii) Assume that $F_{x}$ is strong-$*$ continuous vanishing
  at infinity and let $\eta \in \gamma$. Then the map $\chi
  \mapsto F_{x}(\chi)\eta_{(\chi \circ \rho_{\beta})}$ lies
  in $\Gamma_{0}(\gamma \tr_{\rho_{\beta}} A)$. Hence, there
  exists an $\omega \in \gamma \tr_{\rho_{\beta}} A$ such
  that $F_{x}(\chi)\eta_{(\chi \circ
    \rho_{\beta})}=\omega_{\chi}$ for all $\chi \in \hA$. We
  identify $\gamma \tr_{\rho_{\beta}} A$ with
  $[|\gamma\rangle_{2}A] \subseteq {\cal L}(H,\HfibreK)$ in
  the canonical manner and find that $x|\eta\rangle_{2} =
  \omega$ because $\chi(\langle \eta'|_{2} x|\eta\rangle_{2})
  = \langle \eta'_{(\chi \circ \rho_{\beta})} | \omega_{(\chi
    \circ \rho_{\beta})} \rangle =
  \chi(\langle\eta'|_{2}\omega)$ for all $\chi \in \hat A$,
  $\eta' \in \gamma$. Since $\eta \in \gamma$ was arbitrary,
  we can conclude $x |\gamma\rangle_{2} \subseteq
  [|\gamma\rangle_{2}A]$.  A similar argument, applied to
  $x^{*}$ instead of $x$, shows that $x^{*}
  |\gamma\rangle_{2} \subseteq [|\gamma\rangle_{2}A]$, and
  therefore $x \in \Ind_{\kgamma{2}}(A)$.  Reversing the
  arguments, we obtain the reverse implication.
\end{proof}

Let $X$ be a locally compact Hausdorff space with a
continuous surjection $p \colon X \to Z$ and a family of
Radon measures $\phi=(\phi_{z})_{z \in Z}$  such that
(i) $\supp \phi_{z}=X_{z}:=p^{-1}(z)$ for each $z \in Z$ and
(ii) the
map $\phi_{*}(f)\colon z \mapsto \int_{X_{z}} f d\phi_{z}$
is continuous for each $f\in C_{c}(X)$.  Define a Radon
measure $\nu_{X}$ on $X$ such that $\int_{X} f \rd\nu_{X} =
\int_{Z}  \phi_{*}(f)d\mu$ for all $f\in
C_{c}(X)$.  Then there exists a map $j_{X} \colon
C_{c}(X) \to {\cal L}(L^{2}(Z,\mu),L^{2}(X,\nu_{X}))$ such
that $j_{X}(f)h= fp^{*}(h)$ and
$j_{X}(f)^{*}g = \phi_{*}(\overline{f}g)$ for
all $f,g \in C_{c}(X),h\in C_{c}(Z)$. 
Similarly, let  $Y$ be a locally compact Hausdorff space
with a continuous map $q \colon Y \to Z$ and a family of
measures $\psi=(\psi_{z})_{z\in Z}$ satisfying the same
conditions as $X,p,\phi$, and define a Radon measure
$\nu_{Y}$ on $Y$ and an embedding $j_{Y} \colon C_{c}(Y) \to
{\cal L}(L^{2}(Z,\mu),L^{2}(Y,\nu_{Y}))$  as above.
Let
\begin{align*}
  H&:=L^{2}(X,\nu_{X}), & \beta&:=[j_{X}(C_{c}(X))], &
  A&:=C_{0}(X) \subseteq {\cal L}(L^{2}(X,\nu_{X})={\cal
    L}(H), \\
  K&:=L^{2}(Y,\nu_{Y}), & \gamma&:=[j_{Y}(C_{c}(Y))], &
  B&:=C_{0}(Y) \subseteq {\cal L}(L^{2}(Y,\nu_{Y}))= {\cal L}(K).
\end{align*}
Then $H_{\beta}$, $K_{\gamma}$ are $C^{*}$-$\frakb$-modules
and $A_{H}^{\beta}$, $B^{\gamma}_{K}$ are
$C^{*}$-$\frakb$-algebras, as one can easily
check. Considering $\beta$ and $\gamma$ as Hilbert
$C^{*}$-modules over $C_{0}(Z)$, we can canonically identify
$\beta_{z}\cong L^{2}(X_{z},\phi_{z})$ and $\gamma_{z} \cong
L^{2}(Y_{z},\psi_{z})$. Finally, define a Radon measure
$\nu$ on $ X \ftimes{p}{Z}{q} Y$ such that for all $h \in
C_{c}(X \ftimes{p}{Z}{q} Y)$,
\begin{gather*}
  \int_{X \ftimes{p}{Z}{q} Y} h \rd\nu = \int_{Z}
  \int_{X_{z}} \int_{Y_{z}} h(x,y) \rd\psi_{z}(y)
  \rd\phi_{z}(x) \rd\mu(z).
\end{gather*}
\begin{proposition} \label{proposition:fp-commutative}
  \begin{enumerate}
  \item There exists a  unitary $U \colon \HfibreK \to
    L^{2}(X \ftimes{p}{Z}{q} Y,\nu)$ such that
    $(\Phi(j_{X}(f) \tr h \tl j_{Y}(g)))(x,y) =
    f(x)h(p(x))g(y)$
    for all $f \in C_{c}(X)$, $g \in C_{c}(Y)$, $h \in
    C_{c}(Z)$, $(x,y) \in X \ftimes{p}{Z}{q} Y$.
  \item $\Ad_{U}(\AfibreB)$ is the $C^{*}$-algebra of all $f
    \in L^{\infty}(X \ftimes{p}{Z}{q} Y,\nu)$ that have
    representatives $f_{X},f_{Y}$ such that the maps $ X \to
    \TotL(\gamma)$ and $Y \to \TotL(\beta)$ given by $x
    \mapsto f_{X}(x, \cdot) \in
    L^{\infty}(Y_{p(x)},\psi_{p(x)})$ and $y \mapsto
    f_{Y}(\cdot, y) \in L^{\infty}(X_{q(y)},\phi_{q(y)})$
    respectively, are strong-$*$ continuous vanishing at
    infinity.
  \end{enumerate}
\end{proposition}
\begin{proof}
  The proof of assertion i) is straightforward, and
  assertion ii) follows immediately from Proposition Lemma
  \ref{lemma:fiber-properties} viii) and Lemma
  \ref{lemma:fp-commutative} ii).
\end{proof}
\begin{examples} \label{examples:fp-commutative}
  \begin{enumerate}
  \item Let $X,Y$ be discrete, $Z=\{0\}$, and let
    $\phi_{0}$, $\psi_{0}$ be the counting measures on
    $X,Y$, respectively. Then
    \begin{align*}
      C_{0}(X) \fibre{\beta}{\frakb}{\gamma} C_{0}(Y) \cong
      \{ f \in C_{b}(X \times Y) \mid
      &f(x,\frei) \in C_{0}(Y) \text{ for all } x \in X, \\
      & f(\frei,y) \in C_{0}(X) \text{ for all } y \in Y\}.
    \end{align*}
    This follows from Proposition
    \ref{proposition:fp-commutative} and the fact that for
    each $f\in C_{b}(X \times Y)$, the maps $X \to{\cal
      L}(l^{2}(Y))$, $x \mapsto f(x, \frei)$, and $ Y
    \to{\cal L}(l^{2}(X))$, $y \mapsto f(\frei,y)$, are
    strong-$*$ continuous vanishing at infinity if and only
    if $f(\frei,y) \in C_{0}(X)$ and $f(x,\frei) \in
    C_{0}(Y)$ for each $y \in Y$ and $x \in X$.     
  \item Let $X=\naturals$, $Z=\{0\}$, and let $\phi_{0}$ be the
    counting measure. Then
    \begin{multline*}
      C_{0}(\naturals) \fibre{\beta}{\frakb}{\gamma}
      C_{0}(Y) \cong \{ f \in C_{b}(\naturals \times Y) \mid
      (f(x,\cdot))_{x} \text{ is a sequence in }
      C_{0}(Y) \\ \text{that converges strongly to
      } 0 \in {\cal L}(L^{2}(Y,\psi_{0}))\}
    \end{multline*}
    because for each $f \in L^{\infty}(\naturals \times Y)$,
    the map $Y \to {\cal L}(l^{2}(\naturals))$, $y \mapsto
    f(\frei,y)$, is strong-$*$ continuous vanishing at
    infinity if and only if $f(x,\frei) \in C_{0}(Y)$ for
    all $x \in \naturals$.
    
    \item Let $X=Y=[0,1]$, $Z=\{0\}$, and let
    $\phi_{0}=\psi_{0}$ be the Lebesgue measure. For each
    subset $I \subseteq [0,1]$, denote by $\chi_{I}$ its
    characteristic function.  Then the function $f \in
    L^{\infty}([0,1]\times[0,1])$ given by $f(x,y)=1$ if $y
    \leq x$ and $f(x,y)=0$ otherwise belongs to $C([0,1])
    \fibre{\beta}{\frakb}{\gamma} C([0,1])$ because the
    functions $[0,1] \to L^{\infty}([0,1]) \subseteq {\cal
      L}(L^{2}([0,1]))$ given by $x \mapsto f(x,\frei) =
    \chi_{[0,x]}$ and $y \mapsto f(\frei,y) = \chi_{[y,1]}$
    are strong-$*$ continuous.  In particular, we see that
    $C([0,1]) \fibre{\beta}{\frakb}{\gamma} C([0,1])
    \nsubseteq C([0,1]\times[0,1]) = C([0,1]) \otimes
    C([0,1])$.
  \end{enumerate}
\end{examples}

\def\cprime{$'$}

\end{document}